\documentclass{gtart_h}

\def\ifplaintex{\expandafter\ifx\csname documentclass\endcsname\relax}

\def\ifplaintex{\expandafter\ifx\csname documentclass\endcsname\relax}

\ifplaintex 
\else
\fi

\expandafter\ifx\csname epsfbox\endcsname\relax\input epsf\fi

\def\gt{{\mathsurround=0pt\it $\cal G\mskip-2mu$eometry \&\ 
$\cal T\!\!$opology}}        

\def\gtp{{\mathsurround=0pt\it $\cal G\mskip-2mu$eometry \&\ 
$\cal T\!\!$opology $\cal P\!$ublications}}  

\def\lognumber#1{\def\thelognumber{#1}}
\def\volumenumber#1{\def\thevolumenumber{#1}}
\def\papernumber#1{\def\thepapernumber{#1}}
\def\volumeyear#1{\def\thevolumeyear{#1}}

\def\pagenumbers#1#2{\def\startpage{#1}\def\finishpage{#2}}
\def\published#1{\def\publishdate{#1}}
\def\proposed#1{\def\theproposer{#1}}
\def\seconded#1{\def\theseconders{#1}}
\def\received#1{\def\receiveddate{#1}}

\def\accepted#1{\def\accepteddate{#1}}

\def\coverauthors#1{\def\thecoverauthors{#1}}
\def\asciiauthors#1{\def\theasciiauthors{#1}}
\def\asciiaddress#1{\def\theasciiaddress{#1}}
\def\asciiemail#1{\def\theasciiemail{#1}}

\def\shortauthors#1{\def\theshortauthors{#1}}

\let\\\par\let\thelognumber\relax
\let\thevolumenumber\relax\let\thepapernumber\relax
\let\thevolumeyear\relax\let\thesamplenumber\relax\let\startpage\relax
\let\finishpage\relax\let\publishdate\relax\let\receiveddate\relax
\let\reviseddate\relax\let\accepteddate\relax\let\theasciititle\relax
\let\theasciiauthors\relax\let\theasciiaddress\relax
\let\theasciiabstract\relax
\let\theasciiemail\relax\let\theshortauthors\relax\let\theshorttitle\relax
\let\thecoverauthors\relax

\long\def\maketitlep{   

\count0=\startpage

\gt\hfill      
\hbox to 77pt{\vbox to 0pt{\vglue -15pt\epsfbox{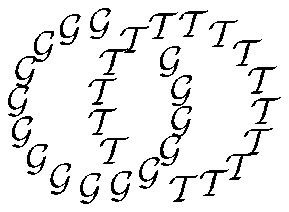}\vss}\hss}
\break
{\small\ifx\thesamplenumber\relax 
Volume \else Sample
\fi\thevolumenumber\ (\thevolumeyear)
\startpage--\finishpage\nl
Published: \publishdate}
\vglue 0.5truein plus 0.4fil minus 0.1truein

{\parskip=0pt\leftskip 0pt plus 1fil\def\\{\par\smallskip}{\ifplaintex\large
\else\Large\fi\bf\thetitle}\par\medskip}   

\vglue 0pt plus 0.1fil 

{\parskip=0pt\leftskip 0pt plus 1fil\def\\{\par}{\sc\theauthors}
\par\medskip}

\vglue 0pt plus 0.1fil 

{\small\parskip=0pt\let\newline\\
{\leftskip 0pt plus 1fil\def\\{\par}{\sl\theaddress}\par}
\expandafter\ifx\theemail\relax    
\relax\else\vglue 5pt plus 0.02fil minus 2pt\def\\{\stdspace{\rm 
and}\stdspace} 
\cl{Email:\stdspace\tt\theemail}\fi
\ifx\theurl\relax                  
\relax\else\vglue 5pt plus 0.02fil minus 2pt\def\\{\stdspace{\rm 
and}\stdspace}
\cl{URL:\stdspace\tt\theurl}\fi\par}

\vglue 7pt plus 0.3fil minus 3pt

{\bf Abstract}
\vglue 5pt plus 0.1fil minus 2pt

\theabstract

\vglue 7pt plus 0.3fil minus 3pt

{\bf AMS Classification numbers}\quad Primary:\quad \theprimaryclass

Secondary:\quad \thesecondaryclass

\vglue 5pt plus 0.3fil minus 2pt

{\bf Keywords:}\quad \thekeywords

\vglue 10pt plus 0.5fil minus 5pt

{\small  Proposed: \theproposer\hfill Received: \receiveddate\nl
Seconded: \theseconders\hfill 
\ifx\reviseddate\relax                         
Accepted: \accepteddate                        
\else
Revised: \reviseddate                          
\fi}
\eject
}       

\font\phead=cmsl9 scaled 950
\font=cmsl9 scaled 1050
\font\pnum=cmbx10 scaled 913
\font\lnum=cmbx10 
\font\pfoot=cmsl9 scaled 950
\font=cmsl9 scaled 1050
\ifplaintex
\headline{\vbox to 0pt{\vskip -4.5mm\line{\small\phead\ifnum
\count0=\startpage ISSN 1364-0380 (on line)
1465-3060 (printed) \hfill {\pnum\folio}\else\ifodd\count0\def\\{ }%
\ifx\theshorttitle\relax\thetitle\else\theshorttitle\fi\hfill{\pnum\folio}
\else\def\\{ and }{\pnum\folio}\hfill\ifx\theshortauthors\relax\theauthors
\else\theshortauthors\fi\fi\fi}\vss}}
\footline{\vbox to 0pt{\vglue 0mm\line{\small\pfoot\ifnum\count0=\startpage
\copyright\ \gtp\hfill\else
\gt, Volume \thevolumenumber\ (\thevolumeyear)\hfill\fi}\vss
}}
\else
\makeatletter
\def\@oddhead{{\small\lhead\ifnum\count0=\startpage ISSN 1364-0380 (on line)
1465-3060 (printed) \hfill {\lnum\number\count0}\else\ifodd\count0
\def\\{ }\ifx\theshorttitle\relax \thetitle \else\theshorttitle\fi\hfill
{\lnum\number\count0}\else\def\\{ and }{\lnum\number\count0}
\hfill\ifx\theshortauthors\relax 
\theauthors\else\theshortauthors\fi\fi\fi}}\def\@evenhead{\@oddhead}
\def\@oddfoot{\small\lfoot\ifnum\count0=\startpage\copyright\ \gtp\hfill\else
\gt, Volume \thevolumenumber\ (\thevolumeyear)\hfill\fi}
\def\@evenfoot{\@oddfoot}
\makeatother
\fi

\newwrite\gtoutfile
\long\gdef\makeheadfile{  
{\def\\{, }\def\s{ }
\immediate\openout\gtoutfile head.xxx
\immediate\write\gtoutfile{Proxy-for: \ifx\theasciiauthors\relax
\theauthors\else\theasciiauthors\fi\s<\ifx\theasciiemail\relax\theemail\else\theasciiemail\fi>}
\immediate\write\gtoutfile{\noexpand\\}
\immediate\write\gtoutfile{Authors: \ifx\theasciiauthors\relax
\theauthors\else\theasciiauthors\fi}
\immediate\write\gtoutfile{Title: \ifx\theasciititle\relax
\thetitle\else\theasciititle\fi}
\immediate\write\gtoutfile{Subj-class: GT or SG or MG etc}
\immediate\write\gtoutfile{MSC-class: \theprimaryclass\ifx\thesecondaryclass\relax\else, \thesecondaryclass\fi}
\immediate\write\gtoutfile{Journal-ref: Geom. Topol. \thevolumenumber
(\thevolumeyear) \startpage-\finishpage}
\immediate\write\gtoutfile{Comments: Published by Geometry and Topology at}
\immediate\write\gtoutfile{\s\s http://www.maths.warwick.ac.uk/gt/GTVol\thevolumenumber/paper\thepapernumber.abs.html}
\immediate\write\gtoutfile{\noexpand\\}
\immediate\write\gtoutfile{}
\ifx\theasciiabstract\relax
\immediate\write\gtoutfile{\theabstract}\else
\immediate\write\gtoutfile{\theasciiabstract}\fi
\immediate\write\gtoutfile{}
\immediate\write\gtoutfile{\noexpand\\}
\immediate\write\gtoutfile{}
\immediate\closeout\gtoutfile}}  

\def\maketitlepage{\maketitlep\makeheadfile}
\let\maketitle\maketitlepage

\lognumber{263}
\volumenumber{8}\papernumber{4}\volumeyear{2004}
\pagenumbers{115}{204}
\received{10 June 2002}
\published{8 February 2004}
\accepted{16 January 2004}
\proposed{Robion Kirby}
\seconded{Vaughan Jones, Joan Birman}

\usepackage{%
amssymb,epic,eepic,epsfig,amscd,pb-diagram,lamsarrow,pb-lams,amsmath}

\theoremstyle{plain}

\newtheorem{theorem}{Theorem}
\newtheorem{proposition}{Proposition}[section]
\newtheorem{lemma}[proposition]{Lemma}
\newtheorem{corollary}[proposition]{Corollary}

\theoremstyle{definition}
\newtheorem{definition}[proposition]{Definition}

\theoremstyle{remark}
\newtheorem{example}[proposition]{Example}

\newtheorem{remark}[proposition]{Remark}

\newcommand{\psdraw}[2]
         {\begin{array}{c} \hspace{-1.3mm}
	\raisebox{-4pt}{\epsfig{figure=draws/#1.eps,width=#2}}
	\hspace{-1.9mm}\end{array}}

\newlength{\standardunitlength}
\newcommand{\tr}{\operatorname{tr}}

\def\lbl#1{\label{#1}}

\def\BN{\mathbb N}
\def\BZ{\mathbb Z}
\def\BQ{\mathbb Q}
\def\BR{\mathbb R}

\def\A{\mathcal A}
\def\B{\mathcal B}

\def\D{\Delta}

\def\cL{\mathcal L}

\def\O{\mathcal O}

\def\N{\mathcal N}

\def\R{\mathcal R}

\def\La{\Lambda}

\def\Ga{\Gamma}
\def\S{\Sigma}

\def\ihs{integral homology 3-sphere}
\def\qhs{rational homology 3-sphere}

\def\fti{finite type invariant}

\def\la{\langle}
\def\ra{\rangle}
\def\lp{(}
\def\rp{)}

\def\w{\omega}
\def\winf{\omega}   
\def\winfp{\omega'} 
\def\wgp{\omega^{\mathrm{gp}}}

\def\logdet{\log'}

\def\a{\alpha}

\def\g{\gamma}

\def\e{\epsilon}
\def\Ga{\Gamma}

\def\b{\beta}
\def\binf{\beta}  
\def\bgp{\beta^{\mathrm{gp}}}
\def\bgpp{\beta^{\mathrm{gp}'}}

\def\Th{\Theta}
\def\s{\sigma}

\def\lk{\mathrm{lk}}

\def\pt{\partial}

\def\sminus{\smallsetminus}
\def\sfm{\mathsf{m}}

\def\ti{\widetilde}

\def\line{\operatorname{\psfig{figure=draws/line.eps,height=0.1in}}}

\def\tcircle{\circlearrowleft}
\def\tline{\uparrow}
\def\ostar{\circledast}

\def\strutb#1#2#3{\overset{#1}{\underset{#2}{ 
\begin{array}{c} \vspace{0.0cm}
\uparrow 
\vspace{-0.40cm} \\        
| \vspace{-0.45cm} \\      
\bullet \vspace{0.00cm}   
\end{array} }}\! #3}

\def\st#1#2{\overset{#1}{\underset{#2}{
\begin{array}{c} \hspace{-1.3mm}
	\raisebox{-4pt}{\psfig{figure=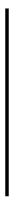,height=0.2in} }
	\hspace{-1.9mm}\end{array} }}}

\def\stbup#1#2#3{\strutb{#1}{#2}{#3}}

\def\stbdown#1#2#3{\overset{#1}{\underset{#2}{
\begin{array}{c} 
\bullet 
\vspace{-0.45cm} \\        
| \vspace{-0.40cm} \\      
\downarrow \vspace{-0.08cm}   
\end{array} }}\! #3}

\def\AS{\mathrm{AS}}
\def\IHX{\mathrm{IHX}}

\def\mat#1#2#3#4{\left[
\begin{matrix}
 #1 & #2  \\
 #3 & #4   
\end{matrix}
\right]}

\def\sfm{\mathsf{m}}
\def\sff{\mathsf{f}}

\def\ygraph{$\mathrm{Y}$-graph}

\def\Herm{\mathrm{Herm}}
\def\inn{\mathrm{Inn}}
\def\Sym{\mathrm{Sym}}
\def\LongHopf{\mathrm{LongHopf}}
\def\CA{\mathrm{String}}

\def\Gl#1{\mathcal G^{\mathrm{null}}_{#1}}

\def\Lhat{\hat\Lambda}
\def\Lloc{\Lambda_{\mathrm{loc}}}
\def\Lnew{\widehat{{\mathcal L}_{\mathrm{loc}}}}
\def\Lmore{\widehat{{\mathcal L\mathcal L}_{\mathrm{loc}}}}
 
\def\loc{\mathrm{loc}}

\def\Zrat{Z^{\mathrm{rat}}}
\def\Zratt{Z^{\mathrm{rat},t}}
\def\Zratgp{Z^{\mathrm{rat,gp}}}
\def\Zratc{\check{Z}^{\mathrm{rat}}}
\def\Z{Z}

\def\Zc{\check{Z}}

\def\NO{\N(\O)}
\def\NXO{\N_X(\O)}

\def\Bla{\B(\La\to\BZ)}
\def\intrat{\int^{\mathrm{rat}}}

\def\hair{\mathrm{Hair}}
\def\hairnu{\mathrm{Hair}^{\nu}}

\def\con{\mathrm{con}}
\def\cov{\mathrm{cov}}
\def\conh{\mathrm{con}_{\{h\}}}
\def\div{\mathrm{div}}
\def\diag{\mathrm{diag}}
\def\wrapping{wrapping}

\def\degh{\mathrm{deg}_h}

\def\phiup#1{\phi^{(#1)}}
\def\heta{\widehat{\eta}}
\def\chiz{\chi_h}
\def\eyes{\operatorname{\psfig{figure=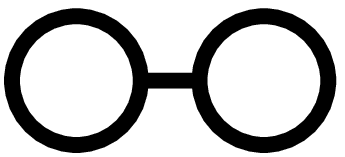,height=0.1in}}}
\def\longto{\longrightarrow}
\def\SO{(S^3,\O)}
\def\sfY{\operatorname{\mathsf{Y}}}

\def\cD{{\mathcal D}^{\times}}
\def\cDarc{{\mathcal D}^{\times,\mathrm{arc}}}

\def\UAO{{\uparrow}\hspace{-0.18cm}\mbox{x}}
\def\DAO{{\downarrow}\hspace{-0.18cm}\mbox{x}}
\def\es{\exp}
\def\ess{\exp}
\def\GA{\A^{\mathrm{gp}}}
\def\GAz{\A^{\mathrm{gp},0}}
\def\IPR{\mathrm{Int}}
\def\IA{\mathrm{Int}_{X'} \A(\star_X,\Lloc)}
\def\IGA{\mathrm{Int}_{X'} \GAXM}
\def\GAX{\GA(\star_X,\Lloc)}
\def\GAXM{\GA(\ostar_{\,X},\Lloc)}

\def\crossed{crossed\ }

\def\Ga{\Gamma}

\def\pt{\partial}

\def\oll{}

\def\gau#1#2#3{\strutb{#2}{#3}{#1}}

\begin{document}


\title{A rational noncommutative invariant of\\boundary links}

\author{\normalsize Stavros Garoufalidis and Andrew Kricker}
\shortauthors{Stavros Garoufalidis and Andrew Kricker}
\asciiauthors{Stavros Garoufalidis\\Andrew Kricker}
\coverauthors{Stavros Garoufalidis\\Andrew Kricker}
\address{School of Mathematics, Georgia Institute of 
Technology\\Atlanta, GA 30332-0160, USA}
\gtemail{\mailto{stavros@math.gatech.edu}\qua{\rm and}\qua\mailto{akricker@math.toronto.edu}}
\secondaddress{Department of Mathematics, University of 
Toronto\\Toronto, Ontario, Canada M5S 3G3}

\asciiaddress{School of Mathematics, Georgia Institute of 
Technology\\Atlanta, GA 30332-0160, USA\\and\\Department of Mathematics, 
University of Toronto\\Toronto, Ontario, Canada M5S 3G3}

\asciiemail{stavros@math.gatech.edu, akricker@math.toronto.edu}

\primaryclass{57N10}\secondaryclass{57M25}
\keywords{Boundary links, Kontsevich integral, Cohn
localization}

\begin{abstract}
In 1999, Rozansky conjectured the existence of a rational presentation
of the Kontsevich integral of a knot. Roughly speaking, this rational 
presentation
of the Kontsevich integral would sum formal power series into rational
functions with prescribed denominators. Rozansky's conjecture was soon
proven by the second author. We begin our paper by reviewing 
Rozansky's conjecture and the main ideas that lead to its proof.
The natural question of extending this conjecture to links leads to
the class of boundary links, and a proof of Rozansky's conjecture in this
case. A subtle issue is the fact that a `hair' map which replaces beads
by the exponential of hair is not 1-1. This raises the question of whether
a rational invariant of boundary links exists in an appropriate space 
of trivalent graphs whose edges are decorated by rational functions in
noncommuting variables. A main result of the paper is to construct such
an invariant, using the so-called surgery view of boundary links and
after developing a formal diagrammatic Gaussian integration.

Since our invariant is one of many rational forms of the Kontsevich integral,
one may ask if our invariant is in some sense canonical. We prove that this
is indeed the case, by axiomatically characterizing our invariant as a 
universal finite type invariant of boundary links with respect to the null 
move.
Finally, we discuss relations between our rational invariant and homology
surgery, and give some applications to low dimensional topology.
\end{abstract}

{\small\maketitlepage}


\section{Introduction}
\lbl{sec.intro}

\subsection{Rozansky's conjecture on a rational presentation of the
Kontsevich integral of a knot}
\lbl{sub.rozansky}

The Kontsevich integral of a knot is a powerful invariant that can be
interpreted to take values in a completed vector space of graphs. The graphs 
in question have univalent and trivalent vertices only 
(so-called {\em unitrivalent graphs}), and are considered modulo some 
well-known relations that include the $\AS$ and $\IHX$ relations. 

Every unitrivalent graph $G$ is the union of a trivalent graph $G^t$ together 
with a number of unitrivalent trees attached on the edges  of $G^t$:
$$
\psdraw{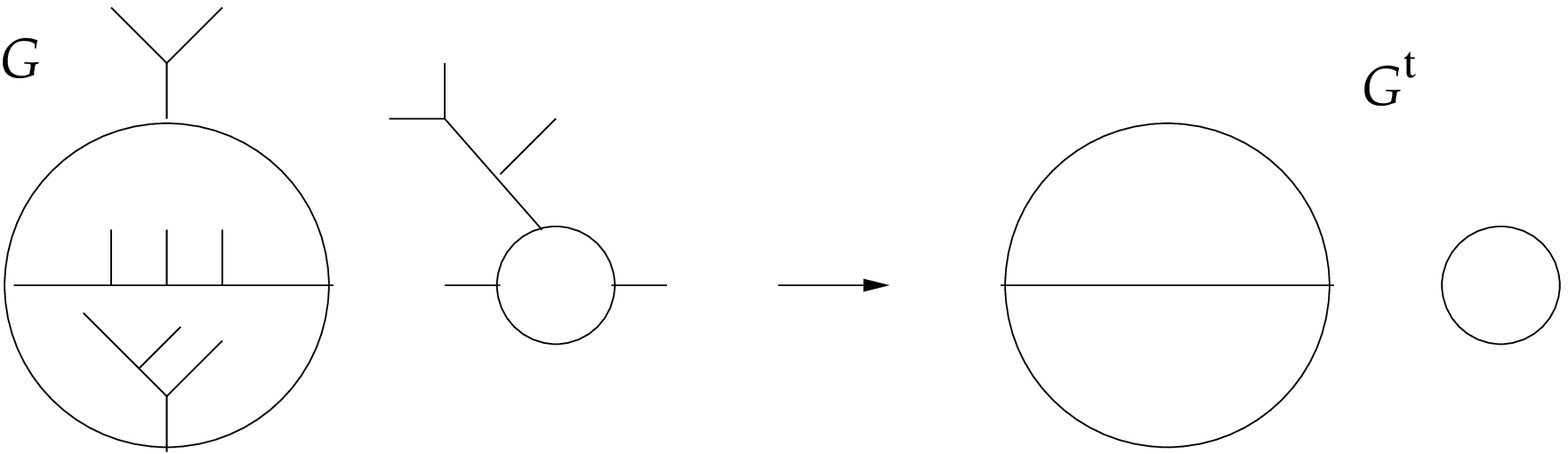}{3in}
$$
This is true provided that no component of $G$ is a tree, and provided 
that we allow $G^t$ to include circles in case $G$ contains components
with one loop (ie, with  betti number $1$).

The $\AS$ relation kills all trees with an internal trivalent vertex, which
are possibly attached on an edge of a trivalent graph. Thus, the only
trees that survive are the {\em hair}, that is the trees with one edge
and two univalent vertices.

As a result, we need only consider trivalent graphs with a number of hair
attached on their edges. (The exceptional case of a single hair unattached
anywhere is excluded since we are silently assuming that the knot is 
zero-framed).
The number of hair may be recorded by a monomial in a variable $h$ attached 
on each edge of a trivalent graph, together with an orientation of the edge
which keeps track of which side of the edge should the hair grow. For
example, we have:
$$
\psdraw{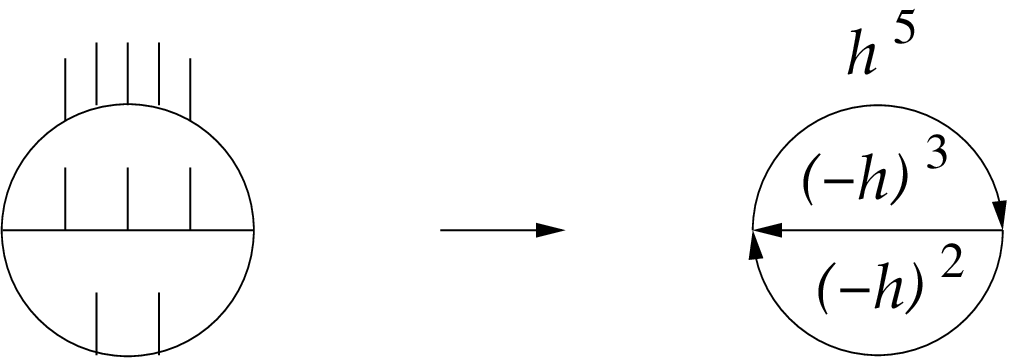}{2.5in}
$$
By linearity, we may decorate edges of trivalent graphs by polynomials,
and even by formal power series in $h$.

In the summer of 1999 Rozansky made the bold conjecture that the Kontsevich
integral of a knot can be interpreted to take values in a space of trivalent
graphs with edges decorated by {\em rational functions in $e^h$}. Moreover,
the denominators of these rational functions ought to be the Alexander
polynomial of a knot.

Rozansky's conjecture did not come out of the blue. It was motivated by
earlier work of his on the colored Jones function; see \cite{R1}. In
that reference, Rozansky proved that the colored Jones function of a knot
(a certain power series quotient of the Kontsevich integral) can be
presented as a power series of rational functions whose denominators
where powers of the Alexander polynomial.

\subsection{Kricker's proof of Rozansky's Conjecture}
\lbl{sub.Krickerproves}

Shortly after Rozansky's Conjecture appeared, the second author gave a proof
of it in \cite{Kr1}. Since the proof contains several ideas that are key to
the results of the present paper, we would like to summarize them here. 

{\bf Fact 1}\qua Untie the knot.

The key idea behind this is the fact that knots (or rather, knot projections)
can be unknotted via a sequence of crossing changes, and that a crossing
change can be achieved by surgery on a $\pm 1$-framed unknot as follows:

\begin{figure}[ht!]
$$ 
\psdraw{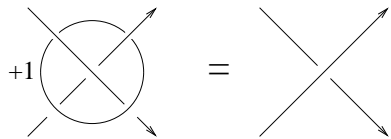}{2in} 
$$
\caption{A crossing change can be achieved by surgery on a unit framed
unknot.}\lbl{split}
\end{figure}

Thus, every knot $K$ in $S^3$ can be obtained by surgery on a framed link $C$
in the complement $S^3\sminus\O$ of a standard unknot $\O$. We will call such 
a link $C$, an {\em untying link} for $K$. For example, an untying link for
the Figure 8 knot is:
$$
\psdraw{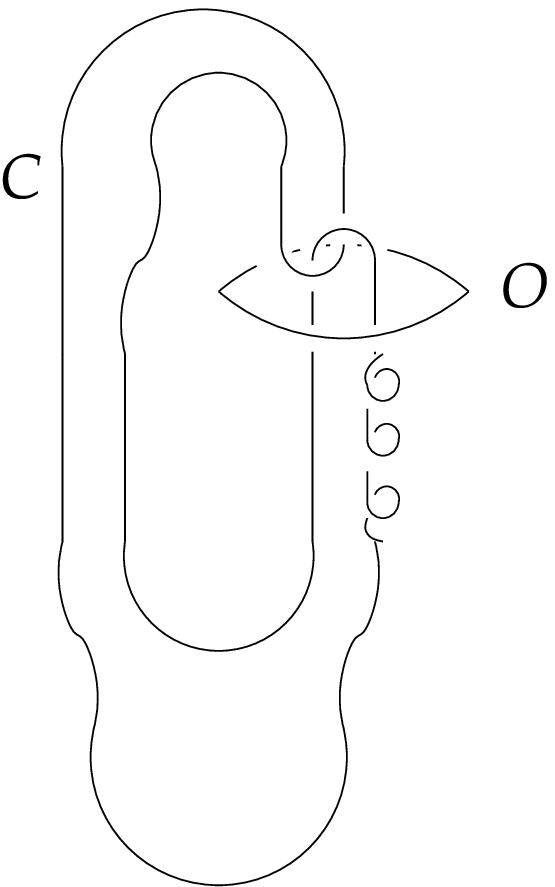}{1in}
$$
Observe that untying links are framed, and null homotopic in $S^3\sminus\O$ 
(in the sense that every component is contractible in $S^3\sminus\O$), the 
interior of a solid torus. Observe further that untying links exist for every 
knot $K$ in an \ihs\ $M$.

{\bf Fact 2}\qua Compute the Kontsevich integral of a knot from the Aarhus
integral of an untying link.

By this, we mean the following. Consider a knot $K$ and an untying link $C$
in $S^3\sminus\O$. Then, we may consider the (normalized) Kontsevich integral 
$\Zc(C \cup \O)$ of the link $C \cup \O$, which can be interpreted to take
values in a completed vector space of unitrivalent graphs with legs decorated
by the components of $C \cup \O$. Then, the Kontsevich integral $Z(K)$ of $K$
can be computed from $\Zc(C \cup \O)$ by 
$$
Z(K)=\int dC \, (\Zc(C \cup \O))
$$
Here $\int dC$ refers to a {\em diagrammatic formal Gaussian integration} 
which roughly speaking glues pairwise the $C$-colored legs of the graphs in
$\Zc(C \cup \O)$ using the negative inverse linking matrix of $C$; see
\cite{A}.

{\bf Fact 3}\qua Compute the Kontsevich integral of an untying link from
the Kontsevich integral of a Long Hopf Link.

By this, we mean the following. The following formula for the Kontsevich 
integral of a {\em Long Hopf Link}, was conjectured in \cite{A0} (in 
conjunction with the so-called {\em Wheels} and {\em Wheeling} Conjectures) 
and proven in \cite{BLT}:
$$
\Z \left(\psdraw{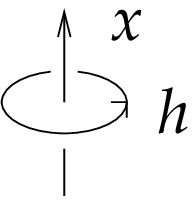}{0.35in} 
\right) = 
\strutb{x}{}{e^h}
\sqcup \nu(h)
$$
where $\nu(h)$ is the Kontsevich integral of the unknot, expressed in terms
of graphs with legs colored by $h$. $\sqcup$ refers to the disjoint union
multiplication of graphs. In this formula a {\em bead} $e^h$
(that is, the exponential of hair), appears explicitly.

Using locality of the Kontsevich integral of a link, we may slice a planar
projection of $C \cup \O$ into a tangle $T$ by cutting $C \cup \O$ along the
meridianal disk in the solid torus $S^3\sminus\O$ that $\O$ bounds.
Then, we can compute $\Zc(C\cup\O)$ from $\Zc(T)$ after we glue in the beads
as instructed by the formula for the Long Hopf Link.

The result of this step is that we managed to write the Kontsevich integral
of an untying link in terms of unitrivalent graphs with edges decorated
by Laurent polynomials in $e^h$, and with univalent vertices colored by $C$.

{\bf Fact 4}\qua Commute the Aarhus integration with the above formula of
the Kontsevich integral of an untying link.

Modulo some care with the 1-loop part of the Aarhus integral, this proves
Rozansky's Conjecture for knots.

\subsection{Rozansky's conjecture for boundary links}
\lbl{sub.rozflinks}

A question arises immediately after Kricker's proof. The Kontsevich integral
is defined not just for knots, but for all links. Is there a more general
class of links that this proof can be applied to?

In order to answer this, let us recall that the Kontsevich integral of a link 
takes values in a completed vector space of unitrivalent graphs, whose
univalent vertices are labeled by the components of the link. The graphs
are considered modulo some well-known relations, that include
a 3-term relation (the $\IHX$ relation) depicted as follows:
$$
\psdraw{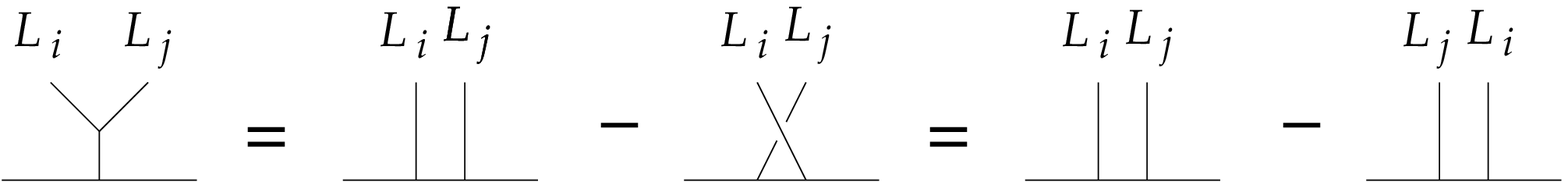}{5in}
$$
Using this relation, it follows that every unitrivalent graph with no tree 
components is a sum of trivalent graphs with {\em hair} attached on their 
edges. Moreover, the hair is now labeled by the components of the link.

Thus, in order to formulate a conjecture for the Kontevich integral of links
along the lines of Rozansky, we need to restrict attention to links whose 
Kontsevich integral has no tree part.
For such links, the Kontsevich integral takes values in a completed vector
space generated by trivalent graphs with hair. The hair are colored by
the components of the link, and we can record this information by placing
monomials in variables $h_i$ (one variable per link component)
on the edges of a trivalent graph as follows:
$$
\psdraw{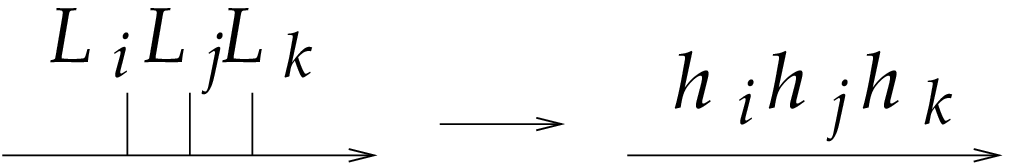}{2.5in} 
$$
By linearity, we can place polynomials, and even formal power series,
in the noncommuting variables $h_i$ on the edges of trivalent graphs.

The question arises: are there any links whose Kontsevich integral has 
vanishing tree-part? Using Habegger-Masbaum (see \cite{HM}), this condition 
is equivalent to the {\em vanishing of all Milnor invariants} of a link.
This class of links contains (and perhaps coincides with) the class of 
{\em sublinks of homology boundary links}, \cite{Le2}. For simplicity, we 
will focus on the class of {\em boundary links}, namely those each component 
of which bounds a surface, such that the surfaces are pairwise disjoint.

{\bf Fact 5}\qua Boundary links can be untied.

Indeed the idea is the following. Choose a Seifert surface as above for a 
boundary link $L$; that is, one surface per component, so that the surfaces 
are pairwise
disjoint. Then, do crossing changes among the bands of the Seifert surface
to unknot each band and unlink them. The result is a standard unlink, and
a framed nullhomotopic link $C$, such that surgery on $C$ transforms the
unlink to the boundary link. We will call such a link $C$, an {\em untying
link} for $L$.

Using this fact, and repeating Facts 2-4 for boundary links, proves Rozansky's
conjecture for boundary links. In Fact 3, the Kontsevich integral of an 
untying link $C \cup\O$ takes values in a completed space of trivalent graphs
whose edges are labeled by Laurent polynomials in noncommuting variables 
$e^{h_i}$, for $i=1,\dots,g$ where $g$ is the number of components of $L$.
In Fact 4, care has to be taken to make sense of rational functions in
noncommuting variables $e^{h_i}$.

\subsection{Is there a rational form of the Kontsevich integral of
boundary links?}
\lbl{sub.isthere}

Our success in proving a rational presentation for the Kontsevich integral
of a boundary link suggests that 

\begin{itemize}
\item
we may try to define an invariant $\Zrat$ of boundary 
links with values in a completed space $\A(\Lloc)$ of trivalent graphs with 
beads, where the beads are rational functions in noncommuting variables
$t_i$, for $i=1,\dots,g$. 
\item
Rational functions in noncommuting variables
ought to form a ring $\Lloc$, which should be some kind of localization
of the {\em group ring} $\La=\BZ[F]$ of the free group $F$ with generators
$t_1,\dots,t_g$.
\item
There ought to be a $\hair$ map from graphs with beads to unitrivalent
graph with hair which replaces $t_i$ by $e^{h_i}$, such that
the Kontsevich integral $Z$ is given by $\hair \circ \Zrat$.
\end{itemize}
 
Let us call the above statement a {\em strong form of Rozansky's Conjecture}
for boundary links, and let us call any such invariant $\Zrat$ {\em a rational
form} of the Kontsevich integral.

Let us point out a subtlety (easy to miss) of this stronger conjecture, even 
in the case
of knots with trivial Alexander polynomial (in which case we may work with
graphs with beads in $\La=\BZ[t^{\pm}]$ only): although it is true that
the map
$$
\La=\BZ[t^{\pm}] \longto \Lhat=\BQ[h]
$$
given by $t \mapsto e^h$ is 1-1, it does not follow in some obvious way
that the $\hair$ map from the space of graphs with beads to the space of 
unitrivalent
graphs is 1-1. This may seem counterintuitive, however a recent paper of
Patureau-Mirand (using Vogel's work on universal algebras and exotic weight
systems that do not come from Lie algebras) proves that this $\hair$ map is 
not 1-1; see \cite{P}.

If the $\hair$ map were 1-1, then Rozansky's Conjecture for boundary links
(which we proved in Section \ref{sub.rozflinks}
would easily imply the existence and uniqueness of a rational form $\Zrat$
of the Kontsevich integral. However, as we discussed above, this is not the 
case.

This raises two problems:
\begin{itemize}
\item
Is there a rational form $\Zrat$ of the Kontsevich integral of a boundary link?
\item
Assuming there is one, is there a canonical (in some sense) form?
\end{itemize}

The purpose of the paper is to solve both problems.

\subsection{The main results of the paper}
\lbl{sub.maini}

In the remainder of this introduction, let us discuss the main
results of this paper, which we will explain in lengthy detail in the 
following 
sections.

{\bf Fact 6}\qua A surgery view of boundary links.

We mentioned already in Fact 5 tat boundary links can be untied. Unfortunately,
an untying link of a boundary link is not unique. In fact, Kirby moves
on an untying link do not change the result of surgery on an untying link,
and therefore give rise to the same boundary link.
It turns out that Kirby moves preserve not only the boundary link but 
its $F$-structure as well. 

By this we mean the following. The choice of Seifert surface is not part of
a boundary link. 
A refinement of {\em boundary links} (abbreviated $\pt$-links) are
{\em $F$-links}, ie, a triple $(M,L,\phi)$ of a link $L$ in an \ihs \
$M$ and an onto map $\phi:\pi_1(M- L)\twoheadrightarrow F$, where $F$ is the 
free group on a set $T=\{t_1,\dots,t_g\}$ 
(in 1-1 correspondence with the components of $L$)
such that the $i$th meridian is sent to the $i$th generator. Two $F$-links
$(M,L,\phi)$ and $(M',L',\phi')$ are equivalent iff $(M,L)$ and $(M',L')$
are isotopic and the maps $\phi$ and $\phi'$ differ by an inner automorphism
of $F$. The underlying link of an $F$-link is a $\pt$-link. Indeed, an
$F$-link gives rise to a map $\tilde\phi: M \longto \vee^g S^1$ which
induces the map $\phi$ of fundamental groups. Choose generic points
$p_i$, one in each circle of $\vee^g S^1$. It follows by transversality
that $\tilde\phi^{-1}(p_i)$ is a surface with boundary component the 
corresponding component of $L$. These surfaces are obviously pairwise disjoint,
thus $(M,L)$ is a boundary link.

There is an action
$$ 
\CA_g \times \text{$F$-links} \longto \text{$F$-links}
$$
whose orbits can be identified with the set of boundary links.
Here $\CA_g$ is the {\em group of group of motions of an unlink in 3-space}, 
and can be identified with the automorphisms of the free group that map
generators to conjugates of themselves:
$$
\CA_g \cong \{ f \in \text{Aut}(F) 
| f(t_i)=\a_i^{-1} t_i \a_i, \quad i=1,\dots,g \}.
$$
The action of $\CA_g$ on an $F$-link $(M,L,\phi)$ is given by composition
with the map $\phi$, as is explained in Section \ref{sub.surgeryb}.

Let $\NO$ denote the set of nullhomotopic framed links $C$ in the complement 
of a standard unlink $\O$ in $S^3$, such that the linking matrix of $C$
is invertible over $\BZ$.

Then, in \cite{GK2} we prove that
the surgery map induces a 1-1 correspondence 
$$
\NO/\la \text{Kirby moves}, \CA_g \ra \leftrightarrow \pt-\text{links}.
$$
This is the so-called surgery view of boundary links.

{\bf Fact 7}\qua Construct an invariant of $\NO$ in a space of graphs with beads.

By analogy with Fact 3, 
in Section \ref{sec.tangles} we define an invariant of links $C \in \NO$ 
as above. It takes values in a completed vector space of univalent graphs 
with beads. The beads are elements of $\La=\BZ[F]$, and record the winding of 
a tangle representative of $C$ in $S^3-\O$. The legs of the graphs are 
labeled by the components of $C$. The strut part of this invariant records 
the equivariant linking matrix of $C$.

{\bf Fact 8}\qua Develop an equivariant version of the Aarhus integral.

Unfortunately, the above described invariant of links $C \in \NO$ is not 
invariant under Kirby moves of $C$. To accommodate for that, in Section 
\ref{sec.aarhus} we construct an integration theory $\intrat$ in the spirit of 
the Aarhus Integral \cite{A}. Let us mention that the Aarhus integration
is a map that cares only about the univalent vertices of graphs, and not
about their internal structure (such as beads on edges, or valency of
internal vertices). By construction, Aarhus integration behaves like 
integration of functions in the sense that it behaves well with changes of
variables. 

In our integration theory $\intrat$, we separate the 
strut part of a diagram with
beads, invert their martix (it is here that a suitable ring $\Lloc$ is 
needed) to construct new struts, and glue the legs
of these new struts to the rest of the diagrams. The result is a formal
linear combination of trivalent graphs whose beads are elements of $\Lloc$.
The main property of this integration theory, is that the resulting 
invariant is independent under Kirby moves and respects the $\CA_g$ action. 
Thus it gives an invariant $\Zrat$ of $\pt$-links which takes values in a 
completed
vector space of trivalent graphs with edges decorated by $\Lloc$, modulo
certain natural relations; see Theorem \ref{thm.1}.

{\bf Fact 9}\qua The ring $\Lloc$ of rational functions in noncommuting variables.

Our integration theory $\intrat$ reveals the need for a ring $\Lloc$.
This ring should satisfy the property that all matrices $W$ over $\La$
which are invertible over $\BZ$ (that is, $\e W$ is invertible over $\BZ$
where $\e: \La=\BZ[F] \longto \BZ$ is the map that sends $t_i$ to $1$)
are in fact invertible over $\Lloc$. This is precisely the defining property
of the {\em noncommutative (Cohn) localization} of $\La=\BZ[F]$. 
Farber-Vogel identified $\Lloc$ with the ring of {\em rational functions in 
noncommuting variables}; see \cite{FV}. 

{\bf Fact 10}\qua $\Zrat$ is a rational form of the Kontsevich integral; see
Theorem \ref{thm.1} in Section \ref{sec.comparison}.

By analogy with Fact 2, we compare the $\Zrat$ invariant of boundary links
with their Kontsevich integral, via the $\hair$ map. The comparison is achieved
using the formula of the Kontsevich integral of the Long Hopf Link.
Section \ref{sec.comparison} completes Fact 10.

So far out efforts give a rational form of the Kontsevich integral of a link.
As we mentioned before, there is potentially more than one rational form
of the Kontsevich integral. How do we know that our construction is in some
sense a natural one?

To answer this question, let us recall that the Kontsevich integral is 
a {\em universal} finite type invariant of links, with respect to the 
Goussarov-Vassiliev crossing change moves. This axiomatically characterizes 
the Kontsevich integral, up to a universal constant which is independent of a 
link.

{\bf Fact 11}\qua The $\Zrat$ invariant is a universal finite type invariant
of $F$-links with respect to the null-move; see Theorem \ref{thm.KOK}
in Section \ref{sub.stateuniversal}.

This so called {\em null-move} is described in terms of surgery on a 
nullhomotopic clasper in the complement of the $F$-link, and generalizes
the null-move on knots. The latter was 
introduced and studied extensively by \cite{GR}.

Section \ref{sec.universal} is devoted to the proof of Theorem \ref{thm.KOK},
namely the universal property of $\Zrat$. This follows from the 
general principle of  locality of the invariant $\Zrat$ (ie, with the 
behavior of $\Zrat$ under sublinks) together with the identification of the 
covariance of $\Zrat$ with an equivariant linking function.

\subsection{Some applications of the $\Zrat$ invariant}
\lbl{sub.applications}

The rational form $\Zrat$ of the Kontsevich integral of a boundary link
sums an infinite series of Vassiliev invariants into rational functions.
An application of the universality property of the $\Zrat$ invariant
is a {\em realization theorem} for these rational functions in Section 
\ref{sec.universal}; see Proposition \ref{prop.MKu2}. Our realization 
theorem is in the same spirit as the realization theorem for the 
values of the Alexander polynomial of a knot that was achieved decades
earlier by Levine; see \cite{Le1}.

An application of this realization theorem for the $2$-loop part of the
$\Zrat$ invariant settles 
a question of low-dimensional topology (namely, to
separate minimal rank knots from Alexander polynomial $1$ knots), and fixes 
an error in earlier work of M. Freedman; see \cite{GT}. This is perhaps the 
first application of finite type invariants in a purely 3-dimensional 
question.

Another application of the $2$-loop part is a formula for the 
Casson-Walker invariant of cyclic branched coverings of a knot in terms of 
the signature of the knot and residues of the $Q$ function, obtained in
joint work \cite{GK1}.

\subsection{Future directions}
\lbl{sub.future}

Rozansky has recently introduced a Rationality
Conjecture for a class of {\em algebraically connected} links, that is links
with nonvanishing Alexander polynomial, \cite{R4}. This class is disjoint 
from ours, since
a boundary link with more than one component has vanishing Alexander 
polynomial. It is an interesting question to extend our results in this
setting of Rozansky. We plan to address this issue in a later publication.

\subsection{Acknowledgements} 
The results of this work were first announced by the second named
author in a meeting on ``Influence of Physics on Topology'' on August
2000, in the sunny San Diego.  We thank the organizers, M Freedman and
P Teichner for arranging such a fertile meeting.  We also wish to
express our gratitude to J Levine who shared with the first author
his expertise on topology, and to J Hillman and T Ohtsuki for
encouragement and support of the second author. Finally, the authors
wish to thank D Bar-Natan and the referee for numerous comments who
improved the presentation of the paper.

The first author was partially supported by an NSF grant DMS-98-00703
and by an Israel-US BSF grant.  The second author was partially
supported by a JSPS fellowship.

\section{A surgery view of $F$-links}
\lbl{sec.surgeryview}

\subsection{A surgery description of $F$-links}
\lbl{sub.surgery}

In this section we recall the surgery view of links, introduced in \cite{GK2},
which is a key ingredient in the construction of the rational invariant
$\Zrat$. Let us recall some important notions.
Fix once and for all a based unlink $\O$ of $g$ components in $S^3$.

\begin{definition}
\lbl{def.NO}
Let $\NO$ denote the set of {\em nullhomotopic links $C$ with 
$\BZ$-invertible linking matrix} in the complement of $\O$.
\end{definition}

Surgery on an element $C$ of $\NO$ transforms $(S^3,\O\,\phi)$ to a $F$-link 
$(M,L,\phi)$. Indeed, since $C$ is nullhomotopic the natural map 
$\pi_1(S^3\sminus\O)\to F$ 
gives rise to a map $\pi_1(M\sminus L)\to F$. Alternatively, 
one can construct disjoint Seifert surfaces for each component of $L$ by 
tubing the disjoint disks that $\O$ bounds, which is possible, since each 
component of $C$ is nullhomotopic. 

Since the linking matrix of $C$ is invertible over $\BZ$, $M$ is an \ihs .
Let $\stackrel{\kappa}\sim$ denote the equivalence relation on $\NO$ 
generated by the moves of {\em handle-slide} $\kappa_1$ 
(ie, adding a parallel of a link component to another component) and 
{\em stabilization} $\kappa_2$ (ie, 
adding to a link an unknot away from the link with framing $\pm 1$). 
It is well-known that $\kappa$-equivalence preserves surgery. A main result 
of \cite{GK2} is the following:

\begin{theorem}{\rm\cite{GK2}}\qua
\lbl{thm.GK2} 
The surgery map gives a 1-1 and onto correspondence
$$
\NO/\la \kappa \ra \to \text{$F$-links} .
$$ 
\end{theorem}

As an application of the above theorem, in \cite{GK2} we constructed a map:
$$
W: \text{$F$-links} \longto \Bla
$$
where $\Bla$ is the set of simple stably congruence class of Hermitian
matrices $A$, invertible over $\BZ$. Here, a matrix $A$ over $\La$ is 
{\em Hermitian} if $A^\star=A$ where $A^\star$ denotes the conjugate transpose
of $A$ with respect to the involution $\La=\BZ[F]\to \La=\BZ[F]$ which
sends $g$ to $\bar g:=g^{-1}$. 
Moreover, two Hermitian matrices $A,B$ are {\em simply stably congruent} iff 
$A\oplus S_1=P(B\oplus S_2)P^\star$ for some diagonal matrices $S_1,S_2$
with $\pm 1$ entries and some {\em elementary} matrix $P$ (ie, one which
differs from the identity matrix on a single non-diagonal entry) 
or a diagonal matrix with entries in $\pm F \subset \La$.
The map $W$ sends an $F$-link to the equivariant linking matrix of $\ti C$,
a lift of a surgery presentation $C$, to the free cover of $S^3\sminus\O$.
It was shown in \cite{GK2} that the map $W$ determines the Blanchfield
form of the $F$-link, as well as a noncommutative version of the Alexander
polynomial defined by Farber \cite{Fa}.

\subsection{A tangle description of $F$-links}
\lbl{sub.tangles}

In this section we give a tangle diagram description of the set $\NO$.
Before we proceed, a remark on notation:

\begin{remark}
\lbl{rem.lara}
Throughout the paper, given an equivalence relation $\rho$ on a set $X$,
we will denote by $X/\la \rho \ra$ the set of {\em equivalence classes}.
Occasionally, an equivalence relation on $X$ will be defined by one of the
following ways: 
\begin{itemize}
\item
Either by the
{\em action} of a group $G$ on $X$, in which case $X/\la \rho \ra$ coincides 
with the set of orbits of $G$ on $X$. 
\item
Or by a {\em move} on $X$, ie, by a subset of $X \times X$, in which case
the equivalence relation is the smallest one that contains the move.
\end{itemize}
\end{remark}

Consider the following surfaces in $\BR^2$ 
$$
\setlength{\unitlength}{0.03\standardunitlength}
	\begin{array}{c}  \hspace{-1.7mm}
         	\raisebox{-8pt}{\input draws/Sigma1.epc }
         	\hspace{-1.9mm}
	\end{array}

$$
together with a distinguished part of their boundary (called a {\em gluing
site}) marked by an arrow
in the figure above. 
By a {\em tangle diagram on a surface $\S_x$} (for $x=a,b,c$) we mean an
oriented, framed, smooth, proper immersions of an
oriented 1-manifold into the surface, up to isotopies rel boundary of
the surface, double points being equipped with crossing information, with
the boundary points of the tangle lying 
on standard points at $t=0$ and $t=1$.

A {\em \crossed tangle diagram on a surface} $\S_x$ is a tangle diagram on 
$\S_x$ possibly with some crosses (placed away from the gluing sites)
that mark those boundary points of the diagram. Here are the possible places 
for a cross
on the surfaces $\S_x$ for $x=a,b$ or $c$:
$$
\psdraw{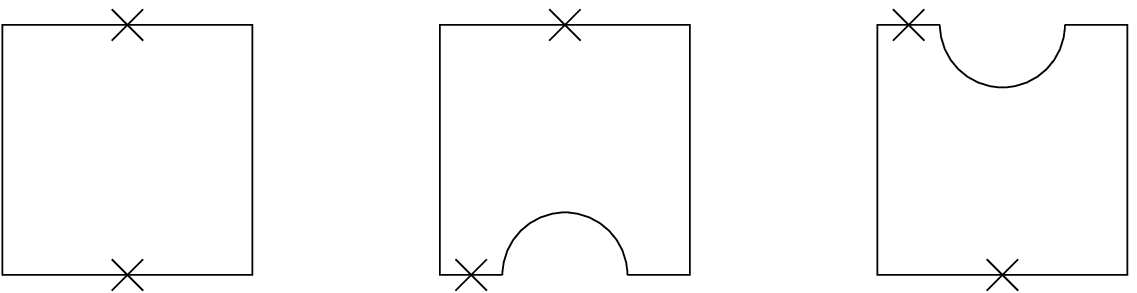}{2in}
$$
The following is an example of a crossed tangle diagram on $\S_b$:
$$
\setlength{\unitlength}{0.03\standardunitlength}
	\begin{array}{c}  \hspace{-1.7mm}
         	\raisebox{-8pt}{\input draws/Sigma2.epc }
         	\hspace{-1.9mm}
	\end{array}

$$
The cross notation evokes a small pair of scissors 
$\setlength{\unitlength}{0.02\standardunitlength}
	\begin{array}{c}  \hspace{-1.7mm}
         	\raisebox{-8pt}{\input draws/scissors.epc }
         	\hspace{-1.9mm}
	\end{array}
$ 
that will be sites where the skeleton of the tangle will be cut.

Let $D_g$ denote a standard disk with $g$ holes ordered from top to bottom
and $g$ gluing sites, shown as horizontal segments below. 
For example, when $g=2$, $D_g$ is:
$$
\psdraw{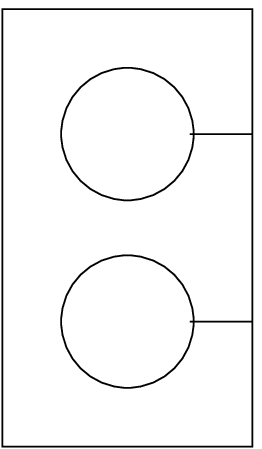}{0.5in}
$$

\begin{definition}
\lbl{def.slicedcrossedlink}
A {\em sliced crossed link} in $D_g$ is 
\begin{itemize}
\item[(a)]
a nullhomotopic link $L$ in 
$D_g \times I$ with $\BZ$-invertible linking matrix such that
\item[(b)]
$L$ is in general position with respect to the gluing sites of $D_g$ times $I$
in $D_g \times I$.
\item[(c)]
Every component of $L$ is equipped with a point (depicted as a cross $\times$)
away from the gluing sites times $I$. 
\end{itemize}
Let $\cL^{\times}(D_g)$ denote the set of isotopy classes of
sliced crossed links.
\end{definition}

\begin{definition}
\lbl{def.slicedcrosseddiagram}
A {\em sliced \crossed diagram} in $D_g$ is a sequence of \crossed
tangles, $\{T_1,\dots,T_k\}$ 
such that
\begin{itemize}
\item[(a)] the boundaries of the surfaces match up, that is, their shape,
distribution of endpoints, and crosses. 
$$
\text{For example, if $T_i$ is}
\setlength{\unitlength}{0.03\standardunitlength}
	\begin{array}{c}  \hspace{-1.7mm}
         	\raisebox{-8pt}{\input draws/Sigma3.epc }
         	\hspace{-1.9mm}
	\end{array}
 \text{ then $T_{i+1}$ is}
\setlength{\unitlength}{0.03\standardunitlength}
	\begin{array}{c}  \hspace{-1.7mm}
         	\raisebox{-8pt}{\input draws/Sigma31.epc }
         	\hspace{-1.9mm}
	\end{array}

$$
\item[(b)] the top and bottom boundaries of $T_1$ and $T_k$
look like:
$$
\setlength{\unitlength}{0.03\standardunitlength}
	\begin{array}{c}  \hspace{-1.7mm}
         	\raisebox{-8pt}{\input draws/Sigma5.epc }
         	\hspace{-1.9mm}
	\end{array}

$$
\item[(c)]
After stacking the tangles $\{T_1,\dots,T_k\}$ from top to bottom, we obtain 
a crossed
link in $D_g \times I$, where $D_g$ is a disk with $g$ holes. 
\end{itemize}
Let $\cD(D_g)$ denote the set of sliced \crossed diagrams  on $D_g$ such that
each component of the associated link is nullhomotopic in $D_g \times I$
and marked with precisely one cross. 
\end{definition}

Here is an example of a sliced crossed diagram
in $D_1$ whose corresponding link in $\NO$ is a surgery presentation for
the Figure 8 knot:
$$
\setlength{\unitlength}{0.025\standardunitlength}
	\begin{array}{c}  \hspace{-1.7mm}
         	\raisebox{-8pt}{\input draws/figure8.epc }
         	\hspace{-1.9mm}
	\end{array}

$$
Clearly, there is a map
$$
\cD(D_g) \longto \cL^{\times}(D_g).
$$
We now introduce some equivalence relations on $\cD(D_g)$ that are important
in comparing $\cD(D_g)$ to $\NO$.

\begin{definition}
\lbl{def.regularisotopy}
{\em Regular isotopy} $\stackrel{r}\sim$ of sliced \crossed diagrams is the 
equivalence relation generated by regular isotopy of individual tangles
and the following moves:
$$
\setlength{\unitlength}{0.03\standardunitlength}
	\begin{array}{c}  \hspace{-1.7mm}
         	\raisebox{-8pt}{\input draws/Sigma6.epc }
         	\hspace{-1.9mm}
	\end{array}

$$
where the glued edge of $T_e$ must have no crosses marking it.
\end{definition}

\begin{proposition}
\lbl{prop.tangles1}
There is a 1-1 correspondence: 
$$
\cD(D_g)/\la r \ra \longleftrightarrow \cL^{\times}(D_g).
$$
\end{proposition}

\begin{proof}
This follows from standard transversality arguments;
see for example \cite{AMR} and \cite[Figure 7]{Tu}. 
\end{proof}

\begin{definition}
\lbl{def.basingrelation}
{\em Basing relation}  $\stackrel{\b}\sim$ of sliced \crossed diagrams  is the 
equivalence relation generated by moving a cross of each component of the
associated link in some other admissible position in the sliced crossed tangle.
Moving the cross of a sliced crossed link can be obtained 
(in an equivalent way) by the local moves $\b_1$ and $\b_2$, where
$\b_1$ (resp. $\b_2$) moves the cross across consecutive tangles in
such a way that a gluing site is not (resp. is) crossed. Here is an
example of a $\b_2$ move:
$$
\{\dots, \setlength{\unitlength}{0.03\standardunitlength}
	\begin{array}{c}  \hspace{-1.7mm}
         	\raisebox{-8pt}{\input draws/Sigma101.epc }
         	\hspace{-1.9mm}
	\end{array}
, \setlength{\unitlength}{0.03\standardunitlength}
	\begin{array}{c}  \hspace{-1.7mm}
         	\raisebox{-8pt}{\input draws/Sigma102.epc }
         	\hspace{-1.9mm}
	\end{array}
, \dots \}
\longleftrightarrow 
\{\dots, \setlength{\unitlength}{0.03\standardunitlength}
	\begin{array}{c}  \hspace{-1.7mm}
         	\raisebox{-8pt}{\input draws/Sigma104.epc }
         	\hspace{-1.9mm}
	\end{array}
, \setlength{\unitlength}{0.03\standardunitlength}
	\begin{array}{c}  \hspace{-1.7mm}
         	\raisebox{-8pt}{\input draws/Sigma103.epc }
         	\hspace{-1.9mm}
	\end{array}
, \dots \}
$$
\end{definition}

A {\em link} $L$ in $D_g$ is one that satisfies condition (a) 
of Definition \ref{def.slicedcrossedlink} only. Let $\N(D_g)$ denote the
set of isotopy classes of sliced links in $D_g$. Proposition 
\ref{prop.tangles1} implies that 

\begin{proposition}
\lbl{prop.tangles2}
There is a 1-1 correspondence: 
$$
\cD(D_g)/\la r, \b \ra \longleftrightarrow \N(D_g).
$$
\end{proposition}

Observe that $D_g \times I$ contained in the complement of an unlink 
$S^3\sminus\O$. Moreover, up to homotopy we have that $S^3\sminus\O$
is obtained from $D_g \times I$ by attaching $g-1$ 2-spheres.
The next equivalence relation on sliced crossed diagrams is precisely
the move of sliding over these 2-spheres.

\begin{definition}
\lbl{def.wrapping}
The {\em \wrapping\ relation} $\stackrel{\w}\sim$ of sliced crossed diagrams
is generated by the following move
$$
\{\dots,\setlength{\unitlength}{0.03\standardunitlength}
	\begin{array}{c}  \hspace{-1.7mm}
         	\raisebox{-8pt}{\input draws/Sigmaw1.epc }
         	\hspace{-1.9mm}
	\end{array}
,\setlength{\unitlength}{0.03\standardunitlength}
	\begin{array}{c}  \hspace{-1.7mm}
         	\raisebox{-8pt}{\input draws/Sigmaw2.epc }
         	\hspace{-1.9mm}
	\end{array}
, \dots \}
\longleftrightarrow
\{\dots, \setlength{\unitlength}{0.03\standardunitlength}
	\begin{array}{c}  \hspace{-1.7mm}
         	\raisebox{-8pt}{\input draws/Sigmaw3.epc }
         	\hspace{-1.9mm}
	\end{array}
, \setlength{\unitlength}{0.03\standardunitlength}
	\begin{array}{c}  \hspace{-1.7mm}
         	\raisebox{-8pt}{\input draws/Sigmaw4.epc }
         	\hspace{-1.9mm}
	\end{array}
, \dots \}
$$
\end{definition}

It is easy to see that a wrapping relation of sliced crossed diagrams
implies that the corresponding links in $S^3\sminus\O$ are isotopic.
Indeed, 
$$
\psdraw{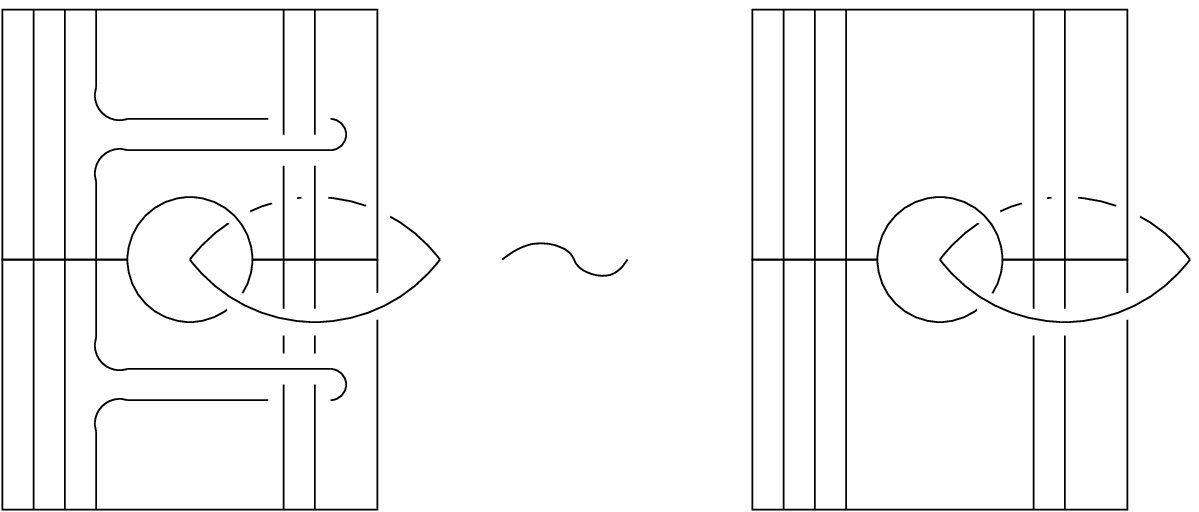}{3in}
$$
The above discussion implies that:

\begin{proposition}
\lbl{prop.tangles3}
There is a 1-1 correspondence: 
$$
\cD(D_g)/\la r, \b, \w \ra \longleftrightarrow \NO.
$$
\end{proposition}

Recall that Theorem \ref{thm.GK2} identifies the set of $F$-links with
a quotient of $\NO$. Since
we are interested to give a description of the set of $F$-links 
in terms of a quotient of $\cD(D_g)$, we need to introduce analogs of
the $\kappa$-relation on $\cD(D_g)$.

\begin{definition}
\lbl{def.handleslide}
The {\em handle-slide move} 
$\stackrel{\kappa}{\sim}$ is generated by the move:
\begin{eqnarray*}
\lefteqn{
\{\dots,T_{i-1},T_{i},\setlength{\unitlength}{0.03\standardunitlength}
	\begin{array}{c}  \hspace{-1.7mm}
         	\raisebox{-8pt}{\input draws/kirby1.epc }
         	\hspace{-1.9mm}
	\end{array}
,T_{i+2},\dots\} 
} & & \\
& \longleftrightarrow &
\{\dots,T_{i-1}',T_{i}',\setlength{\unitlength}{0.03\standardunitlength}
	\begin{array}{c}  \hspace{-1.7mm}
         	\raisebox{-8pt}{\input draws/kirby2.epc }
         	\hspace{-1.9mm}
	\end{array}
,
\setlength{\unitlength}{0.03\standardunitlength}
	\begin{array}{c}  \hspace{-1.7mm}
         	\raisebox{-8pt}{\input draws/kirby3.epc }
         	\hspace{-1.9mm}
	\end{array}
,T_{i+2}',\dots\}.
\end{eqnarray*}
where each tangle in $\{T_k'\}$ is obtained from the corresponding tangle in 
$\{T_k\}$ by taking a framed parallel of every component contributing to
the component $x_j,$. 
\end{definition}

We end this section with our final proposition:

\begin{proposition}
\lbl{prop.tangles}
There is a 1-1 correspondence 
$$
\cD(D_g)/ \la r, \b, \w, \kappa \ra \longleftrightarrow \text{$F$-links}.
$$
\end{proposition}

\begin{proof}
$F$-links are defined modulo inner automorphisms of the free group. One
need only observe that inner automorphisms follow from the basing relation.
\end{proof}

\section{An assortment of diagrams and their relations}
\lbl{sec.relations}

\subsection{Diagrams for the Kontsevich integral of a link}
\lbl{sub.diagramsK}

In this section we explain what are the spaces of diagrams in which our 
invariants will take values. 
Let us begin by recalling well-known spaces of diagrams that the Kontsevich
integral of a knot takes values.

The words ``diagram'' and ``graph'' will be used throughout the paper
in a synonymous way. All graphs will have unitrivalent vertices. The univalent
vertices are attached to the {\em skeleton} of the graphs, which will consist
of a disjoint union of oriented segments (indicated by $\tline_X$), circles
(indicated by $\tcircle_X$) or the empty set (indicated by $\star_X$). 
The skeleton is labeled by a set $X$ in 1-1 correspondence with the 
components of a link. The {\em Vassiliev degree} of a diagram is half the 
number of vertices.

\begin{definition}
\lbl{def.AK}
Let $\A(\tcircle_X, \tline_Y, \star_T, \ostar_S)$ denote the quotient of the 
completed graded $\BQ$-vector space spanned by diagrams with the prescribed 
skeleton,
modulo the $\AS$, $\IHX$ relations, the $\mathrm{STU}$ relation on $X$ and
$Y$ and the $S$-colored {\em infinitesimal basing relation} (the latter
was called $S$-colored link relation in \cite[Part II, Sec. 5.2]{A}) shown
by example for $S=\{x\}$:
$$
\setlength{\unitlength}{0.03\standardunitlength}
	\begin{array}{c}  \hspace{-1.7mm}
         	\raisebox{-8pt}{\input draws/LinkRel.epc }
         	\hspace{-1.9mm}
	\end{array}

$$
where the right hand side, by definition of the relation, equals to zero.
\end{definition}

When $X$, $Y$, $T$ or $S$ are the empty set, then they will be ommitted
from the discussion. In particular, $\A(\phi)$ denotes the completed vector
space spanned by trivalent diagrams, modulo the $\AS$ and $\IHX$ relations.

By its definition, the Kontsevich integral $Z(L)$ of a link $L$ takes values in
$\A(\tcircle_H)$, where $H$ is a set in 1-1 correspondence with the components
of $L$.

There are several useful vector space isomorphisms of the spaces of
diagrams which we now recall.

There is a {\em symmetrization map}
\begin{equation}
\lbl{eq.chisym}
\chi_X: \A(\star_X) \longto \A(\tline_X)
\end{equation}
which is the average of all ways of placing the symmetric legs of a diagram
on a line, whose inverse is denoted by 
\begin{equation}
\lbl{eq.sigmasym}
\s_X: \A(\tline_X) \longto \A(\star_X).
\end{equation}
Moreover, $\chi_X$ induces an isomorphism:
\begin{equation}
\lbl{eq.ostartcircle}
\A(\ostar_X) \longto \A(\tcircle_X)
\end{equation}
Thus, the Kontsevich integral of a link may be interpreted to take values
the quotient of $\A(\star_H)$ modulo the $H$-colored infinitesimal
basing relations.

Notice that $\A(\tline_X)$ is an algebra (with respect to the stacking
of diagrams that connects their skeleton), and that $\A(\star_X)$ is an
algebra with respect to the disjoint union multiplication. However, the
map $\chi_X$ is not an algebra map.

\subsection{Group-like elements and relations}
\lbl{sub.groupA}

In this short section we review the notion of {\em group-like elements}
in the spaces $\A(\tcircle_X)$, $\A(\tline_X)$, $\A(\star_X)$
and $\A(\ostar_X)$. All these spaces have natural Hopf algebra structures,
which are cocommutative and completed.

\begin{definition}
\lbl{def.grouplike}
An element $g$ in a completed Hopf algebra is {\em group-like} if $g=\exp(x)$
for a {\em primitive} element $x$.
\end{definition}

In case of the Hopf algebras of interest to us, the primitive elements
are the span of the {\em connected graphs}.

We we denote the group like elements of 
$$
\begin{cases}
\A(\tcircle_X) \\
\A(\tline_X) \\
\A(\star_X) \\
\A(\ostar_X) 
\end{cases}
\qquad
\text{by}\qquad
\begin{cases}
\GA(\tcircle_X) \\
\GA(\tline_X) \\
\GA(\star_X) \\
\GA(\ostar_X) 
\end{cases}
$$
The isomorphisms of Equations \eqref{eq.chisym}, \eqref{eq.sigmasym} and
\eqref{eq.ostartcircle} induce isomorphisms of the corresponding sets
of group-like elements.

Note that the Kontsevich integral of a link (or more generally, a tangle)
takes values in the set of group-like elements. 

Let us now present a group-like version of the infinitesimal basing
relation. Since this does not appear in the literature, and since it will
be the prototype for basing relations of group-like elements which include
beads, we will discuss it more extensively.

In \cite[part II, Sec.5.2]{A}, the following map
\begin{equation}
\lbl{eq.mxyz}
\vec{m}^{y,z}_x: \GA(\tline_{\{y,z\}}\cup \tline_X,\Lloc)\to
\GA(\tline_{\{x\}}\cup \tline_X,\Lloc)
\end{equation}
was introduced, that glues the end of the $y$-skeleton to the beginning 
of the $z$-skeleton and then relabeling the skeleton by $x$, 
for $x,y,z \not\in X$.

\begin{figure}[ht!]
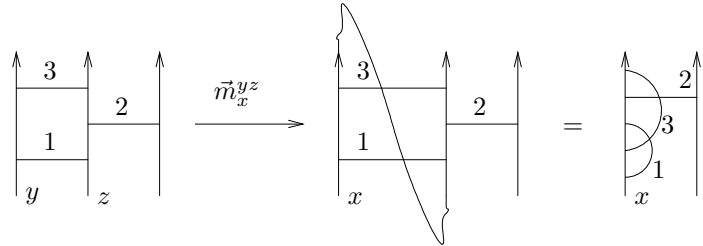
\small
$$\setlength{\unitlength}{0.05\standardunitlength}
	\begin{array}{c}  \hspace{-1.7mm}
         	\raisebox{-8pt}{\input draws/mxyz.epc }
         	\hspace{-1.9mm}
	\end{array}
$$
\caption{
  The map $\vec{m}^{yz}_x$ in degree $2$:
  Connect the strands labeled $x$ and $y$ in a diagram
  in $\A(\uparrow_x\uparrow_y\uparrow)$, to form a new ``long''
  strand labeled $z$, without touching all extra strands.
}
\end{figure}

\begin{definition}
\lbl{def.GAKbasing}
Given $s_1,s_2 \in$ $\GA(\tline_X)$, we say that  \newline
$\bullet$\qua $s_1 \stackrel{\bgp}\sim {s_2}$, 
if there exists an element $s\in \GA(\uparrow_{\{y\}}\uparrow_{\{z\}}
\uparrow_{X-\{x\}})$ with the property that
\begin{equation}
\lbl{eq.GAbasingb}
\vec{m}^{y,z}_x(s) = s_1 \,\,\,\, \mathrm{ and } \,\,
\vec{m}^{z,y}_x(s) = s_2.
\end{equation}
Using the isomorphism $\s_X$, we can define $\stackrel{\bgp}\sim$ on
$\GA(\star_X)$. We define
\begin{equation}
\lbl{eq.GAbasingc}
\GA(\ostar_X) = \GA(\star_X)/\la \bgp \ra .
\end{equation}
\end{definition}

The group-like basing relation can be formulated in terms of pushing
an exponential $e^h$ of hair to group-like diagrams. Since this was not noticed
in the work of \cite{A}, and since in our paper, pushing exponential
of hair is a useful operation, we will give the following reformulation
of the group-like basing relation.

\begin{definition}
Let 
$$
\con: \A(\star_{X\cup \pt X\cup Y}) \rightarrow \A(\star_Y)
$$ 
denote the {\em contraction map} 
defined as follows. For $s \in \A(\star_{X\cup \pt X \cup Y})$,
$\con_X(s)$ denotes the sum of all ways of  pairing {\em all} legs labeled 
from $\pt X$ with {\em all} legs labeled from $X$. 
The contraction map preserves group-like elements.
\end{definition}

\begin{definition}
\lbl{def.GAalternative}
For $s_1,s_2 \in \GA(\star_X)$, we say that
\newline
$\bullet$\qua $s_1 \stackrel{\bgpp}\sim s_2$ iff there exists $s \in
\GA(\star_{X \cup \{\pt h\}})$ such that 
$$
s_1=\conh(s) \quad \text{ and } \quad
s_2=\conh(s|_{X\to Xe^h})
$$
where $s|_{X\to Xe^h}$ is by definition the result of pushing
$e^h$ to each $X$-labeled leg of $s_{12}$.
\end{definition}

\begin{lemma}
\lbl{lem.b1G}
The equivalence relations $\bgpp$ and  $\bgp$ 
are equal on $\GA(\star_X)$.
\end{lemma}

\begin{proof}
Consider $s_1, s_2 \in \GA(\star_X)$ so that 
$\chi_X(s_1)\stackrel{\bgp}\sim \chi_X(s_2)$ and consider the corresponding
element $s \in \GA(\uparrow_{\{y\}}\uparrow_{\{z\}}
\uparrow_{X-\{x\}})$ satisfying 
Equation \eqref{eq.GAbasingb}. We start by presenting $s$ as
a result of a contraction. 
Inserting $\chi_{\{y,z\}}$ and $\s_{\{z,y\}}$ to $s$, we have that
$$
s=\psdraw{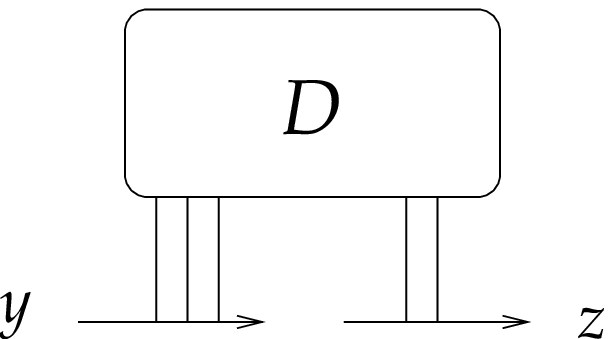}{1in}=
\con_{\{u,v\}} \left(\psdraw{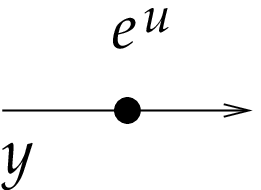}{0.5in} \quad 
\psdraw{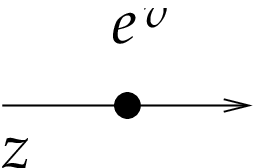}{0.5in} \quad
 \a \right)
$$
where
$$
\a=
(\s_{\{y,z\}}(s))|_{y \to \pt u, z \to \pt v}=
\psdraw{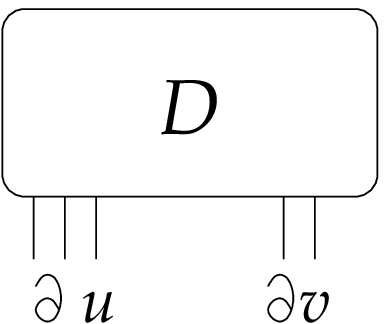}{0.8in}.
$$
Morevoer, we have:
\begin{eqnarray*}
\chi_X(s_1) &=& 
\con_{\{u,v\}} \left(\psdraw{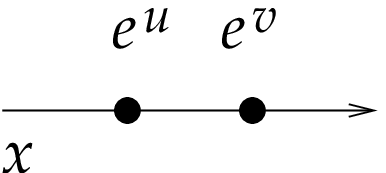}{1in} \,\,\, \a \right) \\
\chi_X(s_2) &=& 
\con_{\{u,v\}} \left(\psdraw{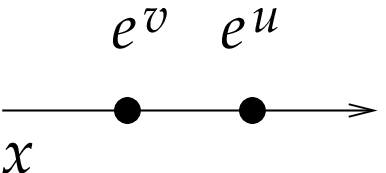}{1in} \,\,\, \a \right) 
\end{eqnarray*}
Now, we bring $e^v e^u$ to $e^u e^v$, at the cost of pushing
$e^v$ hair on the $u$-colored legs. Indeed, we have:
\begin{eqnarray*}
\psdraw{cyclicd5}{1in} &=&
\psdraw{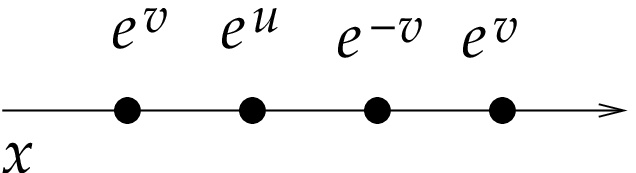}{1.7in} \\ &=&
\sum_{n=0}^\infty 
\psdraw{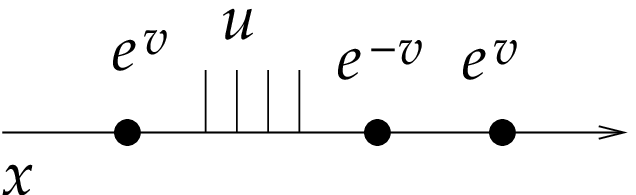}{1.7in} \\ &=&
\sum_{n=0}^\infty 
\psdraw{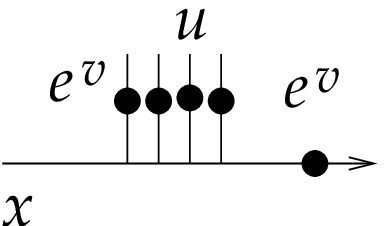}{1.0in}
\end{eqnarray*}
Remembering that we need to contract $(u,\pt u)$ legs and $(v, \pt v)$
legs, we can push the $e^v$-hair on $\a$. Let $\b$ denote the result of
pushing $e^v$ hair on each $\pt u$-colored leg of $\a$.
Then, we have:
\begin{eqnarray*}
\chi_X(s_1) &=& 
\con_{\{u,v\}} \left(\psdraw{cyclicd4}{1in} \,\,\, \a \right) \\
\chi_X(s_2) &=& 
\con_{\{u,v\}} \left(\psdraw{cyclicd4}{1in} \,\,\, \b \right) 
\end{eqnarray*}
Now, we may get $s_2$ from $\chi_X(s_2)$ by attaching $(u,v)$-rooted labeled
forests $T(u,v)$ whose root is colored by $x$; see \cite[part II, prop.5.4]{A}.
Since
$$
\con_{\{u\}} \psdraw{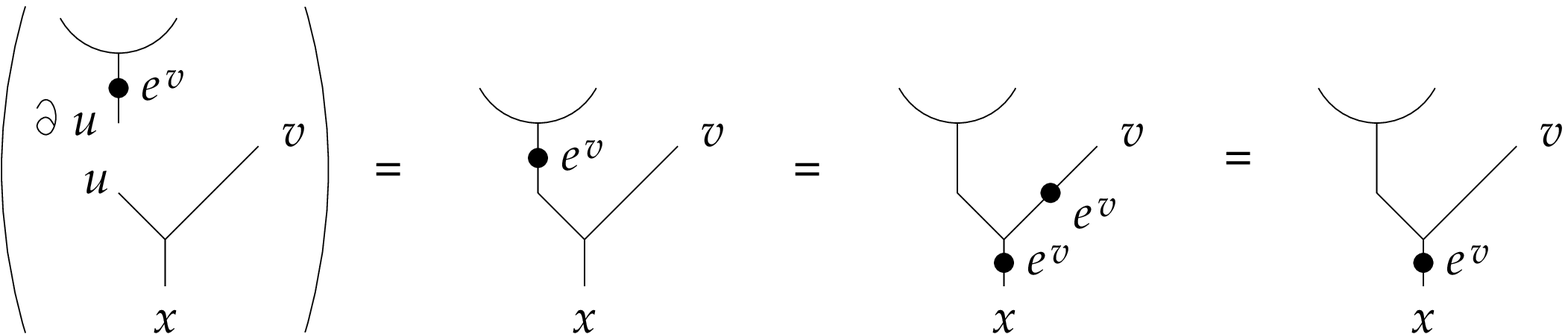}{4.7in}
$$
it follows that 
\begin{eqnarray*}
s_2 &=& \con_{\{u,v\}}(\a \,\, T(u,v)|_{x \to x e^v}) \\
s_1 &=& \con_{\{u,v\}}(\a \,\, T(u,v))
\end{eqnarray*}
Letting $\g=\con_{\{u\}}( \a \,\, T(u,v))$, it follows that
\begin{eqnarray*}
s_2 &=& \con_{\{v\}}(\g|_{X \to X e^v}) \\
s_1 &=& \con_{\{v\}}(\g)
\end{eqnarray*}
In other words, $s_1 \stackrel{\bgpp}\sim s_2$.

The converse follows by reversing the steps in the above
proof.
\end{proof}

\subsection{Diagrams for the rational form of the Kontsevich integral}
\lbl{sub.diagrams}

Let us now introduce diagrams with beads which will be useful in our paper.
The notation will generalize the notion of the $\A$-groups
introduced in \cite{GL1}. 

\begin{definition}
\lbl{def.edgelabel}
Consider a ring $R$ with involution and a 
distinguished group of units $U$. An admissible
labeling of a diagram $D$ with prescribed skeleton is a labeling of the
edges of $D$ and the edges of its skeleton so that:
\begin{itemize}
\item
the labelings on the $\tcircle_X$ and $\tline_Y$ lie in $U$, and satisfy 
the condition that the product of the labelings of the edges
along each component of the skeleton is $1$.
\item
the labelings on the rest of the edges of $D$ lie in $R$.
\end{itemize}
\end{definition}

Labels on edges or part of the skeleton will be called {\em beads}.

\begin{definition}
\lbl{def.AKL}
Consider a ring $R$ with involution and 
a distinguished group of units $U$, and 
(possibly empty sets) $X,Y,T$. Then,
$$
\A(\tcircle_X,\tline_Y,\star_T,R,U)
=\frac{{\mathcal D}(\tcircle_X,\tline_Y \star_T,R,U)}{\lp
\AS, \IHX, \mathrm{STU}, \mathrm{Multilinear}, \text{Vertex Invariance} \rp }
$$
where: 
\begin{itemize}
\item
${\mathcal D}(\tcircle_X,\tline_Y \star_T,R,U)$ is the completed graded
vector space over $\BQ$ of diagrams of prescribed skeleton with oriented
edges, and with admissible labeling.
\item
The {\em degree} of a diagram is the number of trivalent vertices.
\item
$\AS, \IHX$ and $\mathrm{Multilinear}$ are the relations shown in
Figure \ref{relations2} and Vertex Invariance is the relation shown in
Figure \ref{relations3}. Note that all relations are homogeneous, thus
the quotient is a completed graded vector space.
\end{itemize}
\end{definition}

Some remarks on the notation. 
Empty sets will be omitted from the notation, and so will $U$, the selected
group of units of $R$. For example, $\A(\star_Y, R)$, $\A(R)$ and $\A(\phi)$
stands for $\A(\tline_\phi,\star_Y, \ostar_\phi,R,U)$, 
$\A(\tline_\phi,\star_\phi, \ostar_\phi,R,U)$ and $
\A(\tline_\phi,\star_\phi, \ostar_\phi,\BZ,1)$ respectively. 
Univalent vertices of diagrams will often be called {\em legs}. 
Special diagrams, called 
{\em struts}, labeled by $a,c$ with bead $b$ are drawn as follows
$$
\strutb{a}{c}{b}.
$$
oriented from bottom to top.
Multiplication of diagrams $D_1$ and $D_2$, unless otherwise 
mentioned, means their disjoint union and will be denoted by $D_1 \, D_2$ and 
occasionally by $D_1 \sqcup D_2$.

We will be interested only in the rings 
\begin{itemize}
\item
$\La=\BZ[F]$ (the group ring of the 
free group $F$ on some fixed set $T=\{t_1,\dots,t_g\}$ of generators), 
\item
its completion $\Lhat$ (with respect to the powers of the augmentation ideal)
\item
and its {\em Cohn localization}
$\Lloc$ (ie, the localization with respect to the set of matrices over $\La$
that are invertible over $\BZ$).
\end{itemize}
For all three rings $\La,\Lhat$ and $\Lloc$ the selected group of units is 
$F$. All three are rings with involution induced by 
$\bar{g}=g^{-1}$ for $g \in F$, with augmentation over $\BQ$, and with 
a commutative diagram 
$$
\divide\dgARROWLENGTH by2
\begin{diagram}
\node{\La} 
\arrow[2]{e} \arrow{se} 
\node[2]{\Lloc}
\arrow{sw}  \\
\node[2]{\Lhat}
\end{diagram}
$$
where $\La\to\Lhat$ is given by the exponential version of the Magnus expansion
$t_i\to e^{h_i}$. $\Lhat$ can be indentified with the ring of noncommuting 
variables $\{h_1,\dots,h_g\}$ and $\Lloc$ with the ring of rational functions 
in noncommuting variables $\{t_1,\dots,t_g\}$, \cite{FV,BR}. In all cases,
the distinguished group of units is $F$.

\begin{figure}[ht!]
$$ 
\psdraw{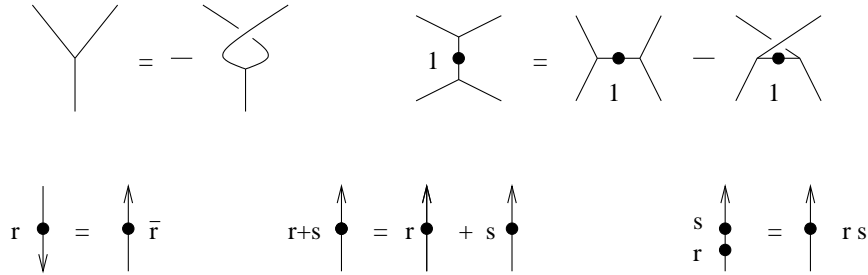}{4.5in} 
$$
\caption{The $\AS$, $\IHX$ (for arbitrary orientations of the edges),
Orientation Reversal and Linearity relations.
Here $r\to \bar{r}$ is the involution of $R$, $r,s \in R$.}\lbl{relations2}
\end{figure}

\begin{figure}[ht!]
$$ 
\psdraw{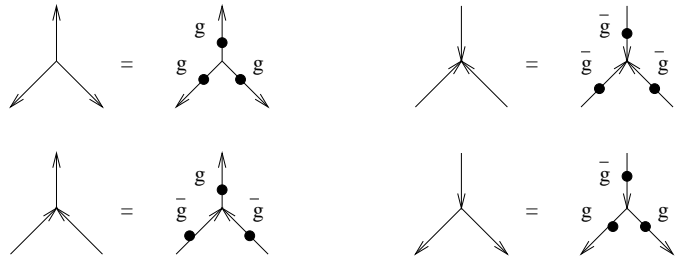}{3.5in} 
$$
\caption{The Vertex Invariance relation that pushes a unit $g \in U$ past a 
trivalent 
vertex.}\lbl{relations3}
\end{figure}

The following is a companion to Remark \ref{rem.lara}:
\begin{remark}
\lbl{rem.lprp}
Throughout the paper, given a subspace $W$ of a vector space $V$,
we will denote by $V/ W $ the quotient of $V$ modulo $W$. 
$W$ will be defined by one of the following ways:
\begin{itemize}
\item
By a subset $S$ of $V$, in which case $W=(S)$ is the subspace spanned by 
$S$.
\item
By a {\em linear action} of a group $G$ on $V$, in which case 
$S=\{ gv-v\}$ defines $W$ and $V/ W $ coincides 
with the space of {\em coinvariants} of $G$ on $V$.
\end{itemize}
\end{remark}

Since labels of the edges of the skeleton are in $U$, and 
the product of the labels around each skeleton component is $1$, the
the Vertex Invariance Relation implies the following analogue
of Equations \eqref{eq.chisym}, \eqref{eq.sigmasym}:

\begin{lemma}
\lbl{lem.chisigma}
For every $X$ and $R=\La,\Lloc$ or $\Lhat$,  there are inverse maps
$$
\chi_X: \A(\star_X,R) \longto \A(\tline_X,R)
\quad \text{and} \quad
\sigma_X: \A(\tline_X,R) \longto \A(\star_X,R). 
$$
\end{lemma}

Our next task is to introduce analogs for the basing and the wrapping 
relations for diagrams with beads. Before we do that, let us discuss
in detail the beads that we will be considering.

\subsection{Noncommutative localization}
\lbl{sub.Cohn}

In this self-contained section we review several facts about the Cohn 
localization that are used in the text. 
Consider the free group $F_g$ on generators $\{t_1,\dots,t_g\}$. 
The Cohn localization $\Lloc$ of its group-ring $\La=\BZ[F_g]$ is 
characterized by a universal property, 
$$ 
\begin{diagram}
\node{\La}       
\arrow{e}\arrow{se,r}{\a}
\node{\Lloc}
\arrow{s,r,..}{\bar\a}    \\
\node[2]{R}      
\end{diagram}
$$
namely that for every $\Sigma$-{\em inverting} ring homomorphism $\a:
\La\to R$ there exists a unique ring homomorphism $\bar{\a}:\Lloc\to R$
that makes the above diagram commute, \cite{C}. Recall that a homomorphism
$\a$ is $\S$-inverting if $\a M$ is invertible over $R$ for every matrix
$M$ over $\La$ which is invertible over $\BZ$.

Farber and Vogel identified $\Lloc$ with the ring of {\em rational functions in
noncommuting variables}, \cite{FV,BR}. An example of such a rational function 
is
$$
(3-t_1^2t_2^{-1}(t_3+1)t_1^{-1})^{-1} t_2 -5.
$$ 
There is a close relation between rational functions in noncommuting variables
and {\em finite state automata}, sedcribed in \cite{BR,FV}.
Farber-Vogel showed that every
element $s \in \Lloc$ can be represented as a solution to a system of equations
$M \vec{x}=\vec{b}$ where $M$ is a matrix over $\La$, invertible over $\BZ$,
and $\vec{b}$ is a column vector over $\La$, \cite[Proposition 2.1]{FV}. 
More precisely, we have that
$$
s=(1,0,\dots,0)M^{-1} \vec{b},
$$
which we will call a {\em matrix presentation} of $s \in \Lloc$.

In the text, we often use substitutions of the form $\varphi_{t_i\to t_ie^h}$,
or $\varphi_{t_i\to e^h t_i}$ or $\varphi_{t_i \to e^{-h} t_i e^h}$, 
where $h$ does not commute with $F_g$. In order to make sense of these
substitutions, we need to enlarge our ring $\Lloc$. This can be achieved
in the following way.
$F_g$ can be included in the free group $F$
with one additional generator $t$. Let $\cL$ denote the 
group-ring of $F$, $\cL_{\loc}$ denote its Cohn localization and 
$\Lnew$ denote the completion of $\cL_{\loc}$ with respect to the
the ideal generated by $t-1$. There is an identification of $\Lnew$ with 
the following ring $\R$ whose elements consist of formal sums of the form 
$\sum_{\sfm,\sff}  \prod_{i=1}^\infty h^{m_i} f_i$
where $f_i \in \Lloc$ and 
\begin{itemize}
\item
$\sfm=(m_1,m_2,\dots):\BN\to\BN$ is eventually $0$ and 
$\sff=(f_1,f_2,\dots):\BN\to\Lloc$ is eventually $1$.
\item
In the above sum, for each fixed $k$, there are finitely many sequences
$\sfm$ with $|\sfm|=\sum_i m_i \leq k$.
\end{itemize}
It is easy to see that $\R$ is a ring with involution (in fact a subring
of the completion of $\cL$ with respect to the augmentation ideal) such that 
$t_i$ and $e^h$ are units, where $h=\log t$.
In particular, $\Lnew$ is an $h$-graded ring. Let $\degh^n:\Lnew\to\Lnew^n$
denote the projection on the $h$-degree $n$ part of $\Lnew$.

\begin{definition}
\lbl{def.vphi}
There is a map $\varphi_i: \La \to \Lnew$ given by substituting $t_i$ for
$e^{-h}t_ie^h$ for $1 \leq i \leq g$. It extends to a map 
$\varphi_i: \Lloc\to\Lnew$. 
\end{definition}
Indeed, using the defining property of the Cohn localization, it suffices
to show that $\varphi_i$ is $\S$-{\em inverting}, \cite{C,FV}. In other 
words, we need to show that if a  matrix $W$ over $\La$ is invertible 
over $\BZ$, then 
$\varphi_i(W)$ is invertible over $\Lnew$. Since $\Lnew$ is a completion
of $\cL_{\loc}$, it suffices to show that $\degh^0(\varphi_i(W))$
is invertible over $\cL_{\loc}$. However, $\degh^0(\varphi_i(W))=W$,
invertible over $\Lloc$ and thus also over $\cL_{\loc}$.
The result follows.

\begin{definition}
\lbl{def.deg1}
Let 
$$
\eta_i:\Lloc\to\Lnew^1
$$ 
denote the $h$-degree $1$ part of $\varphi_i$.
\end{definition}

The following lemma gives an axiomatic definition of $\eta_i$ by properties
analogous to a {\em derivation}. Compare also with the derivations of Fox
differential calculus, \cite{FV}.

\begin{lemma}
\lbl{lem.winf1}
$\mathrm{(a)}$\qua $\eta_i$ is characterized by the following properties:
\begin{eqnarray*}
\eta_i(t_j) & = & \delta_{ij}[t_i,h]=\delta_{ij}(t_ih-ht_i) \\
\eta_i(ab) & = & \eta_i(a)b+a \eta_i(b) \\
\eta_i(a+b) & = & \eta_i(a) + \eta_i(b).
\end{eqnarray*}
$\mathrm{(b)}$\qua 
Moreover, $\eta_i$ can be extended to matrices with entries in $\Lloc$.
In that case, $\eta_i$ satisfies 
$$
\eta_i(AB) = \eta_i(A) B + A \eta_i(B) 
$$
for matrices $A$ and $B$ that can be multiplied, as well as
$$
\eta_i(A^{-1})=-A^{-1} \eta_i(A) A^{-1}
$$
for invertible matrices $A$.
\end{lemma}

\begin{proof}
For (a), 
it is easy to see that $\eta_i$ satisfies the above properties, by looking
at the $h$-degree $0$ and $1$ part in $\varphi_i(a)$. This determines $\eta_i$
uniquely on $\La$. By the universal property, we know that $\eta_i$ extends
over $\Lloc$; a priori the extension need not be unique. Since however,
$\eta_i$ extends uniquely to matrices over $\La$, and since every element
of $\Lloc$ has a matrix presentation (see the discussion in the beginning
of the section), it follows that $\eta_i$ is characterized by the 
derivation properties.

Part (b) follows easily from part (a).
\end{proof}

\subsection{The infinitesimal wrapping relations}
\lbl{sub.winfinite}

We now have all the required preliminaries to introduce an
infinitesimal \wrapping\  relation on 
\begin{equation}
\lbl{eq.A0}
\A^0(\star_X,\Lloc):=\B(\Lloc\to\BZ)\times\A(\star_X,\Lloc).
\end{equation}
Given a vector $v \in \A(\star_{X},\Lloc)$ and a  
diagram $D \in \A(\star_{X},\Lloc)$ with a {\em single} special leg 
colored by $\pt h$ and $M \in \Herm(\Lloc\to\BZ)$, the {\em infinitesimal 
\wrapping\ relation } $\winf_i$ at the $i$th place in
$\Herm(\Lloc\to\BZ) \times \A(\star_{X},\Lloc)$ is generated by the move
$$
(M,v) \stackrel{\winf_i}\sim \left( M, v+
\conh \left( \eta_i(D) -\frac{1}{2} D \sqcup \psdraw{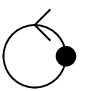}{0.2in}
\tr(M^{-1} \eta_i(M)) \right)
\right).
$$
It is easy to see that the above definition depends only on the image
of $M$ in $\B(\Lloc\to\BZ)$. We define 
$$
\A^0(\star_X,\Lloc)/\la \winf \ra = \A^0(\star_X,\Lloc)/\la \winf_1, \dots,
\winf_g \ra .
$$
We give two examples of infinitesimal wrapping relations:

\begin{example}
\lbl{ex.wrap1}
An important special case of the infinitesimal wrapping relation is
when $M$ is the empty matrix and the beads of a diagram $D$ 
lie in $\La=\BZ[F]$.
Without loss of generality, we may assume that the beads of $D$ are 
monomials $t_j$. If {\em all} the edges labeled by $t_i$ are shown, 
and we define
$$
\psdraw{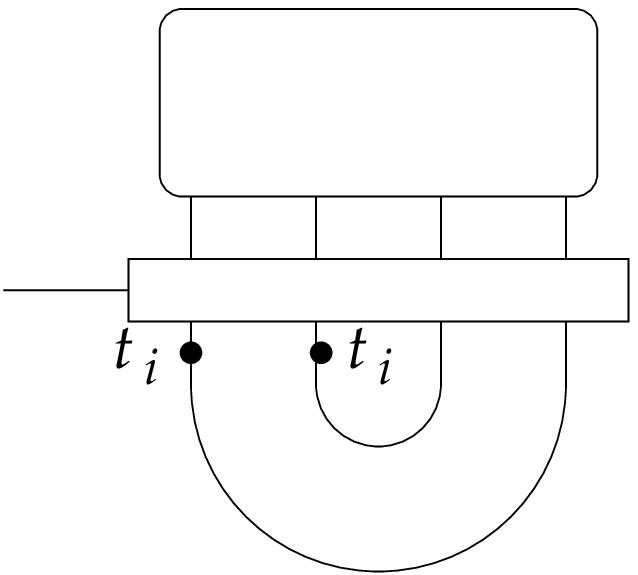}{0.9in} =\psdraw{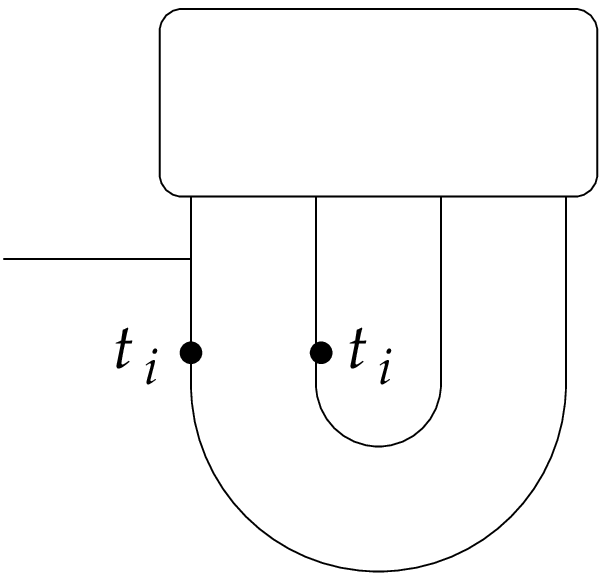}{0.9in}+
\psdraw{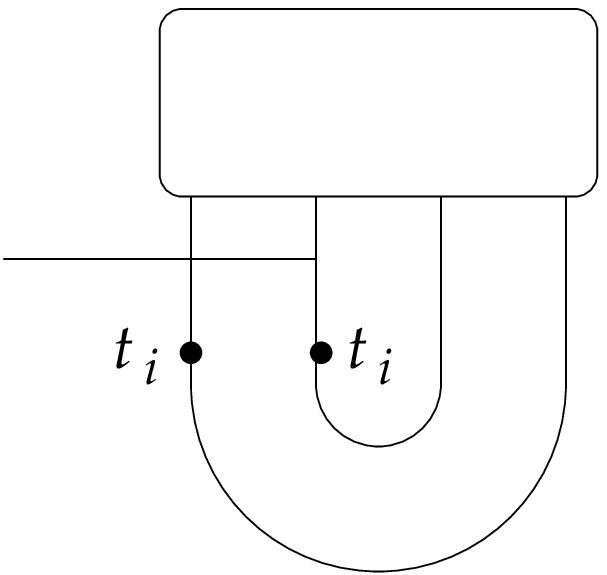}{0.9in}-\psdraw{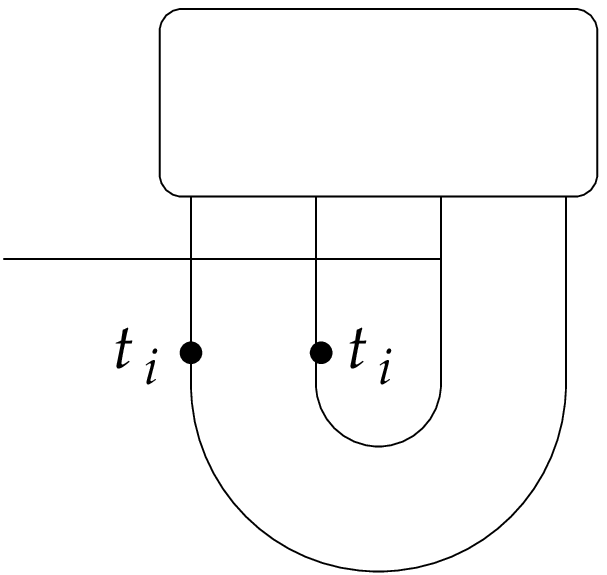}{0.9in}-\psdraw{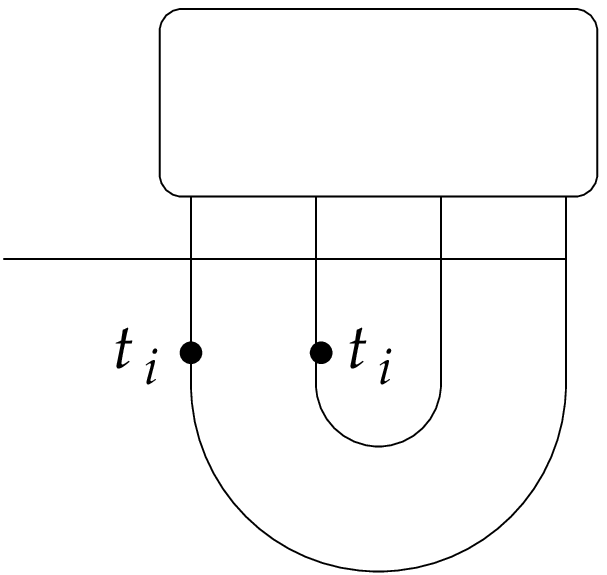}{0.9in}
$$
then we have:
$$
\psdraw{winf0}{0.9in}=0.
$$
\end{example}

Note that in the $\winf$ relation, we allow for $D$ to be a diagram of a circle
(with a bead and with a special $\pt h$ leg). In addition, we allow the
$\pt h$ leg to have a bead in $\Lloc$.

The next example illustrates the infinitesimal wrapping relation in its full
complication.

\begin{example}
\lbl{ex.wrap2}
Consider the matrix
$$
M=
\left[
\begin{array}{ll}
1 & t_3-t_2 t_1^{-1} \\ t_3^{-1}-t_1 t_2^{-1} & 1
\end{array}
\right]
$$
which is Hermitian over $\La$, and invertible over $\BZ$.
Part (b) of Lemma \ref{lem.winf1} implies that
\begin{eqnarray*}
\eta_1(M) &=& 
\mat {0}{-t_2(t_1^{-1}h -ht_1^{-1})}{-(t_1 h - h t_1) t_2^{-1}}{0} \\
&=& \mat{0}{0}{h t_1 t_2^{-1}}{0} +
\mat{0}{-t_2 t_1^{-1} h}{0}{0}
+ \mat{0}{t_2 h t_1^{-1}}{-t_1 h t_2^{-1}}{0}
\end{eqnarray*}
and
\begin{eqnarray*}
\eta_1(M^{-1}) & = &
M^{-1} (t_2 t_1^{-1}, 0)' h (0,1) M^{-1} +
M^{-1} (-t_2,  0)' h (0,t_1^{-1}) M^{-1} \\
&  + &
M^{-1} (0, t_1)' h (t_2^{-1},0) M^{-1} +
M^{-1} (0, -1)' h (t_1 t_2^{-1},0) M^{-1}.
\end{eqnarray*}
where $A'$ denotes the transpose of $A$. 
Further, observe that
\begin{eqnarray*}
-1/2 \tr(M^{-1} \eta_1(M)) & = &
(1/2,0) M^{-1} (0,t_1)' h t_2^{-1} +
(1/2,0) M^{-1} (0,-1)' h t_1 t_2^{-1} \\
& + &
(0,1/2) M^{-1} (t_2 t_1^{-1},0)' h +
(0,1/2) M^{-1} (-t_2,0)' h t_1^{-1}.
\end{eqnarray*}
Now, consider the infinitesimal wrapping relation arising from the diagram

{\small$$ D=\setlength{\unitlength}{0.03\standardunitlength}
	\begin{array}{c}  \hspace{-1.7mm}
         	\raisebox{-8pt}{\input draws/wrapex1.epc }
         	\hspace{-1.9mm}
	\end{array}
 $$}%
That is, we have $(M,v) \stackrel{\winf_1}\sim (M,v+x)$, where $x=x_1+x_2$,
and $x_1$ is given
by gluing the $(h,\pt h)$ legs in each of the following diagrams, 
{\small$$ 
\hspace{-2cm} \setlength{\unitlength}{0.03\standardunitlength}
	\begin{array}{c}  \hspace{-1.7mm}
         	\raisebox{-8pt}{\input draws/wrapex2.epc }
         	\hspace{-1.9mm}
	\end{array}
 \hspace{3cm}, \setlength{\unitlength}{0.03\standardunitlength}
	\begin{array}{c}  \hspace{-1.7mm}
         	\raisebox{-8pt}{\input draws/wrapex3.epc }
         	\hspace{-1.9mm}
	\end{array}

$$
$$ 
\hspace{-2cm} \setlength{\unitlength}{0.03\standardunitlength}
	\begin{array}{c}  \hspace{-1.7mm}
         	\raisebox{-8pt}{\input draws/wrapex4.epc }
         	\hspace{-1.9mm}
	\end{array}
 \hspace{3cm}, \setlength{\unitlength}{0.03\standardunitlength}
	\begin{array}{c}  \hspace{-1.7mm}
         	\raisebox{-8pt}{\input draws/wrapex5.epc }
         	\hspace{-1.9mm}
	\end{array}

$$}%
as well as these diagrams
{\small$$ 
\hspace{-3cm}
\setlength{\unitlength}{0.03\standardunitlength}
	\begin{array}{c}  \hspace{-1.7mm}
         	\raisebox{-8pt}{\input draws/wrapex6.epc }
         	\hspace{-1.9mm}
	\end{array}
 \hspace{4cm}, \setlength{\unitlength}{0.03\standardunitlength}
	\begin{array}{c}  \hspace{-1.7mm}
         	\raisebox{-8pt}{\input draws/wrapex7.epc }
         	\hspace{-1.9mm}
	\end{array}

$$}%
and summing the result up, and $x_2$ is given by gluing the $(h,\pt h)$ legs 
of the disjoint union of $D$ with the following sum:
{\small$$ 
\hspace{-3cm} \setlength{\unitlength}{0.03\standardunitlength}
	\begin{array}{c}  \hspace{-1.7mm}
         	\raisebox{-8pt}{\input draws/wrapex8.epc }
         	\hspace{-1.9mm}
	\end{array}
 \hspace{3cm}+ \setlength{\unitlength}{0.03\standardunitlength}
	\begin{array}{c}  \hspace{-1.7mm}
         	\raisebox{-8pt}{\input draws/wrapex8.epc }
         	\hspace{-1.9mm}
	\end{array}

$$
$$ 
+\setlength{\unitlength}{0.03\standardunitlength}
	\begin{array}{c}  \hspace{-1.7mm}
         	\raisebox{-8pt}{\input draws/wrapex10.epc }
         	\hspace{-1.9mm}
	\end{array}
 \hspace{3cm}+ \setlength{\unitlength}{0.03\standardunitlength}
	\begin{array}{c}  \hspace{-1.7mm}
         	\raisebox{-8pt}{\input draws/wrapex11.epc }
         	\hspace{-1.9mm}
	\end{array}

\hspace{3cm}.$$}%
Notice that $x$ is a sum of diagrams of degree $2$ of the form $\Th$ and 
$\eyes$. It is an interesting problem to understand the infinitesimal
wrapping relation in degree $2$.
\end{example}

\subsection{Group-like relations}
\lbl{sub.gp}

Our invariants (to be constructed) will take place in quotients of
group-like elements, modulo appropriate relations. In this section
we discuss these relations.

In analogy with Section \ref{sub.groupA}, the subscript $^{\mathrm{gp}}$
indicates the set of group-like elements of the appropriate Hopf algebra.

We now define an analog of the basing relations for group-like elements. 
Recall the map $\vec{m}^{yz}_x$ of Equation \eqref{eq.mxyz}.

\begin{definition}
\lbl{def.GAbasing}
Given $s_1,s_2 \in$ $\GAX$, we say that  \newline
$\bullet$\qua $s_1 \stackrel{\bgp_1}\sim {s_2}$, 
if there exists an element $s\in\ \GA(\uparrow_{\{y\}}\uparrow_{\{z\}}
\uparrow_{X-\{x\}},\Lloc)$ with the property that
\begin{equation}
\lbl{eq.b1}
\s_{\{x\}}(\vec{m}^{y,z}_x(s)) = s_1 \,\,\,\, \mathrm{ and } \,\,
\s_{\{x\}}(\vec{m}^{z,y}_x(s)) = s_2.
\end{equation}
$\bullet$\qua $s_1 \stackrel{\bgp_2}\sim {s_2}$ iff $s_2$ is obtained from $s_1$ 
by pushing an element of $F$ on each of the $X$-labeled legs of $s_1$. 
For example, if all $x_1$ legs are shown below, we have:
$$
\psdraw{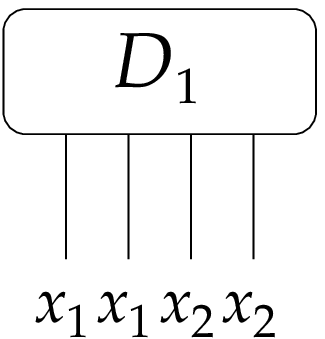}{0.5in} +\psdraw{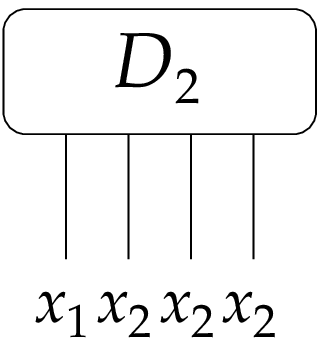}{0.5in} \quad
\stackrel{\binf_2}{\sim}  \quad
\psdraw{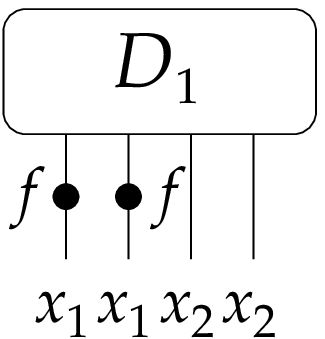}{0.5in} + \psdraw{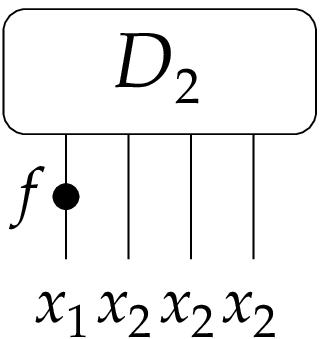}{0.5in}
$$
We define
\begin{equation}
\lbl{eq.GAbasing}
\GA(\ostar_X,\Lloc) = \GA(\star_X,\Lloc)/\la \bgp_1,\bgp_2 \ra .
\end{equation}
\end{definition}

\begin{remark}
\lbl{rem.glue}
When we push an element of the free group $f$ in the legs of a diagram,
depending on the orientation of the leg, we add $f$ or $\bar f$ to the 
leg. See also Figure \ref{relations3} for a convention for all choices
of orientations. Similarly, when we glue two legs of diagrams together,
the label on the bead depends on the orientations of the legs.
For example: 
$$
\psdraw{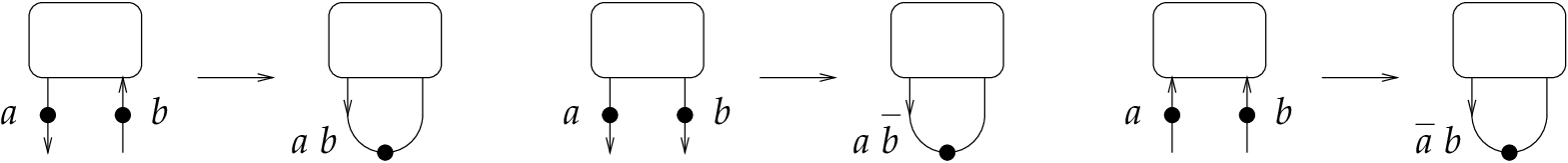}{4.2in}
$$
\end{remark}

We now define a useful reformulation of the $\bgp_1$ relation for diagrams in\break 
$\A(\star_X,\Lloc)$, which makes obvious that:

\begin{itemize}
\item
The $X$-flavored group-like basing relations push a bead, or an exponential of 
hair on all $X$-colored legs.
\end{itemize}

\begin{definition}
For $s_1,s_2 \in \GA(\star_X,\Lloc)$, we say that
\newline
$\bullet$\qua $s_1 \stackrel{\bgpp_1}\sim s_2$ iff there exists $s \in
\GA(\star_{X \cup \{\pt h\}},\Lloc)$ such that 
$$s_1=\conh(s) \quad \text{ and } \quad
s_2=\conh(s|_{X\to Xe^h})
$$
where $s|_{X\to Xe^h}$ is by definition the result of pushing
$e^h$ to each $X$-labeled leg of $s$.
\end{definition}

\begin{lemma}
\lbl{lem.b1}
The equivalence relations $\bgpp_1$ and  $\bgp_1$ are equal on\newline 
$\GA(\star_X,\!\Lloc)$\!.
\end{lemma}

\begin{proof}
It is identical to the proof of Lemma \ref{lem.b1G} and is omitted.
\end{proof}

Now, we define a group-like wrapping relation on group-like elements.
First, we will enlarge the set of group-like elements by a quotient
of a set of Hermitian matrices.
This enlargement might sound artificial, however there are several good 
reasons for doing so.
For example, the null move on the set of $F$-links and its associated 
filtration predicts that in degree $0$, the universal invariant is given by 
the S-equivalence class (of $F$-links), or equivalently, by the set of
Blanchfield pairings of $F$-links. $\GAz(\Lloc)$, defined below, maps
onto $\B(\Lloc\to\BZ)$ which in turn maps onto the set of Blanchfield pairings
of $F$-links.

In addition, the $\A$-groups, studied in relation to Homology Surgery in 
\cite{GL1},
are secondary obstructions to the vanishing of a Witt-type obstruction
which lies in $\B(\BZ[\pi]\to\BZ)$. In the case of $F$-links, this
motivates the fact that the set $\B(\Lloc\to\BZ)$ should be taken together
with $\A(\Lloc)$. With these motivations, 

\begin{definition}
\lbl{def.GAz}
Let us define
$$
\GAz(\star_X,\Lloc)=\B(\Lloc\to\BZ)\times \GA(\star_X, \Lloc) .
$$
\end{definition}

\begin{definition}
\lbl{def.wrappingD}
Two pairs $(M,s_m) \in \Herm(\Lloc\to\BZ) \times \GA(\star_X,\Lloc)$ $m=1,2$ 
are 
related by a {\em \wrapping\ move} at the $j$th gluing site (for some $j$ such 
that $1\leq j \leq g$) if 
there exists an $s \in \GA(\star_{X \cup \{\pt h\}},\Lloc)$ such that:
\begin{eqnarray*}
s_m & = & \conh \left( 
\varphi_m(s)\sqcup \left( 
\exp_\sqcup \left(-\frac{1}{2} \chiz(M^{-1}
\varphi_m M) \right)\right)\right)
\end{eqnarray*}
for $m=1,2$ where $\varphi_1$ is the identity and $\varphi_2$ is the 
substitution $\varphi_{t_j\to e^{-h} t_j e^h}$, and
\begin{equation}
\lbl{eq.chiz}
\chiz(A)=\psdraw{circleb}{0.2in} \,\,\tr\log(A) \in \A(\star_h,\Lloc).
\end{equation}
The \wrapping\ move generates an equivalence relation (the so-called
{\em group-like \wrapping\ relation}) on 
$\GAz(\star_X,\Lloc)$.
\end{definition}

Our first task is to make sure that the above definition makes sense.
Using the $h$-graded ring with involution $\Lnew$ from Section \ref{sub.Cohn}, 
it follows that $\varphi_m(s_{12}) \in \GA(\star_X,\Lnew)$.
Furthermore, $\chiz(M^{-1}\varphi_m M)$ can be thought of either as a circle
with a bead in $\Lnew/(\mathrm{cyclic})$ that has vanishing $h$-degree $0$ 
part, or as a sum
of wheels with $h$-colored legs and beads in $\Lloc$. Thus, the contraction 
$\conh$ in the 
above definition results in an element of $\GA(\star_X,\Lloc)$. It is easy to
see that $\chiz$ satisfies the following properties (compare with
\cite[Proposition 2.3]{GK2}

\begin{proposition}{\rm\cite{GK2}}\qua
\lbl{prop.w} 
{\rm (a)}\qua For Hermitian matrices $A,B$ over $\Lloc$, nonsingular over $\BZ$, we
have in $\Lhat/(\mathrm{cyclic})$ that
$$\chiz(AB)=\chiz(A)+\chiz(B) 
\,\,\,\, \text{ and } \,\,\,\, \chiz(A\oplus B)=\chiz(A)
+\chiz(B).$$
{\rm (b)}\qua $\chiz$ descends to a $\Lhat/(\mathrm{cyclic})$-valued invariant
of the set $\B(\Lloc\to\BZ)$. \newline
{\rm (c)}\qua For $A$ as above, $\chiz(A)=\chiz(A \, \e A^{-1})$ where 
$\e:\La\to\BZ$.
\end{proposition}

From this, it follows that $\chiz$ depends only on the image of $M$ in 
$\B(\Lloc\to\BZ)$ and the above definition makes sense in 
$\GAz(\star_X,\Lloc)$.

\begin{remark}
\lbl{rem.phichoice}
In the above definition \ref{def.wrapping} of the wrapping relation, we could 
have chosen equivalently that $\varphi_1$ and $\varphi_2$ denote the 
substitutions $\varphi_{t_j\rightarrow t_je^h}$ and 
$\varphi_{t_j\rightarrow e^ht_j}$.
\end{remark}

\begin{remark}
\lbl{rem.pth}
In the above definition, the $\pt h$ colored legs are allowed to have beads
in $\Lloc$. Also, we allow the case of the empty matrix $M$, with the
understanding that the $\chiz$ term is absent.
\end{remark}

\begin{remark}
\lbl{rem.specialw}
An important special case of the group-like wrapping relation is when
$M$ is the empty matrix, and the beads lie in $\La$. Without loss of 
generality, we may assume that the beads are monomials in $t_i$. 
In that case, if $s_i \in \GA(\star_X, \La)$ for $i=1,2$, 
and $e^h$ appears before (resp. after) the beads $t_i$, then $s_1$ is
wrapping move equivalent to $s_2$. For example, if all occurrences of $t_j$
are shown below, then we have:
$$
\begin{array}{l}
\stbup{}{}{e^h} \vspace{-0.4cm} \\ 
\stbup{}{}{t_j}
\end{array}
\begin{array}{l}
\stbdown{}{}{e^{-h}} \vspace{-0.4cm} \\ 
\stbdown{}{}{t_j^{-1}}
\end{array}
\begin{array}{l}
\stbdown{}{}{e^{-h}} \vspace{-0.4cm} \\ 
\stbdown{}{}{t_j^{-1}}
\end{array}
\begin{array}{l}
\stbup{}{}{e^h} \vspace{-0.4cm} \\ 
\stbup{}{}{t_j}
\end{array}
\stackrel{\wgp_j}{=}
\begin{array}{l}
\stbup{}{}{t_j} \vspace{-0.4cm} \\ 
\stbup{}{}{e^h}
\end{array}
\begin{array}{l}
\stbdown{}{}{t_j^{-1}} \vspace{-0.4cm} \\ 
\stbdown{}{}{e^{-h}}
\end{array}
\begin{array}{l}
\stbdown{}{}{t_j^{-1}} \vspace{-0.4cm} \\ 
\stbdown{}{}{e^{-h}}
\end{array}
\begin{array}{l}
\stbup{}{}{t_j} \vspace{-0.4cm} \\ 
\stbup{}{}{e^h}
\end{array} .
$$
\end{remark}

\begin{lemma}
\lbl{lem.wg=1}
For $g=1$, the wrapping relation $\wgp$ on $\GAz(\star_X,\Lloc)$ is trivial.
\end{lemma}

\proof
It follows from a swapping argument analogous to the isotopy of the following
figure 
$$ 
\psdraw{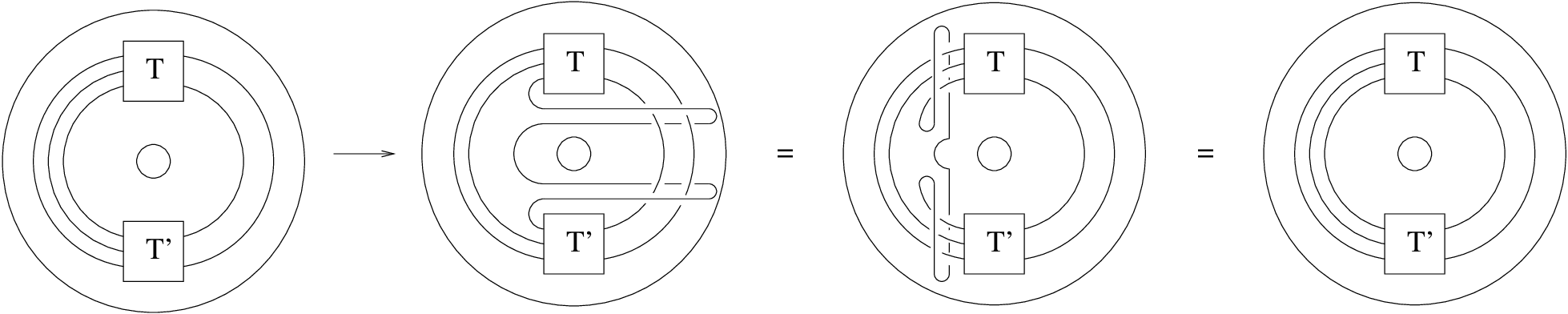}{5in}\eqno{\qed} 
$$

We end this section with a comparison between $\A^0$ and $\GAz$.
There is an obvious inclusion 
$$
\GAz(\star_X,\Lloc) \longrightarrow \A^0(\star_X,\Lloc). 
$$
The next proposition states that 

\begin{proposition}
\lbl{prop.winf3}
The inclusion
$\GAz(\star_X,\Lloc) \longrightarrow \A^0(\star_X,\Lloc)$
maps the $\wgp$-relations to $\winf$-relations, and thus induces a map
$$
\GAz(\star_X,\Lloc)/\la \wgp\ra \longrightarrow \A^0(\star_X,\Lloc)/\la \winf
\ra.
$$
\end{proposition}

\begin{proof}
See Appendix \ref{sec.wversus}. 
\end{proof}

\section{An invariant of links in the complement of an unlink}
\lbl{sec.tangles}

The goal of this section is to define an invariant
$$
\Zratc:\NXO\longto\GA(\ostar_X,\La)
$$ 
as stated in Fact 7 of Section \ref{sub.maini}.

Given a finite set $X$, an {\em $X$-indexed sliced \crossed link on $D_g$} 
is a sliced \crossed link on $D_g$ equipped with bijection of the set of 
components of the presented link with $X$.
All the definitions above have obvious adaptations 
for $X$-indexed objects, indicated by an appropriate $X$ subscript when 
required.

\begin{theorem}
\lbl{thm.Zw}
There exists an invariant 
$$
\Zratc:
\cD_X(D_g)/\la r \ra \longrightarrow \GA(\tline_X,\La) 
$$
that maps the topological basing relations to group-like basing relations
and the topological wrapping relations to wrapping relations. Combined
with Proposition \ref{prop.tangles3}, and symmetrizing using 
Lemma \ref{lem.chisigma}, it induces a map
$$
\s_X\circ\Zratc: \NXO \longrightarrow \GA(\ostar_X,\La).
$$ 
\end{theorem}

\begin{proof}
Our first task is to define $\Zratc$. In order to do so, we need to introduce 
some more notation. Given a word $w$ in the letters $\uparrow$ and 
$\downarrow$, let
\begin{itemize}
\item $\uparrow_w$ denote the corresponding skeleton. 
\item $I_w \in \A(\uparrow_w,\La)$ denote the identity element. 
\item $G_{i,w}$ denote the gluing word at the $i$th gluing site, that is, 
$I_w$ with 
a label of $t_i$ on each component, in the sense described by the following
example: if $w=\uparrow \downarrow \downarrow$, then $G_{2,\uparrow\downarrow
\downarrow}=\strutb{}{}{t_2}\stbdown{}{}{t_2^{-1}}\stbdown{}{}{t_2^{-1}}$. 
\end{itemize}

Given a word $w$ in the symbols $\uparrow,\downarrow,\UAO$ and
$\DAO$, let $w'$ denote $w$ with the crosses forgotten.
Let $I_w^\times$ be $I_{w'}$ with every crossed strand broken.
For example,
$$
I^\times_{\UAO \downarrow \DAO} = \setlength{\unitlength}{0.02\standardunitlength}
	\begin{array}{c}  \hspace{-1.7mm}
         	\raisebox{-8pt}{\input draws/Sigma7.epc }
         	\hspace{-1.9mm}
	\end{array}

$$

\begin{definition}
\lbl{def.Znow}
Take an $X$-indexed element $D = \{T_1,\dots,T_k\} \in \cD_X(D_g)$.
\begin{enumerate}
\item
Select $q$-tangle
lifts of each tangle in the sequence, $\{T'_1,\dots,T'_k\}$,
such that
\begin{enumerate}
\item the strands of each $q$-tangle which lie on boundary lines of the 
following form
$$
\setlength{\unitlength}{0.03\standardunitlength}
	\begin{array}{c}  \hspace{-1.7mm}
         	\raisebox{-8pt}{\input draws/Sigma91.epc }
         	\hspace{-1.9mm}
	\end{array}

$$
are bracketed $(w_1)(w_2)$, where $w_1$ (resp. $w_2$) corresponds to 
the strands going to the left (resp. right) of the removed disc,
and these bracketings match up,
\item boundary words at boundaries of the form
$$
\setlength{\unitlength}{0.03\standardunitlength}
	\begin{array}{c}  \hspace{-1.7mm}
         	\raisebox{-8pt}{\input draws/Sigma92.epc }
         	\hspace{-1.9mm}
	\end{array}

$$
are given the canonical left bracketing.
\end{enumerate}
\item Compose the (usual) Kontsevich invariants of these 
$q$-tangles in the following way:
\begin{enumerate}
\item if two consecutive tangles are as follows, 
$$
\{ \dots, \setlength{\unitlength}{0.03\standardunitlength}
	\begin{array}{c}  \hspace{-1.7mm}
         	\raisebox{-8pt}{\input draws/Sigma93.epc }
         	\hspace{-1.9mm}
	\end{array}
, \setlength{\unitlength}{0.03\standardunitlength}
	\begin{array}{c}  \hspace{-1.7mm}
         	\raisebox{-8pt}{\input draws/Sigma94.epc }
         	\hspace{-1.9mm}
	\end{array}
 \dots \}
$$ 
with a word $w$ in the
symbols $\uparrow,\downarrow,\UAO,\DAO$ describing their matching boundary
word, then they are to be composed 
$$
\dots \circ \Z (T'_i)\circ I^\times_{w} \circ \Z (T'_{i+1}) \circ \dots
$$
\item if, on the other hand, two consecutive tangles meet as
follows,  
$$
\{\dots, \setlength{\unitlength}{0.03\standardunitlength}
	\begin{array}{c}  \hspace{-1.7mm}
         	\raisebox{-8pt}{\input draws/Sigma95.epc }
         	\hspace{-1.9mm}
	\end{array}
, \setlength{\unitlength}{0.03\standardunitlength}
	\begin{array}{c}  \hspace{-1.7mm}
         	\raisebox{-8pt}{\input draws/Sigma96.epc }
         	\hspace{-1.9mm}
	\end{array}
 \dots\}
$$
with their matching boundary word at the $j$th gluing site factored as
$(w_1)(w_2)$, then they are to be composed
$$
\dots \circ \Z (T'_i)\circ (I_{w_1}\otimes G_{j,w_2}) \circ 
\Z (T'_{i+1}) \circ \dots\ .
$$
\end{enumerate}
\end{enumerate}
\end{definition}
This completes the definition of $\Zrat(D)$. A normalized version $\Zratc(D)$
(useful for invariance under Kirby moves) is defined by 
connect-summing a copy of $\nu$ into each 
component of the crossed link associated to $D$. Observe that the result does
not depend on the choice of a place to connect-sum to.

The proof of Theorem \ref{thm.Zw} follows from Lemmas \ref{lem.links1} and
\ref{lem.links2}  below.
\end{proof}

The next lemma, observed by the second author in \cite{Kr1}, follows from a 
``sweeping'' argument.

\begin{lemma}
\lbl{gluecomms} 
Let $s$ denote an element of $\A(\tline_X,\BQ) \subset\A(\tline_X,\La)$, 
(ie, some diagram with labels $1$ at each edge), 
where the top boundary word
of $\tline_X$ is $w_1$, and the bottom boundary word of $\tline_X$ is $w_2$. 
Then
$$
G_{w_1}\circ s = s \circ G_{w_2}.
$$
\end{lemma}

\begin{lemma}
\lbl{lem.links1}
The map 
\begin{equation}
\lbl{eq.Zratctemp}
\Zratc: \cD_X(D_g) \longto \GA(\tline_X, \La)
\end{equation}
is well defined and preserves the isotopy relation $\stackrel{r}\sim$ 
of Definition \ref{def.regularisotopy}.
\end{lemma}

\begin{proof}
For the first part we only need to show that the function does not
depend on the choice of bracketings. This follows because of the following,
where $\psi_1,\psi_2$ are some expressions built from the associator.
\begin{eqnarray*}
& & \dots\circ \Z (T_i') \circ (\psi_1\otimes \psi_2) 
\circ (I_{w_1}\otimes G_{j,w_2})
\circ (\psi_1^{-1}\otimes \psi_2^{-1}) \circ \Z (T_{i+1}')\circ \dots, \\
& = & 
\dots\circ \Z (T_i') \circ (I_{w_1}\otimes G_{j,w_2}) \circ 
(\psi_1\otimes \psi_2)
\circ (\psi_1^{-1}\otimes \psi_2^{-1}) \circ \Z (T_{i+1}')\circ \dots, \\
& = & \dots\circ \Z (T_i') \circ (I_{w_1}\otimes G_{j,w_2}) \circ
\Z (T_{i+1}')\circ \dots.
\end{eqnarray*}
where the first step follows from Lemma \ref{gluecomms}.
The invariance of $\Zratc$ under the $\stackrel{r}\sim$ equivalence relation
follows the same way.
\end{proof}

\begin{lemma}
\lbl{lem.links2}
The map $\s_X \circ \Zratc$ respects basing relations.
\end{lemma}

\begin{proof}
It suffices to show that $s_X \circ \Zrat$ respects basing relations. In other
words, it suffices to show that
$D_1 \stackrel{\b_i}{\sim} D_2$ implies that 
$$\s_X(\Zrat(D_1))\stackrel{\bgp_1}{\sim}\s_X(\Zrat(D_2)),$$ 
for $i=1,2$.

In the first case, consider two diagrams $D_1,D_2 \in \cD_X(D_g)$ such that 
$D_1 \stackrel{\b_1}{\sim} D_2$, and let $D_{12}$ denote the diagram in 
$\cD(D_g)'$ such that $D_1$ is obtained by forgetting some cross, and $D_2$ is 
obtained by forgetting the other cross on that component. 

Enlarge $\cD_X(D_g)$ to $\cD_X(D_g)'$, the set of $X$-indexed 
sliced \crossed diagrams on $D_g$ whose presented links in $D_g\times I$ 
have {\it at least} one cross on every component, thus 
$\cD(D_g)\subset \cD(D_g)'$. The function $\Zrat$ immediately extends to 
$\cD(D_g)'$.
The key point is the fact that the product of labels around
each component of the skeleton equals to $1$; this is true since the
locations of the two crosses cut the component into two arcs, both of which
have trivial intersection word with the gluing sites.

In other words, we have:
$\Zrat(D_{12}) \in
\GA(\uparrow_{\{x\}}\uparrow_{\{x'\}}\uparrow_{X-\{x\}},
\La).
$
It follows by construction that
\begin{eqnarray*}
\Zrat(D_1) & = & \vec{m}^{x,x'}_x(\Zrat(D_{12})), \\
\Zrat(D_2) & = & \vec{m}^{x',x}_x(\Zrat(D_{12})), \\
\Zrat(D_{12}) & \in & \GA(\uparrow_{\{x\}}\uparrow_{\{x'\}}
\uparrow_{X-\{x\}},\Lloc).
\end{eqnarray*}
Thus, $\s_X(\Zrat(D_1))\stackrel{\bgp_1}{\sim}\s_X(\Zrat(D_2))$.

In the second case, consider $D_1 \stackrel{\b_2}{\sim} D_2$. Then, 
for some appropriate expression
$\psi$  built from the associator, we have that: 
\begin{eqnarray*}
\lefteqn{
\s_X(\Zrat(
\{\dots,
\setlength{\unitlength}{0.03\standardunitlength}
	\begin{array}{c}  \hspace{-1.7mm}
         	\raisebox{-8pt}{\input draws/Sigma101.epc }
         	\hspace{-1.9mm}
	\end{array}

,
\setlength{\unitlength}{0.03\standardunitlength}
	\begin{array}{c}  \hspace{-1.7mm}
         	\raisebox{-8pt}{\input draws/Sigma102.epc }
         	\hspace{-1.9mm}
	\end{array}

\dots\} )) } & & \\  &=& 
 \s_X( \dots \circ I^{\times}_w \circ \psi \circ ( I_{w_1}\otimes G_{i,w_2} )
\circ \psi^{-1} \circ \dots )\\
& \stackrel{\bgp_1}{\sim} & 
\s_X( \dots \circ \psi \circ I^{\times}_w \circ ( I_{w_1}\otimes G_{i,w_2} )
\circ \psi^{-1} \circ \dots ) \\
& \stackrel{\bgp_2}{\sim} &
\s_X( \dots \circ \psi \circ ( I_{w_1}\otimes G_{i,w_2} )
\circ I^{\times}_w \circ \psi^{-1} \circ \dots )\\
& \stackrel{\bgp_1}{\sim} &
\s_X( \dots \circ \psi \circ ( I_{w_1}\otimes G_{i,w_2} )
\circ \psi^{-1} \circ I^{\times}_w \circ \dots )\\
&= & \s_X( \Zrat(\{\dots,
\setlength{\unitlength}{0.03\standardunitlength}
	\begin{array}{c}  \hspace{-1.7mm}
         	\raisebox{-8pt}{\input draws/Sigma104.epc }
         	\hspace{-1.9mm}
	\end{array}

,
\setlength{\unitlength}{0.03\standardunitlength}
	\begin{array}{c}  \hspace{-1.7mm}
         	\raisebox{-8pt}{\input draws/Sigma103.epc }
         	\hspace{-1.9mm}
	\end{array}

,\dots\})). 
\end{eqnarray*}
Thus, $\s_X(\Zrat(D_1))\stackrel{\bgp_1}{\sim}\s_X(\Zrat(D_2))$.
\end{proof}

\begin{lemma}
\lbl{lem.winding1}
The map $\s_X \circ \Zratc$ preserves the wrapping relations.
\end{lemma}

\proof
Consider two sliced crossed links $D_1, D_2 \in \cD(D_g)$ related by a 
topological \wrapping\ move at the $j$th gluing site as in Section 
\ref{sub.tangles}. In other words, $D_1$ (shown on the left) and $D_2$
(shown on the right) are given by:
$$
\{\dots,\setlength{\unitlength}{0.03\standardunitlength}
	\begin{array}{c}  \hspace{-1.7mm}
         	\raisebox{-8pt}{\input draws/Sigmaw1.epc }
         	\hspace{-1.9mm}
	\end{array}
,\setlength{\unitlength}{0.03\standardunitlength}
	\begin{array}{c}  \hspace{-1.7mm}
         	\raisebox{-8pt}{\input draws/Sigmaw2.epc }
         	\hspace{-1.9mm}
	\end{array}
, \dots \}
\longleftrightarrow
\{\dots, \setlength{\unitlength}{0.03\standardunitlength}
	\begin{array}{c}  \hspace{-1.7mm}
         	\raisebox{-8pt}{\input draws/Sigmaw3.epc }
         	\hspace{-1.9mm}
	\end{array}
, \setlength{\unitlength}{0.03\standardunitlength}
	\begin{array}{c}  \hspace{-1.7mm}
         	\raisebox{-8pt}{\input draws/Sigmaw4.epc }
         	\hspace{-1.9mm}
	\end{array}
, \dots \}
$$
Let
\begin{eqnarray*}
z_1 &=& \Zratc \left(
\setlength{\unitlength}{0.03\standardunitlength}
	\begin{array}{c}  \hspace{-1.7mm}
         	\raisebox{-8pt}{\input draws/Sigmaw5.epc }
         	\hspace{-1.9mm}
	\end{array}

 \right), \\
z_1' &=& \strutb{}{x'}{e^h} \cup 
(\s_x(z_1))_{x \to \pt h}, \\
z_2 &=& \Zratc \left( 
\setlength{\unitlength}{0.03\standardunitlength}
	\begin{array}{c}  \hspace{-1.7mm}
         	\raisebox{-8pt}{\input draws/Sigmaw6.epc }
         	\hspace{-1.9mm}
	\end{array}

\right) 
\end{eqnarray*}
Now we can compute, for some element $\psi$ built from an associator:
\begin{eqnarray*}
\lefteqn{
(\sigma_X\circ\Zratc)(D_2) } & & \\
& = &
\sigma_X \circ \left( \con_{\{h\}}\left( 
\dots \circ \psi \circ \Delta^{w_2}_{x'}(z_1') \circ \left(I_{w_1}\otimes
G_{j,w_2}\right) \circ \Delta^{w_2}_x(z_2) \circ \psi^{-1} \circ \dots
\right) \right) \\
& \stackrel{\wgp}{\sim} &
\sigma_X \circ \left(  \con_{\{h\}}\left( 
\dots \circ \psi \circ \left(I_{w_1}\otimes
G_{j,w_2}\right) \circ \Delta^{w_2}_{x'}(z_1') \circ 
\Delta^{w_2}_x(z_2) \circ \psi^{-1} \circ \dots
\right) \right), \\
& = & (\sigma_X\circ \Zratc)(D_1). 
\end{eqnarray*}
In the above, $\Delta^{w_2}_{x'}$ means the comultiplication of the diagram
$x$ with pattern given by the gluing word $w_2$. For example,
$$
\D^{\uparrow\downarrow\downarrow\uparrow} \left(\stbup{}{}{e^h}\right) =
\stbup{}{}{e^h} \stbdown{}{}{e^{-h}} \stbdown{}{}{e^{-h}} 
\stbup{}{}{e^h}.\eqno{\qed}
$$

Theorem \ref{thm.Zw} now follows, using the identification of $\NXO$
with the set \newline $\cD_X(D_g)/\la r,\b,\w \ra$ given by Proposition 
\ref{prop.tangles}.

\begin{remark}
\lbl{rem.renormalization}
An application of the $\Zrat$ invariant is a formula for
the LMO invariant of cyclic branched covers of knots, \cite{GK1}.
This application requires a {\em renormalized version} $\Zrat[\a]$ where $\a=
(\a_1,\dots,\a_g)$ and $\a_i \in \GA(\star_{k})$. $\Zrat[\a]$ is defined
using Definition \ref{def.Znow} and replacing the formula in 2(b) by
$$
\dots \circ \Z (T'_i)\circ (I_{w_1}\otimes G_{j,w_2} \Delta(\a_j) ) \circ 
\Z (T'_{i+1}) \circ \dots .
$$
where $\Delta(\a_j)$ places a copy of $\a_j$ in each of the skeleton segments
of the $j$th gluing site. The renormalized invariant $\Zrat[\a]$ still satisfies
the properties of Theorem \ref{thm.Zw}.
\end{remark}

\section{A rational version of the Aarhus integral}
\lbl{sec.aarhus}

\subsection{What is Formal Gaussian Integration?}
\lbl{whatisFGI}
In this section, we develop a notion of Rational Formal Gaussian Integration
$\intrat$, which is invariant under basing, wrapping and Kirby move relations.

Our theory is a cousin of the Formal Gaussian Integration of \cite{A},
and requires good book-keeping but essentially no new ideas.

Before we get involved in details, let us repeat a main idea from \cite{A}:
Formal Gaussian Integration is a theory of contraction of legs of diagrams.
As such, it does not care about the internal structure of diagrams (such as
their internal valency) or the decoration of the internal edges of diagrams.

If diagrams are thought as tensors (as is common in the world of perturbative
Quantum Field Theory) which represent differential operators with polynomial
coefficients, then contractions of legs corresponds to differentiation.

The reader is encouraged to read \cite[part II]{A} for a lengthy introduction
to the need and use of Formal Gaussian Integration.

\subsection{The definition of $\intrat$}
\lbl{sub.inteq}

Recall that $\Herm(\Lloc\to\BZ)$ denotes the set of Hermitian matrices 
over $\Lloc$, invertible over $\BZ$. 
Let $X$ denote a finite set and let $X'=\{x_i\}$ denote an arbitrary
subset of $X$. A diagram $D\in \A(\star_{X},\Lloc)$ is called
$X'$-{\em substantial} if there are no {\em strut} components
(ie, components comprising of a single edge) both of whose
univalent vertices are labeled from the subset $X'$.
An element of $\A(\star_{X},\Lloc)$ is $X'$-substantial if it
is a series of $X'$-substantial diagrams.

\begin{definition}
\lbl{def.integrable}
$\mathrm{(a)}$
An element $s\in \A(\star_X,\Lloc)$ is called {\em integrable with
respect to} $X'$ if there exists an Hermitian matrix
$M\in\Herm(\Lloc\to\BZ)$ such that
\begin{equation}
\lbl{eq.Gaussd}
s = \ess \left( \frac{1}{2} \sum_{i,j}
 \strutb{x_j}{x_i}{M_{ij}} \right) \sqcup R,
\end{equation}
with $R$ an $X'$-substantial element. Notice that $M$, the
{\em covariance matrix $\cov(s)$} of $s$, and $R$, the $X'$-{\em substantial 
part} of $s$, are uniquely determined by $s$. \newline
$\mathrm{(b)}$
Let $\IA\subset \A(\star_X,\Lloc)$ denote
the subset of elements that are {\em integrable with respect} to $X'$.
\newline
$\mathrm{(c)}$
Observe that a basing move preserves integrability. Thus, we can define
$$
\IGA = \left(\GAX \cap \IA \right)/\la \bgp_1,\bgp_2\ra.
$$ 
\end{definition}

For the next definition, recall that every matrix $M \in \Herm(\Lloc\to\BZ)$
is invertible over $\Lloc$.

\begin{definition}
Define a map ({\em Rational Formal Gaussian Integration})
$$
\intrat dX': \IA \rightarrow \A( \star_{X-X'} , \Lloc )
$$
as follows. If $s\in \IA$ has associated decomposition given in Equation 
\eqref{eq.Gaussd}, then
$$
\intrat dX'(s) = 
\left\la  
\exp_{\sqcup}\left( -\frac{1}{2} \sum_{i,j} 
\strutb{\pt x_j}{\pt x_i}{M_{ij}^{-1}} \right)
,R\right\ra_{X'}.
$$
\end{definition}
The gluing of the legs of diagrams occurs according to Remark \ref{rem.glue}.

\begin{remark}
\lbl{rem.Cohn} 
The definition of $\intrat$ for elements $s$
whose covariance matrix $M$ is defined over $\La$, and invertible over $\BZ$, 
requires
beads labeled by the inverse matrix $M^{-1}$. This necessitates the need to
replace $\La$ by an appropriate ring $\Lloc$ that makes all such matrices
$M$ invertible. It follows by definition that the universal choice of
$\Lloc$ is the {\em Cohn localization} of $\La$. 
\end{remark}

\subsection{$\intrat$ respects the basing relations}
\lbl{sub.inteqbasing}

The next theorem shows that $\intrat$ is compatible with the basing relations.

\begin{theorem}
\lbl{thm.cyclic}
$\intrat$ descends to a map:
$$
\intrat dX' : 
\IPR_{X'}\GA(\star_X,\Lloc)/\la \bgp_{X} \ra
\rightarrow
\GA(\star_{X-X'},\Lloc)/\la \bgp_{X-X'} \ra
$$
\end{theorem}

\begin{proof}
Consider  $s_1,s_2\in \IA \cap \GAX$ such that
$s_1 \stackrel{\bgpp_1}{\sim} s_2$ or $s_1 \stackrel{\bgp_2}{\sim} s_2$
(based on some element $x_k\in X'$).  Using Lemma \ref{lem.b1},
we need to show that 
$$
\intrat dX'\left( s_1 \right) = \intrat dX' \left( s_2 \right).
$$

\noindent
{\em The case of $\bgp_2$}.  
Consider the decomposition of $s_1$ 
$$
s_1 = \ess\left( \frac{1}{2}\sum \strutb{x_j}{x_i}{M_{ij}} \right) \sqcup R.
$$
Suppose that an element $f\in F_g$ is pushed onto the legs 
labeled by $x_k\in X'$, to form $s_2$. Now, let $D_{f,k}$ be 
$\diag(1,1,\dots,f,\dots,1)$, where $f$ is in the $k$th entry, and let $R'$ 
denote $R$ with an $f$ pushed onto every leg marked $x_k$. Then,
\begin{eqnarray*}
\intrat dX' \left(s_2\right) & = & \intrat dX' \left( 
\ess\left( \frac{1}{2} 
\sum \strutb{x_j}{x_i}{(D_{f,k}MD_{f^{-1},k})_{ij}} \right) \sqcup R' 
\right) \\
& = &
\left\la  \ess\left( -\frac{1}{2}
\sum \strutb{x_j}{x_i}{(D_{f,k}M^{-1}D_{f^{-1},k})_{ij}} \right)
, R' \right\ra  \\
& = &
\left\la  \ess\left( -\frac{1}{2}
\sum \strutb{x_j}{x_i}{(M^{-1})_{ij}} \right)
, R \right\ra  \\
& = & \intrat dX' \left( s_1 \right)
\end{eqnarray*}
because all the beads labeled $f$ and $f^{-1}$ match up and cancel,
after we pairwise glue the $X'$-colored legs.

\noindent
{\em The case of $\bgpp_1$}. 
Consider $s \in \GA(\star_{X \cup \{\pt h\}},\Lloc)$ such that
$$
s_1=\conh(s) \quad \text{ and } \quad
s_2=\conh(s|_{x_k\to x_k e^h}).
$$
Notice that if $s_1,s_2$ are $X'$-integrable, so is $s$.
Furthermore, $\conh$ commutes with $\intrat$. The result follows
from Lemma \ref{lem.manyhair} below.
\end{proof}

The next lemma says that pushing $e^h$ on $X'$-colored legs commutes
with $X'$-integration. Of course, after $X'$-integration there is
no $X'$-colored leg, thus $X'$-integration is invariant under
pushing $e^h$ on $X'$-colored legs.

\begin{lemma}
\lbl{lem.manyhair}
Let $s\in \IPR_{X'}\GA( \star_{X \cup \{\pt h\}}, \Lloc)$, and 
let $x\in X'\subset X$.
$$
\intrat dX' \left( \left.s\right|_{x\rightarrow xe^h} \right)
=
\intrat dX' \left( s \right) \in \GA(\star_{X \cup \{\pt h\}},\Lloc).
$$
\end{lemma}

\begin{proof}
The idea of this lemma is the same as ``in the case $\bgp_2$'', but 
instead of pushing an element of the free group onto legs labeled by some 
element, we push an exponential of hair. 

The substitution $x \to xe^h$ defines a map
$$
\A(\star_{X \cup \{\pt h\}},\Lloc)\longrightarrow 
\A(\star_{X \cup \{h, \pt h\}},\Lloc)$$
with the property that $s-s|_{x\to xe^h}$ is $X'$-substantial for all
$s \in \A(\star_X,\Lloc)$.

Assume that the canonical decomposition of $s$ is given by Equation 
\eqref{eq.Gaussd}. Let us define
$$
D_{e^h,k}=\diag(1,1,\dots,e^h,\dots,1),
$$ 
where 
$e^h$ is in the $k$th entry. Then, Theorem \ref{thm.varphi} from Appendix
\ref{sub.useful} implies that  
\begin{eqnarray*}
\lefteqn{
\intrat dX'\left( \left.s\right|_{x\rightarrow xe^h} \right)} & & \\
& = &
\es\left(-\frac{1}{2} \chiz \left( 
MD_{e^h,k}M^{-1}D_{e^{-h},k}\right) 
\right)
\sqcup
\left\la \es\left( -\frac{1}{2} 
\sum \strutb{x_j}{x_i}{M^{-1}_{ij}} \right)
,R \right\ra \\
& = &
\left\la 
\es\left( -\frac{1}{2} 
\sum \strutb{x_j}{x_i}{M^{-1}_{ij}} \right)
,R \right\ra 
\end{eqnarray*}
where the last equation follows from Proposition \ref{prop.w}.
\end{proof}

\begin{remark}
An alternative proof of Theorem \ref{thm.cyclic} can be obtained by adapting
the proof of \cite[Part II, Proposition 5.6]{A} to the rational 
integration $\intrat$.
\end{remark}

\subsection{$\intrat$ respects the Kirby moves}
\lbl{sub.kirby}

Recall from Section \ref{sub.tangles} the $\kappa$ handleslide move on the 
set $\cD(D_g)$. 
In a similar fashion, we define an $X'$-{\em handleslide move}
on $\A(\star_X,\Lloc)$, as follows: For two elements $s_1,s_2 \in 
\A(\star_X,\Lloc)$, say that $s_1\stackrel{\kappa_{X'}}{\sim} s_2$, 
if for some $x_i\in X,$ and $x_j\in X',\ i\neq j$,
$$
s_2 = \left(\s_{\{x_i\}} \circ \vec{m}^{y,x_i}_{x_i}\circ \chi_{\{y,x_i\}} 
\right)
\left(s_1\left(x_1,\dots,x_{j-1},x_j+y,x_{j+1},\dots \right)\right).
$$
Notice that 
if $s_1,s_2 \in \GAX$ and $s_1 \stackrel{\kappa_{X'}}{\sim} s_2$, then if
one of
$s_i$ lies in  
$\IA$, so does the other. This is because the covariance matrix
of $s_2$ is obtained from the covariance matrix of $s_1$ by a unimodular 
congruence.

\begin{theorem}
\lbl{thm.kirby1}
The map $\s_X \circ \Zratc$ preserves the Kirby relations.
\end{theorem}

\proof
The argument of \cite[Proposition 1]{LMMO}, adapts to this setting without 
difficulty. There are two possible subtleties that must be observed.
\begin{itemize}
\item[(a)] 
The proof in \cite{LMMO} is based on taking a parallel of the parts of
a component (here, the component $x_j$) away from the site of the cross. 
In our setting, the paralleling operation of the involved components of
tangles behaves as usual (sums over lifts of legs). At a gluing
word, one simply labels both of the two components 
arising from the parallel operation with a bead labeled by
the appropriate letter. Then, when $\s_{\{x_i,y\}}$ is applied,
these two identical beads are pushed onto the legs sitting on those components 
side-by-side. This may alternatively be done by first pushing the original
bead up onto the legs, {\em and then} taking the parallel, as is required in 
this situation.
\item[(b)]
We have defined the bracketing at boundary lines with crosses
to be the standard left-bracketing, while the move in \cite{LMMO} is
given for a presentation where the two involved strands are bracketed 
together. This is treated by preceding and following the $\kappa_{X'}$ with
a $\b_1$ move, to account for the change in bracketing, exactly as in
the proof of the $\b_2$-case of Lemma \ref{lem.links2}.\endproof 
\end{itemize}

\begin{theorem}
\lbl{thm.kirby2}
The map $\intrat dX'$ descends to a map:
$$
\intrat dX' : 
\IPR_{X'}\GA(\star_X,\Lloc)/\la \kappa_{X'}\ra \to \GA(\star_{\,X-X'},\Lloc).
$$
\end{theorem}

This is proved by the following lemma.

\begin{lemma}
Let $s_1 \in \IA$
and let 
$$
s_2 = \left(\s_{\{x_i\}} \circ \vec{m}^{y,x_i}_{x_i}\circ \chi_{\{y,x_i\}} \right)
\left(s_1\left(x_1,\dots,x_{j-1},x_j+y,x_{j+1},\dots \right)\right).
$$
If $x_j \in X'$ then
$$
\intrat dX' (s_1) = \intrat dX' (s_2).
$$
\end{lemma}

\begin{proof}
We prove this in two steps, which usually go by the names of
the ``lucky'' and ``unlucky'' cases.

{\em The lucky case: $s_1,s_2 \in\ \IPR_{\{x_j\}}\A(\star_X,\Lloc)$.}
In this case, we can use the Iterated Integration Lemma \ref{itintlem}.
{\small
\begin{eqnarray*}
\lefteqn{ \intrat dX'(s_2) } & & \\
& = &
\intrat d(X' - \{x_j\}) \left( \intrat d\{x_j\}(s_2) \right), \\
& = &
\intrat d(X' - \{x_j\}) \\ & & 
 \left(
\intrat d\{x_j\}
\left(\s_{\{x_i\}} \circ \vec{m}^{y,x_i}_{x_i}\circ \chi_{\{y,x_i\}} \right)\left( s_1(x_1,\dots,x_{j-1},x_j+y,x_{j+1},\dots) \right) \right), \\
& = & 
\intrat d(X' - \{x_j\}) \\
& & \left(
\left(\s_{\{x_i\}} \circ \vec{m}^{y,x_i}_{x_i}\circ \chi_{\{y,x_i\}} \right)\left( \intrat d\{x_j\} 
\left( s_1(x_1,\dots,x_{j-1},x_j+y,x_{j+1},\dots) \right) \right)\right), \\
& = & 
\intrat d(X' - \{x_j\}) \left(
\left(\s_{\{x_i\}} \circ \vec{m}^{y,x_i}_{x_i}\circ \chi_{\{y,x_i\}} \right)\right. \\
& & \\
& & \left. \left( \intrat d\{x_j\} 
\left( 
\es\left(
\gau{1}{y}{\pt x_j} 
\right)
\flat_{\{x_j\}}
s_1(x_1,\dots,x_{j-1},x_j,x_{j+1},\dots) \right) \right)\right) \\
& & 
\end{eqnarray*}
}
Now, using the Integration by Parts Lemma \ref{intbypartslem}, 
it follows that
\begin{eqnarray*}
\lefteqn{ \intrat d\{x_j\} \left(\es\left(
\gau{1}{y}{\pt x_j} 
\right) \flat_{\{x_j\}} s_1(x_1\dots)\right) } & & \\
& = &
\intrat d\{x_j\} \left( \div_{\{x_j\}} \left( \es\left(
\gau{1}{y}{\pt x_j}  \right) \right) 
\flat_{\{x_j\}} s_1(x_1,\dots) \right) \\
& = &
\intrat d\{x_j\} \left( s_1(x_1,\dots) \right),
\end{eqnarray*}
because 
$\div_{\{x_j\}} \left( \es\left(\strutb{y}{\pt x_j}{1} \right) \right)= 1.$

{\em The unlucky case: $s_1,s_2 \notin\ \IPR_{\{x_j\}} \A(\star_X,\Lloc)$.}
This pathology is treated by deformation. Namely, we ``deform''
$s_1$ by multiplying it by
$
\es\left( \epsilon 
\gau{1}{x_j}{x_j} 
\right),
$
for $\epsilon$ a ``small real,'' commuting with $\Lloc$. Formally, we let
$\epsilon$ be a variable and consider an integration theory based on 
$\BQ(\epsilon)[F_g]$ and its Cohn localization $\Lloc^\epsilon$. 
We denote by $\BQ^{0}(\epsilon) \subset \BQ(\epsilon)$, the subring
of rational functions in $\epsilon$ non-singular at zero, and we denote
by $\Lloc^{\epsilon,0}$ the corresponding noncommutative localization 
of $\BQ^{0}(\epsilon)$.

Letting
\begin{eqnarray*}
s_1^\epsilon & = & \es\left( 
\epsilon
\gau{1}{x_j}{x_j} 
\right) \sqcup s_1 \in \A(\star_X,\Lloc^\epsilon),\\
& & \\
s_2^\epsilon & = &
\left(\s_{\{x_i\}} \circ \vec{m}^{y,x_i}_{x_i}\circ \chi_{\{y,x_i\}} 
\right)
\left(s_1^\epsilon\left(x_1,\dots,x_{j-1},x_j+y,x_{j+1},\dots\right)\right),
\end{eqnarray*}
it follows from the argument of the ``lucky case'' that
$$
\intrat dX'\left( s_2^\epsilon \right) = \intrat dX' 
\left( s_1^\epsilon \right) \in \A(\star_X,\Lloc^{\epsilon,0}).
$$
There is a ring homomorphism $\phi: \Lloc^{\epsilon,0}
\rightarrow \Lloc,$ defined by ``setting $\epsilon$ to zero''. 
Since 
$$
\phi\left(\intrat dX'\left( s_i^\epsilon \right) \right) 
= \intrat dX'\left( s_i \right), 
$$
the result follows.
\end{proof}

\begin{remark}
These last lines are best understood from the realization of
$\Lloc$ and $\Lloc^\epsilon$ as spaces of formal power series.
The elements  $\intrat dX'\left( s_i^\epsilon\right)$ lie
in $\A(\star_X,\Lloc^{\epsilon,0})$
because the determinant of the augmentation of the covariance matrix of $s_i$
is invertible in $\BQ^0(\e)$.
\end{remark}

\subsection{$\intrat$ respects the \wrapping\ relations}
\lbl{sub.wrapping}

In this section we show that $\intrat$-integration respects the group-like
wrapping relations $\wgp$ that were introduced in Definition 
\ref{def.wrapping}. Since the wrapping relations involve pairs of matrices
and graphs, we begin by extending the definition of $\intrat$ to pairs 
$(M,s)$ by setting
$$
\intrat dX'(M,s)=(M\oplus \cov(s), \intrat dX' (s)),
$$
where $\cov(s)$ is the covariance matrix of $s$. The motivation for
this comes from the next theorem which states that $\intrat$ preserves the 
\wrapping\ relations.

\begin{theorem}
\lbl{thm.winding2}
The map $\intrat dX'$ descends to a map
$$
\intrat dX' : 
\IPR_{X'}\GAz(\star_X,\Lloc)/\la \wgp \ra
\rightarrow
\GAz(\star_{X-X'},\Lloc)/\la \wgp \ra
$$
\end{theorem}

\begin{proof}
Consider two pairs $(M,s_1) \stackrel{\wgp}\sim (M,s_2)$ related by
a \wrapping\ move, as in Definition \ref{def.wrappingD}. Observe 
that a wrapping move preserves the covariance matrix. 

Using Theorem \ref{thm.varphi} in Appendix \ref{sub.useful}, it follows
that if $W$ is the common covariance matrix of $s_1,s_2$ and $s$, then 
\begin{eqnarray*}
\lefteqn{ \intrat dX' \left(M, s_m \right) } & & \\
 &=& \left(M\oplus W, \conh \left( 
\intrat dX' 
\left( 
\varphi_m (s) \right)
\sqcup \left( 
\exp \left(-\frac{1}{2} \chiz(M^{-1}
\varphi_m M) \right)\right)
\right)\right) \\
&=&
\left(M\oplus W, \conh \left( 
\varphi_m \left( \intrat dX' 
 (s) \right)
\sqcup \left( 
\exp \left(-\frac{1}{2} \chiz(W^{-1}
\varphi_m W) \right)\right) \right.\right. \\
& &
\left. \left. \sqcup \left( 
\exp \left(-\frac{1}{2} \chiz(M^{-1}
\varphi_m M) \right)\right)
\right)\right) \\
& = &
\left(M\oplus W, \conh \left( 
\varphi_m \left( \intrat dX' 
 (s) \right)
 \right.\right. \\ & &  \left.\left. \sqcup \left( 
\exp \left(-\frac{1}{2} \chiz((W\oplus M)^{-1}
\varphi_m (W\oplus M)) \right)\right)
\right)\right) \\
& \stackrel{\wgp}\sim &
\left(M, \intrat dX' \left(M, s \right)\right) .
\end{eqnarray*}\vspace{-12mm}

\end{proof}

\begin{remark}
\lbl{rem.wrapping}
The group-like \wrapping\ relation is the minimal relation that we need to 
impose on diagrams so that $\intrat \circ \s \circ \Zrat$ is invariant under 
the topological \wrapping\ relation, as follows from the proof of Theorem 
\ref{thm.winding2}.
\end{remark}

\subsection{The definition of the rational invariant $\Zrat$ of $F$-links}
\lbl{sub.defineZ}

We have all the ingredients to define our promised 
invariant of $F$-links. 

Theorems \ref{thm.Zw}, \ref{thm.cyclic}, \ref{thm.kirby2}, 
and \ref{thm.winding2} imply that:

\begin{definition}
\lbl{def.ZratF}
There is a map 
$$
\Zratgp: \text{$F$-links} \longrightarrow \GAz(\Lloc)/\la \wgp \ra
$$
defined as follows: for a $F$-link $(M,L)$ of $g$ components,
choose a sliced crossed link $C$ in $\cD(D_g)$ whose components are in 1-1
correspondence with a set $X$, and define
$$
\Zratgp(M,L)=\frac{\intrat dX \left(\s \circ \Zratc(C) \right)}{
c_+^{\sigma_+(B)} c_-^{\sigma_-(B)}} \in \GAz(\Lloc)/\la \wgp \ra
$$
where $c_\pm=\int dU \Zc (S^3,U_{\pm})$ are some universal constants
of the unit-framed unknot $U_{\pm}$ and where $B$ is the linking matrix of
$C$ and $\s_{\pm}(B)$ is the number of positive (resp. negative) eigenvalues
of the symmetric matrix $B$.
\end{definition}

The slight renormalization given by multiplication by $c_{\pm}$ is needed to 
deal with the stabilization Kirby move on $C$, and is justified as follows:
recall that $\GA(\phi)$ is a group, which acts on the set $\GAz(\Lloc)$ by 
$a \cdot (M,s)=(M, a \sqcup s)$. Since the beads of the graphs in 
$\GA(\phi)$ are rational numbers independent of the $t_i$, 
it follows easily that the above action of $\GA(\phi)$ on $\GAz(\Lloc)$ 
induces an action of $\GA(\phi)$ on the quotient set 
$\GAz(\Lloc)/\la \winf \ra$. This is exactly what appears in the above 
definition of $\Zrat$.

\section{The $\CA$ action and invariants of boundary links}
\lbl{sec.CA}

In this section we show how the invariant of $F$-links of Definition 
\ref{def.ZratF} can be modified in a natural way to give an invariant
of boundary links.

The key geometric observation is that the set of boundary links can be 
defined with the set of orbits of an action
$$
\CA \times \text{$F$-links} \to \text{$F$-links}
$$
for a suitable group $\CA$. 
Starting from this observation, we refine the results of Sections 
\ref{sec.surgeryview}-\ref{sec.aarhus} to take into account the action
of the group $\CA$.

\subsection{A surgery description of boundary links}
\lbl{sub.surgeryb}

We begin by giving a surgery description of boundary links.
Recall that the set of $\pt$-links can be identified
as a quotient of the set of $F$-links
modulo an action of a group of automorphisms $\CA_g$ of the free group,
\cite{K,CS2}, where
$$
\CA_g = \{ f \in \text{Aut}(F) 
| f(t_i)=\a_i^{-1} t_i \a_i, \quad i=1,\dots,g \}.
$$
In particular, $\CA_g$ contains (and for $g=1,2$ coincides with) the group
of inner
automorphisms of $F_g$. It is easy to see that $\CA_g$ is generated by
the special automorphisms $\a_{ij}$ (for $1 \leq i \neq i \leq g$)
acting on the generators of $F_g$ by 
$$
\a_{ij}(t_k)= 
\begin{cases}
t_k & \text{if} \quad k \neq i \\
t_j^{-1} t_i t_j & \text{if} \quad k=i;
\end{cases}
$$ 
see \cite[Lemma 2.4]{K}. The group $\CA_g$ acts on the set of $F$-links 
by changing the basing of the $F$-link, in other words by
sending $(M,L_b,\phi)$ to $\a \circ \phi$ for $\a \in \CA_g$.
An alternative and more geometric way to describe this action is as follows:
identify the set of $F$-links with the set of Seifert surfaces modulo
isotopy and {\em tube equivalence}, see \cite{GL2}. There is an action of 
$\CA_g$ on the set of Seifert surfaces, described on generators $\a_{ij}$
by {\em cocooning} the $i$th surface to the $j$th surface using suitable arcs
to a base point, well-defined modulo isotopy and tube equivalence.

Since the set $\NO/\la \kappa \ra$ is in 1-1 correspondence with the set
of $F$-links and since $\CA$ acts on the latter set, it follows that it
acts on the former set, too.
We now describe this action of $\CA_g$ (or rather, of its generating 
set $\{\a_{ij} \}$) on the set $\NO/\la \kappa \ra$.

Consider a link $L \in \NO$. The action of $\a_{ij}$
on such a link is the following:
``the $i$th unlink component is isotoped through the $j$th unlink component,
and then put back where it started''.

Let us examine what happens when we try to implement such a 
principle. To do this, we take a parallel copy of the $i$th unlink component 
and drag it so that, at the end, the initial link (resp. the link that 
results from the move) is recovered by forgetting the parallel copy (resp. 
forgetting the original copy).

\begin{figure}[ht!]
$$ 
\psdraw{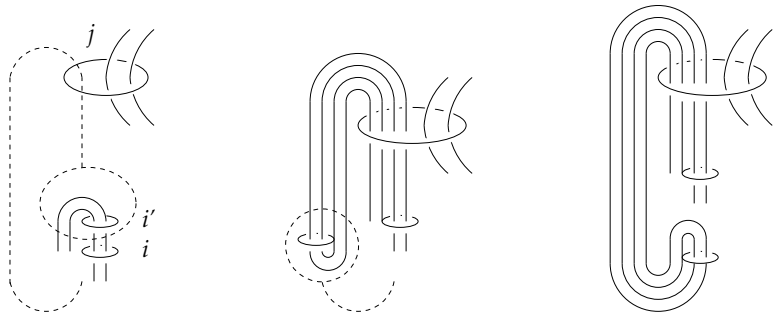}{4in} 
$$
\caption{Dragging the $i$th unlink component  through the $j$th unlink 
component, and then put back where it started.}\lbl{Drag}
\end{figure}

It is easy to see that the $\a_{ij}$ move does not depend on the choice of 
arc, modulo Kirby moves, and that the above discussion together with 
Theorem \ref{thm.GK2} imply that:

\begin{theorem}
\lbl{thm.surgeryb}
There exists an action
$$
\CA \times (\NO/\la\kappa\ra) \longrightarrow \NO/\la\kappa\ra
$$
which induces a 1-1 onto correspondence 
$$
\NO/\la\kappa,\CA\ra\longleftrightarrow \text{$\pt$-links} .
$$
\end{theorem}

\begin{remark}
\lbl{rem.kirbyIII}
It is instructive to compare the above theorem with a result of Roberts 
\cite{Rb}, which states that two 3-manifolds $N_L$ and $N_{L'}$ are 
diffeomorphic relative boundary iff $L$ and $L'$ are related by a sequence
of handle-slides ($\kappa_1$), stabilization ($\kappa_2$), and 
insertion/deletion of an arbitrary knot in $N$ together with a $0$-framed 
meridian of it ($\kappa_3$). In case the knot participating in the $\kappa_3$
move is nullhomotopic, the move is equivalent to a composition of $\kappa_1$
and $\kappa_2$ moves, as was observed by Kirby and Fenn-Rourke.
In general, the set $\NO$ is not closed under the moves $\kappa_3$, however
we can consider the set $\widehat{\NO}$ that consists of all links in $S^3-\O$
that are obtained from nullhomotopic links by $\kappa_{1,2,3}$ moves.
Roberts' theorem implies that there is a 1-1 onto correspondence
$$
\widehat{\NO}/\la \kappa_1,\kappa_2,\kappa_3\ra 
\longleftrightarrow \text{$\pt$-links} .
$$
We leave it as an exercise to show that the links of Figure \ref{Drag}
are related by a sequence of $\kappa_{1,2,3}$ moves. \qed
\end{remark}

In order to give a tangle description of the set of $\pt$-links, we need
to introduce a geometric action of the set $\CA_g$ on the set $\cD(D_g)$.
Recalling the action of $\CA_g$ on the set $\NO$, we
see that this action can be specified by a framed arc, and that as this 
proceeds, then all the strands going between the two copies are 
parallel to each other (that is, dragging the
copy away from the original ``lines-up'' all the strands in between); see
also Figure \ref{Drag}. This will be the idea of our action of $\CA_g$ on 
$\cD(D_g)$. 
To formalize this, we (briefly) enlarge our notion of sliced crossed links, 
to allow for such an arc. That is, we will consider the set 
$\cDarc(D_g)$ 
of {\em arc-decorated sliced crossed diagrams}, represented by sequences of 
tangles 
such that, in addition to all closed components, there is an arc which
starts and finishes at some meridional disc:

$$
\setlength{\unitlength}{0.03\standardunitlength}
	\begin{array}{c}  \hspace{-1.7mm}
         	\raisebox{-8pt}{\input draws/SigmaA.epc }
         	\hspace{-1.9mm}
	\end{array}
 \hspace{2cm} \setlength{\unitlength}{0.03\standardunitlength}
	\begin{array}{c}  \hspace{-1.7mm}
         	\raisebox{-8pt}{\input draws/SigmaA2.epc }
         	\hspace{-1.9mm}
	\end{array}

$$
For an element of this set, we will say that the result 
of a {\em drag} along the arc is the sliced crossed link that arises 
by making the following replacements, accompanied by taking a parallel 
of the arc 
$a$ according to the appropriate word (in the example below, according
to $w_b\overline{w_b}$).
$$
\setlength{\unitlength}{0.03\standardunitlength}
	\begin{array}{c}  \hspace{-1.7mm}
         	\raisebox{-8pt}{\input draws/SigmaA.epc }
         	\hspace{-1.9mm}
	\end{array}
 \stackrel{\mathrm{drag}}\longrightarrow
  \setlength{\unitlength}{0.03\standardunitlength}
	\begin{array}{c}  \hspace{-1.7mm}
         	\raisebox{-8pt}{\input draws/SigmaA3.epc }
         	\hspace{-1.9mm}
	\end{array}
 \hspace{2cm}
\setlength{\unitlength}{0.03\standardunitlength}
	\begin{array}{c}  \hspace{-1.7mm}
         	\raisebox{-8pt}{\input draws/SigmaA2.epc }
         	\hspace{-1.9mm}
	\end{array}
 \stackrel{\mathrm{drag}}\longrightarrow
\setlength{\unitlength}{0.03\standardunitlength}
	\begin{array}{c}  \hspace{-1.7mm}
         	\raisebox{-8pt}{\input draws/SigmaA4.epc }
         	\hspace{-1.9mm}
	\end{array}

$$
Given two sliced crossed links $L_1,L_2 \in \cD(D_g)$, we say that $L_2$ is
obtained from $L_1$ by an $\alpha_{ij}$ move if there exists an arc-decorated
sliced crossed link $L_{12} \in \cDarc(D_{g+1})$ with the following 
properties:
\begin{enumerate}
\item 
there are no intersections of link components with the $(i+1)$st meridianal
disc,
\item
gluing back in the $(i+1)$st removed disc (and forgetting the arc)
gives a diagram which presents $L_1$,
\item
the arc travels from the $i$th meridianal disc, intersects only the
$j$th meridianal disc, and only once, and then ends at the $i$th meridianal
disc.
\item
Dragging $L_{12}$ along the arc and then gluing back the 
$i$th removed unlink component recovers $L_2$. 
\end{enumerate}

\begin{remark}
The above definition implements the procedure described in Figure \ref{Drag}. 
Actually, this $\a_{ij}$ move does not depend on the choice of arc (this
is an instructive 
exercise). In addition, this move translated to the set of $F$-links
coincides with the action of $\CA_g$ by cocooning.
\end{remark}

Proposition \ref{prop.tangles} and Theorem \ref{thm.surgeryb} imply that:

\begin{proposition}
\lbl{prop.tanglesbb}
There is an action
$$
\CA \times \cD(D_g)/\la r, \b, \w, \kappa \ra \longleftrightarrow
\cD(D_g)/\la r, \b, \w, \kappa \ra 
$$
which induces a 1-1 correspondence 
$$
\cD(D_g)/\la r, \b, \w, \kappa, \CA \ra \longleftrightarrow
\text{$\pt$-links}.
$$
\end{proposition}

\subsection{The action of $\CA$ on diagrams}
\lbl{sub.CAinf}

We now define an algebraic action of $\CA_g$ on diagrams.
Recall the group $\CA_g$ that acts on the free group $F_g$ (see Section
\ref{sub.surgeryb}), thus also on $\La$. The action of $\CA_g$ on $\La$
extends to the localization $\Lloc$, since for $\a \in \CA_g$, the induced
ring homomorphism $\La \to \La \subset \Lloc$ is $\Sigma$-inverting, and
thus extends to a ring homomorphism $\Lloc\to\Lloc$. 
This induces an action 
$$\CA_g \times \A^0(\star_X,\Lloc) \to \A^0(\star_X,\Lloc)
$$ 
given by acting on the beads of diagrams and the entries of matrices
by special automorphisms of the free 
group. Following Remark \ref{rem.lprp}, we denote by $\A^0(\star_X,\Lloc)/\la
\CA \ra$ the quotient space.

Similarly, there is an action 
$$
\CA_g \times \GAz(\star_X,\Lloc) \to \GAz(\star_X,\Lloc),
$$ 
which results in a quotient set $\GAz(\star_X,\Lloc)/\la \CA \ra$, too.

The next theorem shows that the action of $\CA_g$ on the set
$\GAz(\star_X,\Lloc)$ 
descends an action on the quotient set 
$\GAz(\star_X,\Lloc)/\la \bgp, \wgp \ra $.

\begin{theorem}
\lbl{thm.CAgcom}
The $\CA_g$ action on $\GAz(\star_X,\Lloc)$ descends to
an action on $\GAz(\star_X,\Lloc)/\la \bgp, \wgp \ra $.
\end{theorem}

\proof
It is an exercise in the definitions to see 
that if $(M,s_1), (M,s_2) \in \GAz(\star_X,\Lloc)$ and
and $(M,s_1) \stackrel{\bgp}{\sim} (M,s_2)$, then $\alpha_{ij}.(M,s_1)
\stackrel{\bgp}{\sim} \alpha_{ij}.(M,s_2)$.
However, this is not so clear for the case of the \wrapping\ move $\wgp$.
Consider pairs $(M,s_1) \stackrel{\wgp_k}{\sim} (M,s_2)$ by a wrapping move 
at the $k$th site. We need to show that $\alpha_{ij}.(M,s_1)
\stackrel{\bgp}{\sim} \alpha_{ij}.(M,s_2)$.

By definition, there exists an
$s \in \GA(\star_{X\cup \{\pt h\}},\Lloc)$ such that
\begin{eqnarray*}
s_1 &=& \conh (s) \\
s_2 &=& 
\conh \left( 
\phi_{t_k\to e^{-h} t_k e^h}(s)\sqcup 
\exp_\sqcup \left(-\frac{1}{2} \chiz(M^{-1}
\phi_{t_k\to e^{-h} t_k e^h} M)\right)\right).
\end{eqnarray*}
Then, we have that
\begin{align*}
\a_{ij} s_1 
&= \conh \left( \a_{ij} (s) \right)
\\
\a_{ij} s_2 
&= \conh
\left( \a_{ij}
\phi_{t_k\to e^{-h} t_k e^h}(s)\sqcup \left( 
\exp_\sqcup \left(-\frac{1}{2} \chiz(M^{-1}
\phi_{t_k\to e^{-h} t_k e^h} M) \right)\right)\right).
\end{align*}
If $k \neq i,j$ then $\a_{ij}$ commutes with 
$\phi_{t_k \mapsto e^{-h}t_k e^h}$,
thus the result follows. It remains to consider the cases of $k=i$ or $k=j$.

{\em The case of $k=i$.} 
Observe that $\a_{ij}$ equals to the substitution $\phi_{t_i \mapsto t_j^{-1}
t_i t_j}$, thus
$$
\a_{ij} \circ \phi_{t_i \mapsto e^{-h} t_i e^h}=
\phi_{t_i \mapsto e^{-h} t_j^{-1} t_i t_j e^h}
$$
We will first compare $\phi_{t_i \mapsto e^{-h} t_j^{-1} t_i t_j e^h}$
and $\phi_{t_i \mapsto t_j^{-1} e^{-h} t_i e^h t_j }$ using a Vertex 
Invariance (in short, VI)
Relation, and then $\phi_{t_i \mapsto t_j^{-1} e^{-h} t_i e^h t_j }$ with the 
identity using the wrapping move to reach the desired conclusion. 
Explicitly, the Vertex Invariance
 Relation of Figure \ref{relations3} implies that 
$$
\conh \left( 
\psdraw{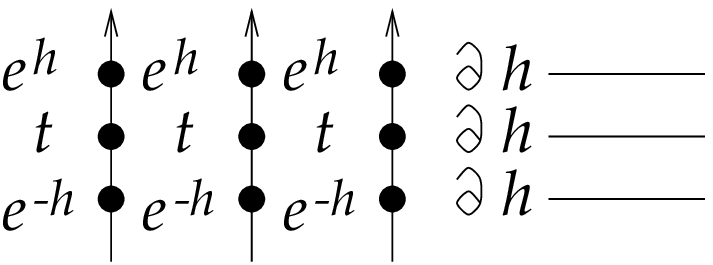}{1.5in} \right)=
\conh \left( 
\psdraw{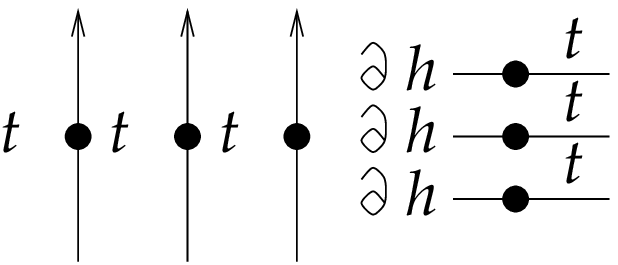}{1.3in} \right)
$$
if all appearances of $t$ and all legs $\pt h$ are drawn. Let $s(\pt h)$ denote
$s$ with the special leg colored explicitly by $\pt h$, and use the following
abbreviation
$$
\Psi(N)= 
\left( 
\exp_\sqcup \left(-\frac{1}{2} \chiz(N) \right)\right).
$$
The above discussion implies that:
\begin{align*}
\a&_{ij}(M,s_2)  \\ 
& =  \left(\a_{ij}M, \conh\left( 
\phi_{t_i \mapsto e^{-h} t_j^{-1} t_i t_j e^h}(s(\pt h)) \sqcup
\Psi(\a_{ij} M^{-1} .\phi_{t_i \mapsto e^{-h} t_j^{-1} t_i t_j e^h}(M)) 
\right)\right) \\
& \stackrel{\text{VI}}= 
\left(\a_{ij}M,
\conh \left(\phi_{t_i \mapsto t_j^{-1} e^{-h}  t_i e^h t_j} s(t_j \pt h) 
\sqcup
\Psi\left(\a_{ij}M^{-1}. \phi_{t_i \mapsto t_j^{-1} e^{-h}  t_i e^h t_j}M
\right) 
\right)\right) \\
& \stackrel{\wgp_i}\sim 
\left(\a_{ij}M,
\conh \left(\phi_{t_i \mapsto t_j^{-1}  t_i t_j} s(t_j \pt h) 
\sqcup 
\Psi\left(\a_{ij}M^{-1}. \phi_{t_i \mapsto t_j^{-1} t_i  t_j}M\right)
\right) \right) \\
& = 
\left(\a_{ij}M,
\left(\a_{ij} s(t_j \pt h)\left|_{\pt h \mapsto 0}
\right. 
\right)\right) \\
& = 
\left(\a_{ij}M,
\conh (\a_{ij} s)
\right) \\
& = 
\a_{ij} (M,s_1) .
\end{align*}

{\em The case of $k=j$.} 
Observe that 
$$
\a_{ij} \circ \phi_{t_j \mapsto e^{-h} t_j e^h}=
\phi_{t_j \mapsto e^{-h} t_j e^h, 
                    t_i \mapsto t_j^{-1} t_i t_j}
$$
First, we will compare the substitutions $\phiup {0}$ with $\phiup {4}$ 
(defined below) using a wrapping move $\wgp_j$, the following $\Psi$-identity
\begin{equation}
\lbl{eq.PsiIdentity}
\Psi(A^{-1}B) \sqcup \Psi(B^{-1}C)=\Psi(A^{-1}C)
\end{equation}
(which is a consequence of Proposition \ref{prop.w}) and Equation 
\eqref{identh} of Lemma \ref{lem.identities}.\nl
Then, we will compare $\phiup {4}$ with $\phiup {8}$ using the wrapping move
$\wgp_i$.\nl
Finally, we will compare $\phiup {8}$ with the identity using the case
of $k=i$ above.

Here, $\phiup {a}$ denote the substitutions:
\begin{align*}
\phiup {0} & = & t_j \mapsto e^{-h} t_j e^h, 
           &   & t_i \mapsto t_j^{-1} t_i t_j, \\
\phiup {2} & = & t_j \mapsto e^{-h} e^{h'} t_j e^{-h'} e^h, 
           &   & t_i \mapsto e^{h'} t_j^{-1} e^{-h'} t_i e^{h'} t_j e^{-h'}, \\
\phiup {3} & = & t_j \mapsto e^{h'} t_j e^{-h'}, 
           &   &    \\
\phiup {4} & = & t_j \mapsto t_j, 
           &   & t_i \mapsto e^{h} t_j^{-1} e^{-h} t_i e^{h} t_j e^{-h}, \\
\phiup {6} & = & t_j \mapsto t_j, 
           &   & t_i \mapsto e^{h} t_j^{-1} e^{-h} e^{h'} t_i e^{-h'} 
           e^{h} t_j e^{-h}, \\    
\phiup {7} & = & t_i \mapsto e^{h'} t_i e^{-h'}, 
           &   &  \\
\phiup {8} & = & t_j \mapsto t_j, 
           &   & t_i \mapsto e^{h} t_j^{-1} t_i t_j e^{-h}. 
\end{align*}
Then, we have that
\begin{align*}
\a_{ij}&(M,s_2)  \\ &=
\left( \a_{ij}M, \conh \left(\phiup {0} s(\pt h) \sqcup \Psi(\a_{ij} M^{-1} 
\phiup {0} M \right) 
\right)  \\
& \stackrel{\wgp_j}\sim 
 \left( \a_{ij}M, \con_{\{h,h'\}}
\left(\phiup {2} s(\pt h+\pt h') \sqcup \Psi(\phiup {3} \a_{ij} M^{-1} 
\phiup {2} M) \right.\right. \\ & \qquad \left.\left. \sqcup \Psi(\a_{ij} M^{-1} 
\phiup {3} \a_{ij} M) \right)  
\right) \\
& \stackrel{\text{\eqref{eq.PsiIdentity}}}= 
\left( \a_{ij}M, \con_{\{h,h'\}}
\left(\phiup {2} s(\pt h+\pt h') \sqcup \Psi(\a_{ij} M^{-1} 
\phiup {2} M) 
\right)\right) \\
& \stackrel{\text{\eqref{identh}}}= 
\left( \a_{ij}M, \conh
\left(\phiup {4} s(\pt h) \sqcup \Psi(\a_{ij} M^{-1} 
\phiup {4} M) 
\right)\right) \\ 
& \stackrel{\wgp_i}\sim 
\left( \a_{ij}M, \con_{\{h,h'\}}
\left(\phiup {6} s(\pt h+\pt h') \sqcup \Psi(\phiup {7} \a_{ij} M^{-1} 
\phiup {6} M) \right.\right. \\ & \qquad \left.\left.\sqcup \Psi(\a_{ij} M^{-1} 
\phiup {7} \a_{ij} M) \right)  
\right) \\
& \stackrel{\text{\eqref{eq.PsiIdentity}}}= 
\left( \a_{ij}M, \con_{\{h,h'\}}
\left(\phiup {6} s(\pt h+\pt h') \sqcup \Psi(\a_{ij} M^{-1} 
\phiup {6} M) 
\right)\right) \\
& \stackrel{\text{\eqref{identh}}}= 
\left( \a_{ij}M, \conh
\left(\phiup {8} s(\pt h) \sqcup \Psi(\a_{ij} M^{-1} 
\phiup {8} M) 
\right)\right) \\
& \stackrel{\text{$k=i$}}\sim 
\left( \a_{ij}M, \conh
\left(\a_{ij} s(\pt h) 
\right)\right) \\
& = 
\a_{ij}(M,s_1) .\tag*{\qed}
\end{align*}\eject

\subsection{$\Zrat$ commutes with the $\CA$ action}
\lbl{sub.Zcommutes1}

Recall the invariant $\Zrat$ of sliced crossed links of Section 
\ref{sec.tangles}.

\begin{theorem}
\lbl{thm.CAgcommute}
The map 
$$\s_X\circ\Zrat: \cD_X(D_g) \to \GA(\star_X,\La)/\la \wgp \ra$$ 
is $\CA_g$-equivariant, and induces a well-defined map 
$$
\s_X\circ\Zrat: \cD_X(D_g)/\la \CA_g \ra \to \GA(\star_X,\La)/\la \wgp,
\CA_g \ra.
$$
\end{theorem}

\begin{proof}
We will show that for $D\in \cD_X(D_g)$, we have
$$
\alpha_{ij}.\Zrat(D) = \Zrat(\alpha_{ij}.D).
$$
Given sliced crossed links $L_1,L_2$ such that $L_2=\a_{ij}.L_1$, consider
an accompanying link $L_{12}$ as in the previous section. $L_{12}$ 
may be presented so that
all the strands going through the $i$th meridianal disc proceed, parallel
to each other, through the $j$th meridianal disc, and then onto the $(i+1)$st
meridianal disc. That is, we can find a sequence of tangles, containing
some arc '$x$' (not exactly the arc '$a$' of Section \ref{sub.tangles}),
such that the link $L_{12}$ is presented by taking a parallel the arc '$x$'
according to some word (in the example at hand, $w_b$). 
Let us define $T_p,T_q,T_r,T_s,T_t,T_u$  by
$$
T_p=\setlength{\unitlength}{0.03\standardunitlength}
	\begin{array}{c}  \hspace{-1.7mm}
         	\raisebox{-8pt}{\input draws/SigmaTp.epc }
         	\hspace{-1.9mm}
	\end{array}
 \hspace{0.5cm}
T_q=\setlength{\unitlength}{0.03\standardunitlength}
	\begin{array}{c}  \hspace{-1.7mm}
         	\raisebox{-8pt}{\input draws/SigmaTq.epc }
         	\hspace{-1.9mm}
	\end{array}
 \hspace{0.5cm}
T_r\setlength{\unitlength}{0.03\standardunitlength}
	\begin{array}{c}  \hspace{-1.7mm}
         	\raisebox{-8pt}{\input draws/SigmaTr.epc }
         	\hspace{-1.9mm}
	\end{array}
 \hspace{0.5cm}
$$
and 
$$
T_s=\setlength{\unitlength}{0.03\standardunitlength}
	\begin{array}{c}  \hspace{-1.7mm}
         	\raisebox{-8pt}{\input draws/SigmaTs.epc }
         	\hspace{-1.9mm}
	\end{array}
 \hspace{0.5cm}
T_t=\setlength{\unitlength}{0.03\standardunitlength}
	\begin{array}{c}  \hspace{-1.7mm}
         	\raisebox{-8pt}{\input draws/SigmaTt.epc }
         	\hspace{-1.9mm}
	\end{array}
 \hspace{0.5cm}
T_u=\setlength{\unitlength}{0.03\standardunitlength}
	\begin{array}{c}  \hspace{-1.7mm}
         	\raisebox{-8pt}{\input draws/SigmaTu.epc }
         	\hspace{-1.9mm}
	\end{array}

$$
Continuing to work with the example of the above paragraph, 
then the sequence will look as follows
(where the $j$th, the $i$th and then the $(i+1)$st meridianal
discs are displayed):
$$
\{\dots, \Delta_x^{w_b}(T_p), \Delta_x^{w_b}(T_q),
\dots,
\Delta_x^{w_b}(T_r),
T_s,
\dots,
T_t,
\Delta_x^{w_b}(T_u),
\dots \}
$$
Now, the link $L_1$ (resp. $L_2$) 
is presented by the diagrams recovered by gluing
back in the $(i+1)$st (resp. $i$th) removed disc.

So, on the one hand we have:
\begin{eqnarray*}
\lefteqn{
\Zrat(L_1) } & & \\
& = & \dots \circ \Delta^{w_b}_x(\Z(T_p))
\circ (I_{w_d}\otimes G_{j,w_bw_e})
\circ \Delta^{w_b}_x(\Z(T_q))\circ
\dots \\
& & \dots
\circ \Delta^{w_b}_x(\Z(T_r))
\circ (I_{w_a}\otimes G_{i,w_b}) \circ
\Z(T_s)\circ \dots
\circ \Z(T_t)
\circ \Delta^{w_b}_x(\Z(T_u)) \circ \dots \\
& \stackrel{\wgp}{\sim} &
\dots \circ \Delta^{w_b}_x(\Z(T_p))
\circ (I_{w_d}\otimes G_{j,w_bw_e}) \circ (I_{w_d} \otimes G_{i,w_b} \otimes I_{w_e})
\circ \Delta^{w_b}_x(\Z(T_q))\circ
 \\
& & \dots
\circ \Delta^{w_b}_x(\Z(T_r)) \circ
\Z(T_s)\circ \dots
\circ \Z(T_t)
\circ \Delta^{w_b}_x(\Z(T_u)) \circ \dots 
\end{eqnarray*}
A $\wgp$-move has been used here to move the gluing word
$G_{i,w}$ up to (just underneath) the gluing word arising from
the $j$th meridional slice. This is legitimate 
because the meaning of the $\wgp$-move is that such
an element (if it involves {\it every} appearance of $t_i$) 
can be slid along a {\it paralleled} arc (here arising from
the component $x$).

On the other hand, we have:
\begin{align*}
&\hspace{-3mm}\Zrat(L_2) \\
& =  \dots \circ \Delta^{w_b}_x(\Z(T_p))
\circ (I_{w_d}\otimes G_{j,w_bw_e})
\circ \Delta^{w_b}_x(\Z(T_q))\circ
\dots \\
& \quad \dots
\circ \Delta^{w_b}_x(\Z(T_r))
\circ 
\Z(T_s)\circ \dots
\circ \Z(T_t)
\circ (I_{w_c\overline{w_{b}}}\otimes G_{i,w_b}) \circ
\Delta^{w_b}_x(\Z(T_u)) \circ \dots \\
& \stackrel{\wgp}{\sim} 
\dots \circ \Delta^{w_b}_x(\Z(T_p))
\circ (I_{w_d} \otimes G_{i,w_b} \otimes I_{w_e})
\circ (I_{w_d}\otimes G_{j,w_bw_e}) 
\circ \Delta^{w_b}_x(\Z(T_q))\circ
 \\
& \quad \dots
\circ \Delta^{w_b}_x(\Z(T_r)) \circ
\Z(T_s)\circ \dots
\circ \Z(T_t)
\circ \Delta^{w_b}_x(\Z(T_u)) \circ \dots 
\end{align*}
Observe that this line results from the line ending the previous 
paragraph of equations by the algebraic action of $\alpha_{ij}$,
as is required for the proof.
\end{proof}

\subsection{$\intrat$ commutes with the $\CA_g$ action}
\lbl{sub.CAg}

The next theorem says that the action of $\CA_g$ on diagrams commutes 
with $\intrat$-integration.

\begin{theorem}
\lbl{thm.CAgcommutes}
The map
$$ 
\intrat dX': \IPR_{X'} \GAz(\star_X,\Lloc) \to \GAz(\star_{X-X'},\Lloc)
$$
is $\CA_g$ equivariant and induces a well-defined map:
$$ 
\intrat dX': 
\IPR_{X'} \GAz(\star_X,\Lloc)/\la \CA \ra \to \GAz(\star_{X-X'},\Lloc)
/\la \CA \ra .
$$
\end{theorem}

\begin{proof}
Let the canonical decomposition of $s \in \IPR_{X'} \GA(\star_X,\Lloc)$ 
with respect to $X'$ be
$$
s = \es \left( \frac{1}{2} \sum_{kl} \gau{M_{kl}}{x_k}{x_l}{} \right) \sqcup R,
$$
for some matrix $M\in \Herm(\Lloc\to\BZ)$. Then,
the canonical decomposition of $\alpha_{ij}.s$ with respect to $X'$
is
$$
\alpha_{ij}.s = \es \left( \frac{1}{2} \sum_{kl} \gau{\alpha_{ij}.M_{kl}}{x_k}{x_l}{} \right) \sqcup \alpha_{ij}.R.
$$
Observe that $(\alpha_{ij}.M)^{-1}=\alpha_{ij}.M^{-1} \in \Herm(\Lloc\to\BZ)$
which implies that $\a_{ij}.s$ is $X'$-integrable. Furthermore, we have that

\begin{eqnarray*}
\intrat dX' \left( \alpha_{ij}.s \right) & = &
\intrat dX'\left( \es \left( \frac{1}{2} \sum_{kl} \gau{ \alpha_{ij}.M_{kl} }
{x_k}{x_l}{} \right) \sqcup (\alpha_{ij}.R) \right)  \\
& = &
\left<
\es \left( -\frac{1}{2} \sum_{kl} \gau{ \alpha_{ij}.(M_{kl}^{-1}) }
{x_k}{x_l}{} \right) , \alpha_{ij}.R\right> \\
& = &
\alpha_{ij}. \intrat dX' \left( s \right).
\end{eqnarray*}
Note that we could have used Theorem \ref{thm.varphi} to deduce the
above equality.
\end{proof}

We end this section with the following 

\begin{lemma}
\lbl{lem.CAg=1}
{\rm (a)}\qua For all $g$, the identity map $\A(\star_X,\Lloc)\to
\A(\star_X,\Lloc)$ maps the $\inn_g$-relation to the basing relation.

{\rm (b)}\qua For $g=1,2$, the identity map $\A(\star_X,\Lloc)\to
\A(\star_X,\Lloc)$ maps the $\CA_g$-relation to the basing relation.
\end{lemma}

\begin{proof}
Recall that $\inn_g$ is the subgroup of $\CA_g$ that consists of inner 
automorphisms of the free group $F$.

For the first part, an inner automorphism by an element $g \in F$
maps a bead $b \in \Lloc$ to $g^{-1}b g \in \Lloc$. The relations of Figure
\ref{relations2} cancel the $g$ at every trivalent vertex of a diagram. 
The remaining vertices change a diagram by a basing move $\bgp_2$.

The second part follows from the fact that $\inn_g=\CA_g$ for $g=1,2$,
see \cite{CS2,K}.
\end{proof}

\subsection{The definition of the rational invariant $\Zrat$ of $\pt$-links}
\lbl{sub.defineZb}

We can now define an invariant of boundary links as was promised
in Fact 10 of Section \ref{sub.maini}.

Definition \ref{def.ZratF}, together with
Theorems \ref{thm.surgeryb}, \ref{prop.tanglesbb}, \ref{thm.CAgcommute}
and \ref{thm.CAgcommutes} imply that:

\begin{definition}
\lbl{def.Zratb}
There is a map 
$$
\Zratgp: \text{$\pt$-links} \to \GAz(\Lloc)/\la \wgp,\CA \ra
$$
defined as follows: for a $\pt$-link $(M,L)$ of $g$ components,
choose a sliced crossed link $C$ in $\cD(D_g)$ whose components are in 1-1
correspondence with a set $X$, and define
$$
\Zratgp(M,L)=\frac{\intrat dX \left(\s \circ \Zratc(C) \right)}{
{c_+^{\sigma_+(B)} c_-^{\sigma_-(B)}}} \in \GAz(\Lloc)/\la \wgp, \CA \ra
$$
where $c_\pm=\int dU \Zc (S^3,U_{\pm})$ are some universal constants
of the unit-framed unknot $U_{\pm}$ and where $B$ is the linking matrix of
$C$ and $\s_{\pm}(B)$ is the number of positive (resp. negative) eigenvalues
of $B$.
\end{definition}

The above definition makes sense and uses the fact that as in Section 
\ref{sub.defineZ}, the set $\GAz(\Lloc)/\la \wgp \ra$ has a well-defined 
action of the group $\GA(\phi)$, which commutes with the $\CA$ action 
(where $\CA$ 
acts trivially on $\GA(\phi)$).

\subsection{Some more appearances of the $\CA$ group}
\lbl{sub.somemore}

Let us finish this section with some more information on the $\CA$ group,
which we will not need in our paper.
We wish to thank K. Vogtmann for bringing to our attention the references
in this section. 

The group $\CA_g$ (and its quotient $\CA_g/\mathrm{Inn}_g$)
appears in three different contexts: 

$\bullet$\qua
$\CA_g/\mathrm{Inn}_g$ is the group that acts on the set of $F$-links whose set
of orbits is identified with the set of boundary links. 

$\bullet$\qua
$\CA_g$ is the {\em group of motions} of a standard unlink in $\BR^3$,
as explained by Goldsmith following unpublished results of Dahm, 
\cite[Theorem 5.4]{Go}.
Let us explain a bit more. 
Using terminology from \cite{Go}, recall that a given a submanifold $N$ in
a noncompact manifold $M$, a {\em motion} of $N$ in $M$ is a one-parameter
family $f_t$ of diffeomorphisms of $M$ with compact support, for $0 \leq t
\leq 1$, such that $f_0$ and $f_1$ pointwise fix $N$. A motion is
{\em stationary} if it is homotopic to a motion that pointwise fixes $N$
at all times $t$. The set of equivalence classes of motions (modulo stationary 
ones) is a group. In the case of interest, $N$ is a standard oriented
unlink of ordered components in $\BR^3$ (resp. $S^3$), 
and the group of motions is identified with the group $\CA_g$ (resp.
$\CA_g/\mathrm{Inn}_g$). That is, $\CA_g$ is the group of motions of a unknotted
unlinked {\em string} in 3-space. In this context, see also 
Figure \ref{Drag}.

$\bullet$\qua
$\CA_g$ is the subgroup of the automorphism group of the free group
that sends every generator $t_i$ to a conjugate of itself. Geometric group 
theory tells us quite a bit about this algebraic group. McCool has given
a presentation in terms of the generators $\a_{ij}$, for
 $1 \leq i \neq j \leq g$
$$
\CA_g = \la \,\, \a_{ij} \,\, | [\a_{ij},\a_{kl}]=1, \,
[\a_{ik},\a_{jk}]=1, \, [\a_{ij},\a_{ik} \a_{jl}]=1 \, \ra
$$
where the indices occurring in each relation are assumed to be distinct.
Bounds for the cohomological dimension of $\CA_g$ are known, as well as normal
forms for the associated language of the elements in $\CA_g$, \cite{GuK}.

We will not use explicitly the presentation of $\CA_g$, however the reader 
may keep it in mind in the proof of \ref{thm.CAgcom}.
The reader should consult \cite{BrL,Go,GuK,Mc} for further information
on the $\CA$ group.

\section{A comparison between $\Z$ and $\Zrat$}
\lbl{sec.compaarhus}

Fact 10 of Section \ref{sub.maini} asks for an invariant $\Zrat$ which
is a rational form of the Kontsevich integral. 
The goal of this section is to show that the invariant of Definition
\ref{def.Zratb} is indeed a rational form of the Kontsevich integral
of a boundary link.

This will be obtained by comparing the construction of the Aarhus integral
$\int$ with that of $\intrat$.

\subsection{The $\hair$ map}
\lbl{sub.hair}

First, we need to define a $\hair$ map which replaces beads by the
exponential of hair.

Given generators $\{t_1,\dots, t_g\}$ of the free group, let 
$\{h_1,\dots,h_g\}$ be noncommuting variables such that $t_i=e^{h_i}$.

\begin{definition}
\lbl{def.hairmap}
Let
$$
\hair: \A(\star_X,\Lloc)\to\A(\star_{X \cup H})
$$ 
be the map that replaces each variable $t_i$ by an exponential of $h_i$-colored
hair, as in Figure \ref{attachlegs}.

Observe that the image of the $\hair$ map lies in the span of all 
diagrams that do not contain a $H$-labeled tree as a connected component,
and that it maps group-like elements to group-like elements.

We will extend the $\hair$ map on $\GAz(\star_X,\Lloc)$ (see Definition 
\ref{def.GAz}) by sending $(M,D)\in\Bla\times\GA(\star_X,\Lloc)$ to
$$
\hairnu(M,D)=\exp_\sqcup\left(-\frac{1}{2} \psdraw{circleb}{0.2in} \,\,
\chi(M) \right) \sqcup
\hair(D) \sqcup \nu(h_1) \dots \sqcup \nu(h_g),
$$
where $\nu(h)=Z(S^3,U)$ whose legs are colored by $h$ and 
$$
\chi(M)=\tr\logdet(M):=\sum_{n=1}^\infty\frac{(-1)^n}{n}\tr
\left((M (\e M)^{-1}-I)^n\right)
$$
where $\e:\La\to\BZ$. Note that the $\exp_\sqcup$ term in the above 
expression is a disjoint union of $H$-colored {\em wheels}. 

Finally, let $\hairnu$ be a slightly renormalized version $\hairnu$ of
the $\hair$ map defined by
$$
\hairnu(M,D)=\hair(M,D) \sqcup 
\nu(h_1)\sqcup \dots \sqcup \nu(h_g).
$$
\end{definition}

\begin{figure}[ht!]
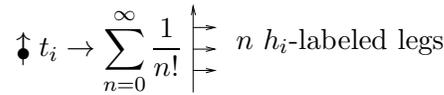

$$
\strutb{}{}{t_i} \to \sum_{n=0}^\infty \frac{1}{n!} \setlength{\unitlength}{0.03\standardunitlength}
	\begin{array}{c}  \hspace{-1.7mm}
         	\raisebox{-8pt}{\input draws/attacchn.epc }
         	\hspace{-1.9mm}
	\end{array}
 
$$
\caption{From beads to legs}\lbl{attachlegs}
\end{figure}

\begin{lemma}
\lbl{lem.hair}
There exist maps that fit in a commutative diagram
$$
\divide\dgARROWLENGTH by2
\begin{diagram}
\node{\GAz(\ostar_X,\Lloc)/\la \wgp  \ra}    
\arrow{s}\arrow[2]{e}
\node[2]{\A^0(\ostar_X,\Lloc)/\la \winf \ra}
\arrow{s}                          \\
\node{\GAz(\ostar_X,\Lloc)/\la \wgp, \CA_g \ra}  
\arrow[2]{e}\arrow{se,l}{\hair} 
\node[2]{\A^0(\ostar_X,\Lloc)/\la \winf, \CA_g \ra}
\arrow{sw,l}{\hair} \\
\node[2]{\A(\ostar_{X \cup H})} 
\end{diagram}
$$
\end{lemma}

\begin{proof}
Proposition \ref{prop.cbasing} implies that the horizontal maps are
 well-defined.
It is easy to show that the $X$-flavored group-like and infinitesimal
basing relations (defined in Section \ref{sec.binf}) 
are mapped by the $\hair$-map to $X$-flavored infinitesimal basing relations.
It is also easy to see that the $\wgp$ and $\CA_g$ relations are mapped to
$H$-flavored basing relations. 
\end{proof}

\subsection{A comparison between $\Zrat$ and $Z$}
\lbl{sub.equivariant1}

Consider the set $\cL_X(\O)$ of links $L \cup \O$ in $S^3$ where $\O$ is
an unlink whose components are in 1-1 correspondence with the set
$H=\{h_1,\dots,h_g\}$ and $L$ is a link whose components are in 1-1 
correspondence with a set $X$. Define a map
$$
\varphi_{\cup}: \cD_X(\O) \to \cL_X(\O)
$$
that maps the associated link $L$ of a sliced crossed diagram 
in $D_g$ to the link $L \cup \O \subset S^3$
according to a standard inclusion of $D_g \times I$ in $S^3$, by adding
a $g$ component unlink $\O$, one component for each gluing site of $D_g$ as 
follows:
$$
\varphi_{\cup} \left(
\{\dots,
\setlength{\unitlength}{0.03\standardunitlength}
	\begin{array}{c}  \hspace{-1.7mm}
         	\raisebox{-8pt}{\input draws/Sigma95.epc }
         	\hspace{-1.9mm}
	\end{array}
, \setlength{\unitlength}{0.03\standardunitlength}
	\begin{array}{c}  \hspace{-1.7mm}
         	\raisebox{-8pt}{\input draws/Sigma96.epc }
         	\hspace{-1.9mm}
	\end{array}
,\dots\}
\right)
\longrightarrow  \psdraw{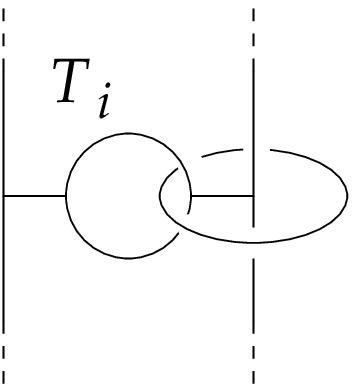}{0.5in}.
$$
Consider the Kontsevich integral 
$$
\s_{X \cup H}\circ\Z:\cL_X(\O)\to\A(\ostar_{X \cup H})
$$ 
which takes values in the completed space of unitrivalent graphs whose
legs are colored by $X \cup H$, modulo the $\AS, \IHX$ relations and the
{\em infinitesimal basing relations} discussed in Appendix \ref{sec.wversus}.
The purpose of this section is to show that 

\begin{theorem}
\lbl{thm.equiv1}
The following diagram commutes
$$
\begin{diagram}
\node{\NXO} 
\arrow{e,t}{\s\circ \Zrat} \arrow{s,l}{\varphi_{\cup}} 
\node{\GA(\ostar_X,\La)/\la \wgp \ra}
\arrow{s,r}{\hairnu}          \\
\node{\cL_X(\O)}
\arrow{e,t}{\s \circ \Z}
\node{\A(\ostar_{X \cup H})}
\end{diagram}
$$
\end{theorem}

The proof will utilize the following {\em magic formula} for the Kontsevich
integral of the Long Hopf Link, conjectured in \cite{A0} (in conjunction with
the so-called Wheels and Wheeling conjectures) and proven in \cite{BLT}.

\begin{theorem}{\rm\cite{BLT}}\lbl{BLT}
$$
\left(\s_{\{x,y\}}\circ\Z \right)\left(\psdraw{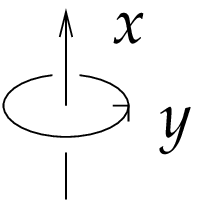}{0.35in} 
\right) = 
\strutb{x}{}{e^h}
\sqcup \nu(h)
\in \A(\tline_{\{x\}} \ostar_{\,\{y\}}). $$
\end{theorem}

{\bf Proof of Theorem \ref{thm.equiv1}}\qua
We will show first that the following diagram commutes: 
\begin{equation}
\lbl{eq.AXH}
\begin{diagram}
\node{\cD_X(\O)} 
\arrow{e,t}{\s\circ \Zrat} \arrow{s,l}{\varphi_{\cup}}
\node{\GA(\ostar_X,\La)}
\arrow{s,r}{\hairnu}          \\
\node{\cL_X(\O)}
\arrow{e,t}{\s \circ \Z}
\node{\A(\ostar_{X \cup H})}
\end{diagram}
\end{equation}
Given a word $w$ in the symbols $\uparrow$ and $\downarrow$, let
\begin{itemize}
\item
$\LongHopf_{j,w}$ denote the long Hopf link cabled according to the word $w$,
with the closed component labeled $h_j$.
For example: 
$$
\LongHopf_{j,\downarrow\uparrow\uparrow} = 
\psdraw{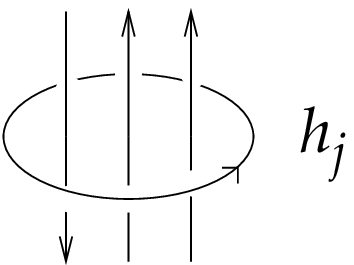}{0.5in}
$$
\item
$G^h_{j,w}$ denote the $j$th hairy gluing word corresponding to $w$.
This is definition by example: 
$$
G^h_{j,\downarrow\uparrow\uparrow} = 
\stbdown{}{}{e^{-h_j}} \stbup{}{}{e^{h_j}} \stbup{}{}{e^{h_j}} 
$$
\end{itemize}

We use the following corollary to Theorem \ref{BLT}.

\begin{corollary}
\lbl{nicecor}
$$
\left(\s_{\{h_j\}}\circ\Z \right)(\LongHopf_{j,w})
= G^h_{j,w} \sqcup \nu(h_j) \in \A(\uparrow_w \ostar_{\,\{h_j\}}).
$$
\end{corollary}

The commutativity of diagram \eqref{eq.AXH} now follows by the definition
of $\Zrat$.
Theorem \ref{thm.equiv1} follows from the commutativity of diagram 
\eqref{eq.AXH}, after considering the quotient $\NXO$ of $\cD_X(\O)$
(given by Proposition \ref{prop.tangles}) and 
symmetrizing (ie, composing with
$\s:\GA(\tline_X,\La)\to\GA(\ostar_X,\La)$) together with the fact that
$\hairnu$ commutes with symmetrization.\qed


\subsection{A comparison between $\intrat$ and $\int$}
\lbl{sub.equivariant2}

Recall that the Aarhus integral respects the group-like basing relations 
as shown in \cite[Part II, Proposition 5.6]{A} and also the infinitesimal
basing relations as was observed by D. Thurston, see \cite[Proposition 
2.2]{BL}. In other words, there is a well-defined map 
$$
\int dX':
\IPR_{X'} \A(\ostar_{X \cup H}) \to  \A(\ostar_{(X-X') \cup H}),
$$
where $H=\{h_1,\dots,h_g\}$. The purpose of this section is to show that

\begin{theorem}
\lbl{thm.equiv2} 
The following diagram commutes 
$$
\begin{diagram}
\node{\mathrm{Int}_{X'} \GAz(\ostar_X,\Lloc)} 
\arrow{e,t}{\intrat dX'}\arrow{s,l}{\hairnu} 
\node{\GAz(\ostar_{X-X'},\Lloc)}
\arrow{s,r}{\hairnu}  \\
\node{\mathrm{Int}_{X'}\A(\ostar_{X \cup H})}
\arrow{e,t}{\int dX'}
\node{\A(\ostar_{\,(X-X')\cup H})}
\end{diagram}
$$
\end{theorem}

\begin{proof}
It follows from Remark \ref{rem.apply} of Appendix \ref{sub.useful}.
\end{proof}





\subsection{$\Zrat$ is a rational form of the Kontsevich integral of a 
boundary link}
\lbl{sec.comparison}

Let us now formulate a main result of the paper:

\begin{theorem}
\lbl{thm.1}
There exists an invariant of $F$-links
$$
\Zrat: \text{$F$-links} \longto \A^0(\Lloc)/\la \winf \ra
$$
which 
\begin{itemize}
\item[\rm(a)]
descends to an invariant of boundary links
$$
\Zrat: \pt-\text{links} \longto \A^0(\Lloc)/\la \winf,\CA \ra
$$
\item[\rm(b)]
determines the Kontsevich integral $Z$ of a boundary link by
$$
Z=\hairnu \, \circ \Zrat.
$$
\item[\rm(c)]
determines the Blanchfield pairing of an $F$-link.
\end{itemize}
\end{theorem}

\begin{proof}
The invariant $\Zrat$ that satisfies (a) was defined in Section 
\ref{sub.defineZb}.

Let us show (b). 
We begin by observing that the following diagram commutes, by the definition
of the invariant $\Zratgp$, using the work of Section \ref{sec.CA}:
$$
\divide\dgARROWLENGTH by1
\begin{diagram}
\node{\text{$F$-links}}
\arrow{e,t}{\Zratgp}\arrow{s}
\node{\GAz(\Lloc)/\la \wgp \ra} 
\arrow{s} \\
\node{\text{$\pt$-links}}
\arrow{e,t}{\Zratgp}
\node{\GAz(\Lloc)/\la \wgp, \CA \ra}
\end{diagram}
$$
We now claim that the following diagram commutes:
\begin{equation}
\lbl{eq.diag}
\divide\dgARROWLENGTH by1
\begin{diagram}
\node{\text{$\pt$-links}}                
\arrow{e,t}{\Zratgp}\arrow{s,=} 
\node{\GAz(\Lloc)/\la \wgp, \CA \ra}
\arrow{s,r}{\hairnu}           \\
\node{\text{$\pt$-links}}              
\arrow{e,t}{\Z}
\node{\A(\ostar_H)}
\end{diagram}
\end{equation}
The above two diagrams, together with Lemma \ref{lem.hair} Theorem 
\ref{thm.1}.In order to show the commutativity of \eqref{eq.diag}, we will
assemble several commutative diagrams in one with two layers.
The top layer consists of the rational invariants $\Zrat$ and its
integration theory $\intrat$; the bottom layer consists of the 
LMO-Aarhus-Kontsevich integral $\Z$ and its integration theory $\int$. 
{\tiny
$$
\divide\dgARROWLENGTH by2
\begin{diagram}
\node[2]{\NXO}                
\arrow[2]{e,A} \arrow{s,-,..} \arrow{sw}
\node[2]{\NXO/\la \kappa_X, \CA \ra}       
\arrow[2]{s,..} \arrow{sw,r}{\Zratgp} \\
\node{\GA(\ostar_X,\La)/\la \wgp \ra}                 
\arrow[2]{e,t,3}{\intrat dX} \arrow[2]{s,..}
\node{}
\arrow{s,..}
\node{\GAz(\Lloc)/\la \wgp,\CA \ra} 
\arrow[2]{s,..} \\
\node[2]{\cL_X(\O)}         
\arrow{e,-} \arrow{sw}
\node{}
\arrow{e}
\node{\text{$\pt$-links}}   
\arrow{sw,r}{\Z}  \\   
\node{\A(\ostar_{X \cup H})}   
\arrow[2]{e,t}{\int dX}
\node[2]{\A(\ostar_H)} 
\end{diagram}
$$
}
Theorems \ref{thm.Zw}, \ref{thm.equiv1} and \ref{thm.equiv2} imply that
all faces of the above cube except the side right one 
commute. Since the map 
$\NXO\to\cL_X(\O)/\la \kappa,\CA \ra$ is onto, it follows that the remaining
face of the cube commutes.

Using the identification $\A(\tcircle_H) \cong \A(\ostar_H)$ (see 
Equation \eqref{eq.ostartcircle}) it follows that diagram
\eqref{eq.diag} commutes. This concludes the proof (b) of Theorem \ref{thm.1}.
(c) follows from \cite{GK2}; see also the discussion of Section 
\ref{sub.surgery}.
\end{proof}

\subsection{Comments on Theorem \ref{thm.1}}
\lbl{sub.comments.thm1}

Let us make some complimentary remarks on Theorem \ref{thm.1}.

\begin{remark}
\lbl{rem.empty}
First a remark about the empty set. For the unlink $\O$ of $g$ components, 
we have that 
$\Zrat(\O)=(\text{empty matrix}, 1) \in \A^0(\Lloc)/\la \winf, \CA \ra$.
On the other hand, the Kontsevich intregal of $\O$ is given by $\Z(\O)=
\nu(h_1) \sqcup \dots \sqcup \nu(h_g)$ (where $\nu(h_i)$ is the Kontsevich
integral of the unknot, whose diagrams are colored by $h_i$), and
$$
\hairnu(\text{empty matrix}, 1)=\nu(h_1) \sqcup \dots \sqcup \nu(h_g).
$$
This is how $\Zrat(\O)$ determines $Z(\O)$, as stated in the above theorem.
\end{remark}

\begin{remark}
\lbl{rem.augmentation}
Erasing a component $L_i$ of a boundary link corresponds to substituting
$t_i=1$. In other words,
$$
\Zrat(M,L-L_i)=\Zrat(M,L)|_{t_i=1}.
$$
\end{remark}

\begin{remark}
\lbl{rem.ww}
The $\Zrat$ invariant of $F$-links takes values in a  quotient
of a set of pairs of matrices and vectors $\A^0(\Lloc)$ by an equivalence
relation. 
If we restrict $\Zrat|_{\mathcal S}$ to a class $\mathcal S$ of $F$-links with 
the same Blanchfield pairing (or the same $S$-equivalence class), then
$\Zrat|_{\mathcal S}$ is equivalent to an invariant
$$ \mathcal S \longrightarrow \A(\Lloc)/\lp \w \rp
$$
that takes values in a {\em graded vector space}, where the degree of a graph
is the number of trivalent vertices. In addition, the edges of the trivalent
graphs are decorated by elements of a subring of $\Lloc$ that depends on 
$\mathcal S$. 
In particular, for acyclic $F$-links (ie, for links whose free cover is 
acyclic, or equivalently, whose Blanchfield pairing is trivial)
$\Zrat$ is equivalent to an invariant with values in the more manageable graded
vector space $\A(\La)/\lp \winf,\CA \rp$; see Section \ref{sub.defineZ}.
\end{remark}

\begin{proof}
Recall that the matrix part $W$ of $\Zrat$ is an element of $\B(\Lloc\to\BZ)$
and determines the $S$-equivalence class of the $F$-link, as follows from
\cite{GK2}. Recall also that the wrapping relations 
$\wgp$ and $\winf$ depend only on $W$. Thus, it follows that $S$-equivalent 
$F$-links
have equal matrix-part, and that the matrix-free part of $\Zrat$ makes
sense as an element of the graded vector space $\A(\Lloc)/\lp \winf \rp$.

\cite[Theorem 3]{GK2} implies that $(M,L)$ is 
a good boundary link, (ie, $H_1(X^{\w},\BZ)=0$) iff its associated matrix
$W$ is invertible over $\La$. In that case, the construction of the $\Zrat$
invariant implies that $\Zrat(M,L) \in 
\A^0(\La)/\la \winf, \CA \ra$. 
\end{proof}

\begin{remark}
\lbl{rem.g=1}
For $g=1$, the sets of 1-component $F$-links and $\pt$-links coincide with 
the set of 
knots in \ihs s, and the relations $\winf$ and $\CA_1$ are trivial, 
see Lemmas \ref{lem.wg=1} and \ref{lem.CAg=1}. In this case, $\Zrat$ takes
values in $\B(\Lloc\to\BZ) \times \A(\Lloc)$, and the matrix-free part
of $\Zrat$ takes values in the {\em graded algebra} $\A(\Lloc)$, where
the multiplication in $\A(\Lloc)$ is given by the disjoint union of graphs.
It turns out that the matrix-free part of $\Zrat$ is a group-like element
in $\A(\Lloc)$, ie, it equals to the exponential of a series of 
{\em connected} graphs.
  
For $g=2$ (ie, for two component $\pt$-links) the relation $\CA_g$ is 
trivial, see Lemma \ref{lem.CAg=1} and the set of 2-component $F$-links 
coincides with the one of $\pt$-links.
\end{remark}

\begin{remark}
\lbl{rem.compare}
In a sense, $\Z$ is a graph-valued $\Lhat$-invariant of $\pt$-links, where
$\Lhat$ is the completion of $\La \otimes \BQ$ 
with respect to the augmentation ideal
whereas $\Zrat$ is a graph-valued $\Lloc$-invariant. It often
happens that $\Lhat$-invariants can be lifted to $\Lloc$-invariants,
see for example the discussion of Farber-Ranicki, \cite{FRa}, on a
problem of Novikov about circle-valued Morse
theory. It would be nice if an
extension of this circle-valued Morse theory to infinite dimensions 
(using the circle-valued Chern-Simons action on the space of gauge equivalence
classes of connections on a principal bundle over a 3-manifold) would
lead to an independent construction of the $\Zrat$-invariant.
\end{remark}

\begin{remark}
\lbl{rem.equivariant}
As will become apparent from the construction of the $\Zrat$ invariant, it
is closely related to equivariant linking matrices of nullhomotopic links
in the complement of an $F$-link, see Proposition \ref{prop.later}. It is
an interesting question to ask whether the $\Zrat$ invariant of an $F$-link
can be defined in terms of an equivariant construction on the free cover
of the $F$-link.
\end{remark}

\section{A universal property of $\Zrat$}
\lbl{sec.universal}

\subsection{Statement of the result}
\lbl{sub.stateuniversal}

In this section we will show that the rational invariant $\Zrat$, evaluated
on acyclic $F$-links, is characterized by a universal property.
Namely, it is the universal \fti \ of acyclic $F$-links with respect to
the null filtration, \cite{Ga}.

Using the terminology of \cite{Ga}, a {\em null-move} on an $F$-link 
$(M,L,\phi)$ can be described by surgery on a clover $G$ in $M-L$ whose leaves 
lie in the kerner of $\phi$, \cite{Ga}. The result $(M,L,\phi)_G$
of surgery is another $F$-link. As usual, one can define a notion of \fti s 
on the set of $F$-links using null moves; we will call the corresponding
filtration the {\em Euler filtration}. The number of invariants of degree 
$n$ can be identified with a quotient vector space $\Gl {n}$. 
For simplicity, we will confine ourselves to the case of
{\em acyclic $F$-links}, that is $F$-links whose free cover $X^\omega$
satisfies $\bar H_\star(X^\omega,\BZ)=0$. In \cite{Ga}, it was shown that 
there is an onto map
$$
\A(\La)/\la \winf \ra \longto \Gl {}
$$
which preserves degrees. The next theorem constructs an inverse to this
map and characterizes $\Zrat$ by a universal property.

\begin{theorem}
\lbl{thm.KOK}
$\Zrat$ is a universal $\A(\La)/\la \winf \ra$-valued 
\fti \ of acyclic $F$-links with respect to the Euler filtration. 
In particular, we have that 
$\A(\La)
\cong \Gl {}$, over $\BQ$.
\end{theorem}

The arguments in this section are minor modifications of arguments given in
\cite{GR} for the case of knots. Since the work of \cite{GR} predated the
existence of the rational form $\Zrat$, and since we prefer to be 
self-sufficient,  we will repeat the arguments here.

\subsection{A relative version of the $\Zrat$ invariant}
\lbl{sub.relativeZrat}

We already discussed in the introduction (before the statement of Theorem
\ref{thm.KOK}) the notion of \fti s of $F$-links, based on the null move.
The null move is described in terms of surgery on a clasper whose leaves
lie in the kernel of the basing map $\phi$ of an $F$-link.

Our task here is to give a {\em relative version} of the $\Zrat$ invariant
for null links in the complement of an $F$-link.
Given an $F$-link $(M,L,\phi)$, where $\phi:\pi_1(M\sminus L)\to F$, 
we will call a 
link $C \in M\sminus L$ {\em null relative to $(M,L,\phi)$} if each component 
of $C$ is mapped to $1$ under $\phi$.
Notice that surgery under a null link preserves the underlying $F$-structure. 
In what follows, we will omit $\phi$ from the notation, for typographical
reasons.
 
Given a null link $C$ relative to an $F$-link $(M,L)$, choose a sliced crossed
link $L' \cup C' \in \cD(D_g)$ such that after surgery on $L' \subset 
S^3\sminus\O$ 
it gives rise to the $F$-link $(M,L)$ with $C' \subset S^3\sminus\O$ 
transformed to
a null link $C$ relative to $(M,L)$. We define
\begin{eqnarray*}
\Zratgp((M,L),C) &=&
\frac{\intrat dX( \s \circ \Zratc(L'\cup C'))}{(\mathrm{normalization})} 
\in \GAz(\ostar_Y,\Lloc) /\la \wgp \ra \\
\Zratgp((M,L)_C) &=&
\frac{\intrat dY( \Zrat((M,L),C))}{(\mathrm{normalization})} 
\in \GAz(\ostar,\Lloc)/\la \wgp \ra
\end{eqnarray*}
where $X$ (resp. $Y$) is in 1-1 correspondence with the components of $L$ 
(resp. $C$), and where $(\mathrm{normalization})$ refers to factors inserted
to deal with the Kirby move of adding a distant union of unit-framed unlinks,
as in Section \ref{sub.defineZ}.
Theorem \ref{thm.1} has the following relative version (as is revealed
by its proof in Section \ref{sub.defineZ}):
 
\begin{theorem} 
\lbl{thm.RC}
There exists an invariant 
$$
\Zratgp: \text{$F$-links} \cup \mathrm{null} \,\, 
C \to \GAz(\ostar_C,\Lloc)/\la \winf \ra
$$
which is local (ie, compatible with respect to considering sublinks of C)
and fits in a commutative diagram
$$
\divide\dgARROWLENGTH by2
\begin{diagram}
\node{\text{$F$-links} \cup \mathrm{null} \,\, C}
\arrow{e,t}{\Zratgp} \arrow{s,l}{C-\mathrm{surgery}}
\node{\GAz(\ostar_Y,\Lloc)/\la \wgp \ra}
\arrow{s,r}{\intrat dY} \\
\node{\text{$F$-links}} 
\arrow{e,t}{\Zratgp}
\node{\GAz(\Lloc)/\la \wgp \ra}
\end{diagram}
$$
\end{theorem}

As in Section \ref{sub.defineZ}, we denote by $\Zrat$ the composition of
$\Zratgp$ with the map 
$$
\GAz(\ostar_Y,\Lloc)/\la \wgp \ra
\to \A^0(\ostar_Y,\Lloc)/\la \winf \ra.
$$
Given a null link $C$ relative to $(M,L)$, one can define an {\em equivariant
linking matrix} $\lk_{(M,L)}(C)$, by considering the linking numbers of the
components of the lift $\ti C$ of $C$ in the free cover of $M-L$. 
$\lk_{(M,L)}(C)$ is a Hermitian matrix over $\Lloc$, well defined up to 
conjugation by a diagonal matrix with elements in $F$. The
details of the definition, together with the following proposition will
appear in a subsequent publication, \cite{Ga}:

\begin{proposition}
\lbl{prop.later}
For $C$ and $(M,L)$ as above, the covariance matrix $\cov$ of $\Zrat$ satisfies
$$\cov(\Zrat((M,L),C))=\lk_{(M,L)}(C).$$
\end{proposition}

The discussion in this section simplifies considerably when the $F$-links
in question are acyclic. In this case, there is no need to localize, and
the wrapping relation becomes trivial. In other words, we get an invariant
$$
\Zrat: \text{Acyclic} \, \text{$F$-links} \cup \text{null links} 
\longto \A(\ostar_Y,\La)/\la \winf \ra.
$$

\subsection{Surgery on claspers}
\lbl{sub.surgeryc}

Let us review in brief some elementary facts about surgery on claspers.
For a more detailed discussion, we refer the reader to \cite{GGP}
and also \cite[Section 3]{GR}.

Consider a null clasper $G$ of degree $1$ in $S^3-\O$ (a so-called \ygraph ). 
Surgery on 
$G$ can be described by surgery on a six component link $E \cup L$ associated 
to $G$, where $E$ (resp. $L$) is the three component link that consists of 
the edges (resp. leaves) of $G$. 

\begin{lemma}
\lbl{lem.lkG}
The equivariant linking matrix of $E\cup L$
and its negative inverse are given as follows:
\begin{equation}
\lbl{eq.arms}
\left(\begin{array}{cc}
0 & I \\ I & \lk_{\SO}(L_i,L_j) \\
\end{array}\right) \hspace{0.3cm}
\text{ and } \hspace{0.3cm}
\left(\begin{array}{cc}
\lk_{\SO}(L_i,L_j) & -I \\ -I & 0 \\
\end{array}\right) .
\end{equation}
\end{lemma}

\begin{proof}
The $0$ block follows from the fact that $\{E_i,E_j,\O\}$ is an unlink.
The $I$ block follows from the fact that $E_i$ is a meridian of $L_i$,
using the formula for the Kontsevich integral of the Long Hopf Link, together
with the fact that the linking number between $L_i$ and $\O$ is zero.
\end{proof}

The six component link $E \cup L$ is partitioned in three blocks
of two component links $A_i=\{E_i,L_i\}$ each for $i=1,2,3$, the {\em arms}
of $G$. A key feature of surgery on $G$ is the fact that surgery on any 
proper subset of the set of arms does not alter $M$. In other words,
alternating $(S^3,\O)$ with respect to surgery on $G$ equals to alternating 
$(S^3,\O)$ with 
respect to surgery on all nine subsets of the set of arms 
$A=\{A_1,A_2,A_3 \}$. That is,
\begin{equation}
\lbl{eq.altZ}
\Zrat([(S^3,\O),G])=\Zrat([(S^3,\O), A]).
\end{equation}
Due to the locality property of the $\Zrat$ invariant (ie, Theorem 
\ref{thm.RC}) the nontrivial 
contributions to the right hand side of Equation \eqref{eq.altZ} come from
the {\em strutless part} $Z^{\text{rat},t}(\SO,A_1 \cup A_2 \cup A_3)$
of $Z^{\text{rat},t}(\SO,A_1 \cup A_2 \cup A_3)$ that consists of graphs with
legs {\em touch} (ie, are colored by) all three arms of $G$.
The above discussion generalizes to the case of an arbitrary 
disjoint union of claspers.

\subsection{Counting above the critical degree}
\lbl{sub.countabove}

\begin{proposition}
\lbl{prop.MKu1}
The Euler degree $n$ part of $\Zrat$ is a type $n$ invariant 
of acyclic $F$-links with values in $\A_n(\La)$.
\end{proposition}

\begin{proof}
Since $\A_{\text{odd}}(\La)=0$, it suffices to consider the case of 
even $n$. Suppose that $G=\{G_1, \dots, G_{m} \}$ 
(for $m \geq 2n+1$) is a collection of null claspers in $S^3-\O$ each of 
degree $1$, and let $A$ denote the set of arms of $G$.
Equation \eqref{eq.altZ} and its following discussion implies that
$$
\Zrat([\SO,G])=\Zrat([\SO, A])
$$
and that the nonzero contribution to the right hand side come from
diagrams in $\Zratt(\SO,A)$ that touch all arms. Thus, contributing
diagrams have at least $3(2n+1)+1=6n+4$ $A$-colored legs, to be glued 
pairwise.

Notice that the diagrams in $\Zratt(\SO,A)$ contain no struts. Thus, at most 
three $A$-colored
legs meet at a vertex, and after gluing the $A$-colored legs we obtain
trivalent graphs with at least $(6n+4)/3=2n+4/3$ trivalent vertices, in other
words of Euler degree at least $2n+2$. Thus,
$\Zrat_{2n}([\SO,G])=0 \in \A_{2n}(\La)$, which implies that
$\Zrat_{2n}$ is a invariant of acyclic $F$-links of type $2n$ with values in 
$\A_{2n}(\La)$.
\end{proof} 

Sometimes the above vanishing statement is called {\em counting above
the critical degree}.

\subsection{Counting on the critical degree}
\lbl{sub.counton}

Our next statement can be considered as {\em counting on the critical degree}. We need a preliminary definition. 

\begin{definition}
\lbl{def.cc}
Consider a null clasper $G$ in $S^3-\O$ of degree $2n$, and let 
$G^{break}=\{G_1,\dots, G_{2n}\}$ 
denote the collection of degree $1$ claspers $G_i$ which are obtained
by inserting a Hopf link in the edges of $G$. Let 
$G^{nl}=\{G^{nl}_1,\dots, G^{nl}_{2n}\}$ denote the collection of 
abstract unitrivalent graph obtained by removing the leaves of the 
$G_i$ (and leaving one leg, or univalent vertex, for each leave behind).
Choose arcs from a fixed base point to the trivalent vertex of each 
$G^{nl}_i$, which allows us to define the equivariant linking numbers
of the leaves of $G^{break}$. Then the {\em complete contraction} 
$\la G\ra\in\A(\La)$ of $G$ is defined to be the sum over all ways of gluing 
pairwise the legs of $G^{nl}$, and placing as a bead the
equivariant linking number of the corresponding leaves.
\end{definition}

The result of a complete contraction of a null clasper $G$ is a 
well-defined element of $\A(\La)$. Changing the arcs is taken care by the
Vertex Invariance relations in $\A(\La)$. 
The next proposition computes the {\em symbol} of $\Zrat$:

\begin{proposition}
\lbl{prop.MKu2}
If $G$ is a null clasper of degree $2n$ in $S^3\sminus\O$, then
$$
\Zrat_{2n}([\SO,G])=\la G \ra \in \A_{2n}(\La)/\la \winf \ra 
$$
\end{proposition}
 
\begin{proof}
It suffices to consider a collection
$G=\{G_1, \dots, G_{2n} \}$ of claspers in $S^3$ each of 
degree $1$. Let $A$ denote the set of arms of $G$. The counting argument
of the above Proposition shows that the contributions
to $\Zrat_{2n}([\SO,G])=\Zrat_{2n}([\SO, A])$ come from complete contractions
of a disjoint union  $D=Y_1 \cup \dots \cup Y_{2n}$ of $2n$ {\em vortices}. 
A vortex is the diagram $\sfY$, the next simplest unitrivalent graph after
the strut. Furthermore, the $6n$ legs of $D$ should touch all
$6n$ arms of $G$. In other words, there is a 1-1 correspondence between
the legs of such $D$ and the arms of $G$.

Consider a leg $l$ of $D$ that touches an arm $A_l=\{E_l,L_l\}$ of $G$.
If $l$ touches $L_l$, then due to the restriction of the negative inverse
linking matrix of $G$ (see Lemma \ref{lem.lkG}), it needs to be contracted
to another leg of $D$ that touches $E_l$. But this is impossible, since
the legs of $D$ are in 1-1 correspondence with the arms of $G$.

Thus, each leg of $D$ touches presicely one edge of $G$. In particular, each 
leg of $D$ is colored by three edges of $G$. 

Consider a vortex colored by three edges of $G$ as part of the $\Zrat$
invariant. Since the three edges of $G$ is an unlink in a ball disjoint from
$\O$, it follows that the beads on the edges of the vortex are $1$, if
the three edges belong to the same $G_i$, and $0$ otherwise.

Thus, the diagrams $D$ that contribute are a disjoint union of $2n$ vortices 
$Y=\{Y_1,\dots,Y_{2n}\}$ and these votrices are in 1-1 correspondence with
the set of claspers $\{G_1,\dots,G_{2n}\}$, in such a way that the
legs of each vortex $Y_i$ are colored by the edges of a unique clasper $G_j$.

After we glue the legs of such $Y$ using the negative inverse linking matrix 
of $G$, the result follows.
\end{proof}

Let us mention that the discussion of Propositions \ref{prop.MKu1} and 
\ref{prop.MKu2} really applied to the unnormalized rational Aarhus integral; 
however since we are counting above the critical
degree, we need only use the degree $0$ part of the normalization which
equals to $1$; in other words we can forget about the normalization.

The above proposition is useful in realization properties of the $Z_{2n}$
invariant, but also in proving the following {\em Universal Property}:

\begin{proposition}
\lbl{prop.MKu3}
For all $n$, 
the composite map of \cite{Ga} with that of Proposition \ref{prop.MKu2}
$$
\A_{2n}(\La)/\la \winf \ra
\longto \Gl {2n}  \stackrel{\Zrat_{2n}}\longto \A_{2n}(\La)/\la \winf \ra
$$
is the identity. Since the map on the left is onto, it follows that
the above maps are isomorphism. 
\end{proposition}

\begin{remark}
In the above propositions \ref{prop.MKu1}-\ref{prop.MKu3}, $S^3$ can be 
replaced by any \ihs\ (or even a \qhs s) $M$.  \qed
\end{remark}

Clearly, Propositions \ref{prop.MKu1}-\ref{prop.MKu3} imply Theorem
\ref{thm.KOK}.

\section{Relations of the $\Zrat$ invariant with Homology \newline Surgery}
\lbl{sec.HSurgery}

In our paper, we constructed an invariant $\Zrat$ of $F$-links. The 
construction of this invariant leads naturally to the noncommutative 
localization of $\La=\BZ[F]$; see Fact 8 in Section \ref{sub.maini}.
Thus, our paper leads to a natural relation (from the point of view of
Qunatum Topology) between $F$-links and noncommutative localization.

Two previous relations between $F$-links and noncommutative localization are 
 known through the work of Cappell-Shaneson and Farber-Vogel;
see \cite{CS1,CS2,Vo,FV}. Let us sidcuss those briefly, pointing out that
this subject, well-known among senior topologists, is not as well-known
to quantum topologists.

Cappell-Shaneson reduced the problems of {\em Homology Surgery} to the 
computation of $Ga$-groups; \cite{CS1,CS2}. 
A typical problem of Homology Surgery the
following: given a knot $K$ in $S^3$, consider the manifold $M=S^3_{K,0}$ 
obtained by $0$-surgery on $K$. Of course it is true that $H_\star(M,\BZ)
=H_\star(S^2 \times S^1, \BZ)$, but the problem is to decide when is the case
that the equivariant homology of $M$ and $S^2 \times S^1$ coincides:
$H_\star(M,\La)=H_\star(S^2 \times S^1, \La)$ for $\La=\BZ[Z]$.

The $\Ga$ groups $\Ga(\La\to\BZ)$ 
of Cappell-Shaneson (interpreted in terms of matrices, a la
Levine) ask for a classification of cobordism of Hermitian matrices over 
$\La$ which are invertible over $\BZ$.

These $\Ga$ groups are hard to compute, but Cappell-Shaneson reduced 
the problem of classification of $F$-links modulo cobordism (in high enough
dimensions) to the computation of their $\Ga$-groups.

Vogel identified these $\Ga$-groups in terms of Wall surgery $L$-groups of a 
suitable localization of $\La$, \cite{Vo}. Later on, Vogel and Farber 
identified this localization with the ring of rational functions in 
noncommuting variables, \cite{FV}. 

This explains the two previously known relations among $F$-links and
the noncommutative localization.

Later on, Farber constructed an invariant of $F$-links with values in 
$\Lloc$, \cite{Fa}. This invariant was reinterpreted in terms Seifert
surfaces (and the surgery view of $F$-links) in \cite{GL2}.

This reinterpretation of Farber's invariant, together with the results
of \cite{GK2} imply that the ``matrix part'' of our $\Zrat$-invariant
of $F$-links determines the Blanchfield pairing of them.

One last comment: if one tries to develop a theory of finite type
invariants for degree $1$ maps, then one is naturally lead to the notion
of beads, which lie in $\BZ[\pi]$ for the fundamental group $\pi$ of the
target manifold. This was discovered in joint work of Levine and the
first author; see \cite{GL1}. Presumably, one could
try to construct an invariant for degree $1$ maps in the spirit of the
$\Zrat$ invariant, which would take values in a completed space of 
trivalent graphs with beads in the noncommutative localization
of $\BZ[\pi]$. We will develop this in a later publication.

\section*{Appendices}

\appendix\small

\section{A brief review of diagrammatic calculus}
\lbl{sec.diag}

\subsection{Some useful identities}
\lbl{sub.identities}

In this section we establish some notations and collect some lemmas.
Let $X$ denote some finite set. Given a sum of diagrams $s$ with legs
colored by a set $X$, and a subset $X'=\{x_1,\dots,x_r\}$ of $X$, we
will often denote $s$ by $s(x_1,\dots,x_r)$ or simply by $s(x)$. This notation
allows us to extend linearly to diagrams whose legs are labeled by formal
linear combinations of elements of $X$.
For example, for $a,b \in \BQ$, $v,w \in X$,
$s(\dots,x_{i-1},av + bw,x_{i+1},\dots )$ denotes the element obtained
from $s$ by replacing each diagram appearing in $s$ by the sum of all
diagrams obtained by relabeling each leg labeled $x_i$ by either $v$ (with
a multiplicative factor of $a$) or $w$ (with a multiplicative factor
of $b$). Given diagrams $A(x), B(x)$, we let
$$
\left\la A(\pt x), B(x) \right\ra_X \,\, \text{ and } \,\,
A(\pt x) \, \flat_X B(x)
$$
denote respectively the sum of all ways of pairing all $X$-colored legs of 
$A$ with all (resp. some) $X$-colored legs of $B$.

Furthermore, for some $q\in \Lloc$, $s(\dots,x_{i-1},x_i q,x_{i+1},
\dots)$ denotes the element obtained by pushing a bead labeled by $q$ onto
each leg labeled $x_i$ in the sense specified by the following diagram:
$$
\strutb{}{x_i}{1} \to \strutb{}{x_i}{q}
$$
Taken together, these conventions 
give meaning to such expressions as $s(xM)$, for $M\in \Herm(\Lloc\to\BZ)$.

The following lemma collects a number of useful identities, which the
reader would be well-advised to understand.

\begin{lemma}
\lbl{lem.identities}
{\rm (a)}\qua
Let $A\in \A(\star_{X},\Lloc), B\in \A(\star_{X\cup Y},\Lloc)$ (for $Y$ a 
bijective copy
of $X$), and let $X'$ and $X''$ denote bijective copies of $X$.
\begin{equation}
\lbl{identa}
\la A(\pt x), B(x,y)|_{x=y} \ra_{X} = \la A(\pt x' + \pt x''), 
B(x',x'')\ra_{X' \cup X''}
\end{equation}

{\rm (b)}\qua
Let $A_1,A_2,B\in \A(\star_{X},\Lloc)$. 
\begin{equation}
\lbl{identb}
\la A_1(\pt x) A_2(\pt x) , B(x)\ra_X 
= \la A_2(\pt X), A_1(\pt x) \,\flat_{X}\,B(x)\ra_X.
\end{equation}

{\rm (c)}\qua
Let $A,B\in\A(\star_X,\Lloc),$ and let $Y$ be
a bijective copy of $X$.
\begin{equation}
\lbl{identc}
A(\pt x)\,\flat_{X}\,B(x) = \la  A(\pt y), B(x+y) \ra_{Y} 
\end{equation}

{\rm (d)}\qua
Let $A=\es(a) \in \A(\star_X,\Lloc)$ and $B=\es(b) \in \A(\star_X,\Lloc)$. 
Then 
\begin{equation}
\lbl{identd}
\es(a)\, \flat_{X'}\, \es(b) = \es(c),
\end{equation}
where $c$ is the sum of all {\em connected} diagrams that arise by pairing
{\em all} legs labeled from $X'$ of some diagram (not necessarily connected)
appearing in $A$ to {\em some} legs of 
some diagram (not necessarily connected) appearing in 
$B$.

{\rm (e)}\qua
Let $A \in \A(\star_X,\Lloc)$, let $Y$ be a bijective copy of $X$,
and let $M\in \Herm(\Lloc\to\BZ)$.
\begin{equation}
\lbl{idente}
\left(\es\left( \gau{M_{ij}}{\pt x_j}{y_i}  \right)\right)\,
\flat_{X}\, A(x) = A(yM).
\end{equation}

{\rm (f)}\qua
In particular, 
\begin{eqnarray}
\lbl{identf}
\left\la \exp\left( \st{\pt x}{y} \right), B(x) \right\ra_x
& = & B(y) \\
\lbl{identg}
\exp\left( \st{\pt x}{y} \right) \flat_x \exp(B(x)) 
& = & \exp\left( \sum \text{ \parbox{1.8in}{{\rm all ways of replacing some 
of the $x$-legs of B by $y$-legs}}} \right) 
\end{eqnarray}

{\rm (g)}\qua
In addition, if $s(\pt z) \in \A(\star_{X \cup \{h', h'',\pt z\}},
\Lloc)$ then
\begin{eqnarray}
\lbl{identh}
\con_{\{h',h''\}}(s(\pt h' + \pt h'')) = \conh(s(\pt h)|_{h', h''\mapsto h}) .
\end{eqnarray}
\end{lemma}

\subsection{Iterated Integration}
\lbl{sub.inteq1}

This section proves an Iterated Integration property of the $\intrat$ integral. 
Fix labeling sets $X''\subset X' \subset X$.

\begin{lemma}[Iterated integration]
\lbl{itintlem}
Given $s\in \IPR_{X'}\A(\star_X,\Lloc) \cap \IPR_{X''}\A(\star_X,\Lloc)$, 
then,
\begin{eqnarray*}
\intrat dX''(s) & \in & \IPR_{X'-X''}\A(\star_{X-X''},\Lloc) \\
\intrat d(X'-X'') \left(
\intrat dX''(s)\right) &=& \intrat dX' (s).
\end{eqnarray*}
\end{lemma}

\begin{proof}
Denote the components of $X'$ and $X''$ by $a_i$ and $b_j$ respectively. 
Consider the covariance matrix of $s$ with respect to $X'$, presented with
respect to this basis:
$$
M=
\left[
\begin{array}{ll}
A & C \\ C^\star & B
\end{array}
\right],
$$
where $A=A^\star$ and $B=B^\star$, where here and below $A^\star$ denotes
the transpose of $C$, followed by the conjugation $t_i \to t_i^{-1}$.
If the canonical decomposition of $s$ with respect to $X'$ is then,
$$
\es\left(\frac{1}{2} \gau{A_{ij}}{a_j}{a_i}  
       + \frac{1}{2}\gau{B_{ij}}{b_j}{b_i}  
       + \gau{C_{ij}}{b_j}{a_i}  \right)\sqcup R,
$$
\noindent
(suppressing the summations) then the canonical decomposition of $s$ with 
respect to $X''$ is
$$
\es\left(\frac{1}{2} 
\gau{A_{ij}}{a_j}{a_i} 
\right)\sqcup \left(
\es\left(
\frac{1}{2}
\gau{B_{ij}}{b_j}{b_i}  + \gau{C_{ij}}{b_j}{a_i} \right) 
\sqcup R \right).
$$
Using the identities of Lemma \ref{lem.identities}, we have
\begin{eqnarray*}
\lefteqn{\intrat  dX''(s) } & & \\
& = & 
\left\la 
\es\left( -\frac{1}{2} \gau{A^{-1}_{ij}}{\pt a_j}{\pt a_i}{\ } \right),
\es\left( \frac{1}{2} 
\gau{B_{ij}}{b_j}{b_i}  + \gau{C_{ij}}{b_j}{a_i} 
\right) \sqcup R \right\ra_{X''} \\
& = &
\left\la 
\left.
\left(
\es\left(
\gau{C_{ij}}{b_j}{\pt a_i} 
\right)
\flat_{X'}\
\es
\left(
-\frac{1}{2}
\gau{A_{ij}^{-1}}{a_j}{a_i}  \right)
\right)
\right|_{a_i \rightarrow \pt a_i}, \right. \\ & & \left.
\es
\left(
\frac{1}{2}
\gau{B_{ij}}{b_j}{b_i} 
\right)\sqcup R\right\ra_{X''} \\
& = &
\es\left(
\frac{1}{2}
\gau{(B-C^\star A^{-1}C)_{ij}}{b_j}{b_i} 
\right)
\sqcup   \\  & & 
\left(
\left\la 
\es\left(
-\frac{1}{2}
\gau{A_{ij}^{-1}}{\pt a_j}{\pt a_i}{\ }
-\gau{(A^{-1}C)_{ij}}{b_j}{\pt a_i}  \right),R \right\ra_{X''}
\right)
\end{eqnarray*}
\noindent
It is easy to see that if $M$ and $A$ are invertible over $\BZ$,
then so is $B-C^\star A^{-1}C$ over $\BZ$, see 
\cite[Part II, Proposition 2.13]{A}. Thus, the first part of the theorem 
follows. Continuing,
\begin{eqnarray*}
\lefteqn{
\intrat d(X''-X') \left(
\intrat dX' ( s ) \right)} & & \\
& = &
\left\la 
\es
\left(
-\frac{1}{2}
\gau{(B-C^\star A^{-1}C)^{-1}_{ij}}{\pt b_j}{\pt b_i} 
\right),
 \right.  \\  & & \quad
\left.
\left\la 
\es
\left(
-\frac{1}{2}
\gau{A_{ij}^{-1}}{\pt a_j}{\pt a_i}{\ }
-
\gau{(A^{-1}C)_{ij}}{b_j}{\pt a_i}  \right),
R\right\ra_{X''}\right\ra_{X'-X''} \\
& = &
\left\la 
\left.
\left(
\es
\left(
-\gau{(A^{-1}C)_{ij}}{\pt b_j}{\pt a_i} 
\right)
\flat_{X'-X''}\
\es
\left(
-\frac{1}{2}
\gau{(B-C^\star A^{-1}C)^{-1}_{ij}}{b_j}{b_i} 
\right)
\right)
\right|_{b_i \rightarrow \pt b_i}\right. 
\left. \right. \\
& & \quad \left.
\es \left(
-\frac{1}{2}
\gau{A_{ij}^{-1}}{\pt a_j}{\pt a_i}{\ }
\right)
, R \right\ra_{X'} \\
& = &
\left\la 
\es
\left(
-\frac{1}{2}
\gau{(A^{-1}+ A^{-1}C(B-C^\star A^{-1}C)^{-1}C^\star 
A^{-1})_{ij}}{\pt a_j}{\pt a_i} 
\right. \right. \\
& & \quad \left. \left.
+
\gau{(A^{-1}C(B-C^\star A^{-1}C)^{-1})_{ij}}{\pt b_j}{\pt a_i}  
-\frac{1}{2}
\gau{(B-C^\star A^{-1} C)^{-1}_{ij}}{\pt b_j}{\pt b_i} \right)
, R\right\ra_{X'}. 
\end{eqnarray*}
The identification of the term on the left of the pairing with 
$\es\left( -\frac{1}{2} \gau{M^{-1}_{ij}}{\pt x_j}{\pt x_i}{\ } \right)$ 
completes the proof.
\end{proof}

\begin{remark}
\lbl{rem.iterated}
With the notation of the previous lemma, the matrices $A \oplus
(B-C^\star A^{-1} C)$ and $M$ are congruent by a product of
elementary matrices,
as follows from the following identity
$$
\left[
\begin{array}{ll}
A & 0 \\ 0 & B-C^\star A^{-1} C
\end{array}
\right] = 
\left[
\begin{array}{ll}
I & -AC^\star \\ 0 & I
\end{array}
\right]^\star
\left[
\begin{array}{ll}
A & C \\ C^\star & B
\end{array}
\right]
\left[
\begin{array}{ll}
I & -AC^\star \\ 0 & I
\end{array}
\right] .
$$
and from the fact that the conjugating matrix is a product of elementary 
matrices.
It follows from this, that the Iterated Integration 
Lemma \ref{itintlem} is valid on $\GAz(\star_X,\Lloc)$.
\end{remark}

\subsection{Integration by parts}
\lbl{sub.inteq2}

Fix labeling sets $X$ and $X'\subset X$. Recall from Section 
\ref{sub.inteq} the notion of $X'$-substantial diagrams, that is diagrams
that do not contain a strut labeled from $X'\cup \pt X'$.
The {\em divergence} of an $X'$-substantial diagram (extended by linearity
to sums of diagrams) with respect to $X'\subset X$ is defined to be the sum 
of diagrams obtained by pairing {\em all} legs labeled from $\pt_{X'}$ with 
{\em some} legs labeled from $X'$.

\begin{lemma}[Integration by Parts]
\lbl{intbypartslem}
Let $s$ be an $X'$-substantial diagram and 
$t\in \IPR_{X'}\A(\star_X,\Lloc)$. Then
\begin{eqnarray*}
s\,\flat_X t & \in & \IPR_{X'}\A(\star_X,\Lloc) \\
\intrat dX'( s\,\flat_X  t ) &=&
\intrat dX' (\div_{X'}(s)\,\flat_{X-X'}\, t).
\end{eqnarray*}
\end{lemma}

\begin{proof}
For the first statement, 
assume, for convenience, that $X'=\{f_i\}$ and $X-X'=\{g_j\}$, and rename
$X'$ and $X-X'$ by $F$ and $G$ respectively.
Applying Lemma \ref{lem.identities}.(\ref{identc}),
introducing a labeling set $\tilde{F}$ (resp. $\tilde{G}$), bijective with
$F$ (resp. $G$)
(and associated sets $\pt_{\tilde{F}}$ and $\pt_{\tilde{G}}$),
gives the expression in question as:
$$s(\oll{f},\oll{g},\oll{\pt f},\oll{\pt g})\,
\flat_{F\cup G}\, t(\oll{f},\oll{g}) = 
\left\la s(\oll{f},\oll{g},\oll{\pt \tilde{f}},\oll{\pt \tilde{g}})
, t(\oll{f}+\oll{\tilde{f}},\oll{g}+\oll{\tilde{g}})\right
\ra_{\tilde{F}\cup \tilde{G}}
$$
If the canonical decomposition of $t$ with respect to $F$ is
$$
t(f,g) = \es\left( \frac{1}{2}\gau{M_{ij}}{f_j}{f_i}  \right) \sqcup 
T(\oll{f},\oll{g}),
$$
\noindent
then the previous expression is equal to
$$
\es \left( \frac{1}{2} \gau{M_{ij}}{f_j}{f_i}  \right) \sqcup 
\left\la 
s(\oll{f},\oll{g},\oll{\pt \tilde{f}},\oll{\pt \tilde{g}}),
\es
\left(
\gau{M_{ij}}{\tilde{f}_j}{f_i} +\frac{1}{2}\gau{M_{ij}}{\tilde{f}_j}{\tilde{f}_i} \right)
\sqcup T(\oll{f}+\oll{\tilde{f}},\oll{g}+\oll{\tilde{g}})\right\ra_{\tilde{F}\cup \tilde{G}}.
$$
Now, because $s$ is an $X'$-substantial operator, this gives the
required decomposition of $s\,\flat_{F\cup G}\, t$ with respect to $F$.

For the second statement, we compute
\begin{eqnarray*}
\lefteqn{\intrat dF (s\,\flat_{F\cup G}\, t)} \\ 
& = &
\left\la  
\es \left( -\frac{1}{2} \gau{M^{-1}_{ij}}{\pt f_j}{\pt f_i}{\ } 
\right),\right.  \\
& & 
\left.\left\la 
s(\oll{f},\oll{g},\oll{\pt \tilde{f}},\oll{\pt \tilde{g}}),
\es
\left(
\gau{M_{ij}}{\tilde{f}_j}{f_i} +\frac{1}{2}\gau{M_{ij}}{\tilde{f}_j}{\tilde{f}_i} \right)
\sqcup T(\oll{f}+\oll{\tilde{f}},\oll{g}+\oll{\tilde{g}})\right\ra_{\tilde{F}\cup \tilde{G}}\right\ra_F. 
\end{eqnarray*}
To proceed, we must adjust this calculation so that it can distinguish
struts, appearing in the left of the pairing, according to which factors,
appearing in the right of the pairing, they are glued to. This adjustment
is made with Lemma \ref{lem.identities}.(\ref{identa}), and leads to
the following
(introducing sets $F^a$, $F^b$ and $F^c$, bijective copies of $F$):
\begin{eqnarray*}
\lefteqn{
\left\la  
\es \left( 
-\frac{1}{2} \gau{M^{-1}_{ij}}{\pt f_j^a}{\pt f_i^a}{\ \ }
-\frac{1}{2} \gau{M^{-1}_{ij}}{\pt f_j^b}{\pt f_i^b}{\ \ }
-\frac{1}{2} \gau{M^{-1}_{ij}}{\pt f_j^c}{\pt f_i^c}{\ \ } \right. \right. 
} & &
\\
& & 
\left.\left. -\gau{M^{-1}_{ij}}{\pt f_j^a}{\pt f_i^b}{\ \ }
-\gau{M^{-1}_{ij}}{\pt f_j^a}{\pt f_i^c}{\ \ }
-\gau{M^{-1}_{ij}}{\pt f_j^b}{\pt f_i^c}{\ \ }
 \right),\right. \\
& & \\
& &
\left.\left\la 
s(\oll{f^a},\oll{g},\oll{\pt \tilde{f}},\oll{\pt \tilde{g}}),
\es
\left(
\gau{M_{ij}}{\tilde{f}_j}{f_i^b} +\frac{1}{2}\gau{M_{ij}}{\tilde{f}_j}{\tilde{f}_i} \right)
\right.\right. \\
& & \left.\left.
\sqcup T(\oll{f^c}+\oll{\tilde{f}},\oll{g}+\oll{\tilde{g}})\right\ra_{\tilde{F}\cup \tilde{G}}\right\ra_{F^a\cup F^b
\cup F^c}.  
\end{eqnarray*}
The calculation proceeds by removing factors appearing in 
left hand side of the pairing, at the expense of operating with them on the
right hand side, follow Lemma \ref{lem.identities}.(\ref{identb}).

The first term to be treated in this way will be
$
\es\left( -\frac{1}{2} \gau{M_{ij}^{-1}}{\pt f_j^b}{\pt f_i^b}{\ \ }
\right),$ giving (following an application of Lemma 
\ref{lem.identities}.(\ref{identd})):
$$
\es\left( -\frac{1}{2} \gau{M_{ij}^{-1}}{\pt f_j^b}{\pt f_i^b}{\ \ }
\right)\flat_{F^b}
\es\left(
\gau{M_{ij}}{\tilde{f}_j}{f_i^b}{\ } \right) = 
\es\left(
\gau{M_{ij}}{\tilde{f}_j}{f_i^b}{\ } \right) \sqcup
\es\left(
-\frac{1}{2}\gau{M_{ij}}{\tilde{f}_j}{\tilde{f}_i}{\ } \right).
$$
The next term to be so treated will be
$\es\left(
-\gau{M^{-1}_{ij}}{\pt f_j^c}{\pt f_i^b}{\ } \right),\
$
which gives:
\begin{eqnarray*}
\lefteqn{
 \es\left(
-\gau{M^{-1}_{ij}}{\pt f_j^c}{\pt f_i^b}{\ } \right)
\flat_{F^b\cup F^c}
\left(
\es\left(
\gau{M_{ij}}{\tilde{f}_j}{f_i^b}{\ } \right) 
\sqcup
T( \oll{f^c}+\oll{\tilde{f}},\oll{g}+\oll{\tilde{g}} ) \right) 
} & & \\
& &=
\es\left(
\gau{M_{ij}}{\tilde{f}_j}{f_i^b}{\ } \right) \sqcup
T( (\oll{f^c}-\oll{\tilde{f}})+\oll{\tilde{f}},\oll{g}+\oll{\tilde{g}} ). 
\end{eqnarray*}
The result of these operations is
\begin{eqnarray*}
\lefteqn{ \left\la  
\es \left( 
-\frac{1}{2} \gau{M^{-1}_{ij}}{\pt f_j^a}{\pt f_i^a}{\ \ }
-\frac{1}{2} \gau{M^{-1}_{ij}}{\pt f_j^c}{\pt f_i^c}{\ \ }
-\gau{M^{-1}_{ij}}{\pt f_j^a}{\pt f_i^b}{\ \ }
-\gau{M^{-1}_{ij}}{\pt f_j^a}{\pt f_i^c}{\ \ }
 \right),\right. } \\
& & \\
& &
\hspace{1cm}\left.\left\la 
s(\oll{f^a},\oll{g},\oll{\pt \tilde{f}},\oll{\pt \tilde{g}}),
\es
\left(
\gau{M_{ij}}{\tilde{f}_j}{f_i^b} \right)
\sqcup T(\oll{f^c},\oll{g}+\oll{\tilde{g}})\right\ra_{\tilde{F}\cup \tilde{G}}\right\ra_{F^a\cup F^b
\cup F^c}. \\
& & 
\end{eqnarray*}
Performing the pairing with respect to the variables
$\tilde{F}$, leads to the replacement of the the inner 
$\la \cdot, \cdot \ra $ with
$$
\left\la 
s(\oll{f^a},\oll{g},\oll{f^b}M,\oll{\pt \tilde{g}}),
T(\oll{f^c},\oll{g}+\oll{\tilde{g}})\right\ra_{\tilde{G}}.
$$
Trading the factor 
$
\es\left(
-\gau{M_{ij}^{-1}}{\pt f_j^b}{\pt f_i^a}{\ \ }
\right),
$  gives:
\begin{eqnarray*}
\lefteqn{\es\left(
-\gau{M_{ij}^{-1}}{\pt f_j^b}{\pt f_i^a}{\ \ }
\right) \flat_{F^a\cup F^b}
s(\oll{f^a},\oll{g},\oll{f^b}M,\oll{\pt \tilde{g}}) } & & \\
 &=&
\left.
\div_{F^a}\left( s(\oll{f^a},\oll{g},\oll{f^bM + \pt f^a},\oll{\pt g}) 
\right)\right|_{\oll{\pt g}\rightarrow \oll{\pt \tilde{g}}}. 
\end{eqnarray*}
Performing the pairing with respect to $F^b$ sets $f^b$ to zero, 
and thus we arrive at the following expression.
\begin{eqnarray*}
\lefteqn{
\left\la  
\es \left( 
-\frac{1}{2} \gau{M^{-1}_{ij}}{\pt f_j^a}{\pt f_i^a}{\ \ }
-\frac{1}{2} \gau{M^{-1}_{ij}}{\pt f_j^c}{\pt f_i^c}{\ \ }
-\gau{M^{-1}_{ij}}{\pt f_j^a}{\pt f_i^c}{\ \ }
 \right),  \right. } & & \\ & & \left.
\left\la 
\left.\div_{F}\left(s(\oll{f^a},\oll{g},\oll{\pt f^a},\oll{\pt g})
\right)
\right|_{\oll{\pt g}\rightarrow \oll{\pt \tilde{g}}},
T(\oll{f^c},\oll{g}+\oll{\tilde{g}})\right\ra_{\tilde{G}}
\right\ra_{F^a \cup F^c}
.
\end{eqnarray*}
The proof follows by inspection.
\end{proof}

\section{The Wheels Identity}
\lbl{sec.wheels}

The goal of this section is to show the following {\em Wheels Identity}:

\begin{theorem}
\lbl{thm.wheels}
For every $M \in \Herm(\Lloc\to\BZ)$, we have that
$$
\int  dX \left( \es\left( \frac{1}{2} \sum
\strutb{x_j}{x_i}{M_{ij}(e^{h_1},\dots,e^{h_g})} \right)
\right) =
\exp_{\sqcup} \left( -\frac{1}{2} \psdraw{circleb}{0.2in} \chi(M) \right),
$$
where we think of $\chi(M)$ as taking values in $\Lhat/(\mathrm{cyclic})$,
which is identified via the map $t_i\to e^{h_i}$
with the set $\BQ\la\la h_1,\dots,h_g \ra\ra/(
\mathrm{cyclic})$ of formal power series in noncommuting variables modulo
cyclic permutations.
\end{theorem}

\begin{proof}
Choose generators $\{t_1,\dots,t_g\}$ for the free group $F$ and identify
$\Lhat$ with the ring of noncommutative power series $\{h_1,\dots,h_n\}$.
Recall the map $\La \to \Lhat$ given by $t_i\to e^{h_i}$.

The Aarhus integral is calculated by splitting the quadratic part from the
rest of the integrand as follows:
\begin{eqnarray*}
\lefteqn{ \int  dX \left( \es\left( \frac{1}{2} \sum
\gau{M_{ij}(e^{h_1},\dots,e^{h_g})}{x_j}{x_i}  \right)
\right) } \\
& = & 
\int  dX \left(
\es\left( \frac{1}{2} \sum
\gau{M_{ij}(1,\dots,1)}{x_j}{x_i}  \right) \right. 
\\ & & \quad \left. 
\sqcup\, 
\es\left( \frac{1}{2} \sum \left(
\gau{M_{ij}(e^{h_1},\dots,e^{h_g})}{x_j}{x_i}  
-
\gau{M_{ij}(1,\dots,1)}{x_j}{x_i}  
\right)
\right)
\right) \\
& = &
\left\la  
\es\left( -\frac{1}{2} \sum
\gau{M_{ij}(1,\dots,1)^{-1}}{\pt x_j}{\pt x_i}  \right), \right.
\\ & & \quad  
\left.
\es\left( \frac{1}{2} \sum \left(
\gau{M_{ij}(e^{h_1},\dots,e^{h_g})}{x_j}{x_i}  
-
\gau{M_{ij}(1,\dots,1)}{x_j}{x_i}  
\right)
\right)
\right\ra_X.
\end{eqnarray*}
With the following notation for ``Red'' and ``Green'' struts,
\begin{eqnarray*}
\gau{R}{}{} & = & \sum \gau{M_{ij}(1,\dots,1)^{-1}}{\pt x_j}{\pt x_i}{}, \\
\gau{G}{}{} & = & \sum
\left(
\gau{M_{ij}(e^{h_1},\dots,e^{h_g})}{x_j}{x_i} 
-
\gau{M_{ij}(1,\dots,1)}{x_j}{x_i} 
\right).
\end{eqnarray*}
it follows that the desired integral equals to the following formal series
(which makes sense, since finitely many terms appear in each degree)
\begin{eqnarray*}
\lefteqn{
\left\la  \es\left( -\frac{1}{2} 
\gau{R}{}{} \right),
\es\left( \frac{1}{2} 
\gau{G}{}{} \right) \right\ra  } & & \\
&=&
\sum_{n=0}^{\infty}(-1)^n \frac{1}{2^n} \frac{1}{2^n} \frac{1}{n!}
\frac{1}{n!} 
\left\la  
\left(
\strutb{}{}{R} \right)^n,
\left(
\strutb{}{}{G} \right)^n \right\ra .
\end{eqnarray*}
Specifically, the answer is recovered 
if each ``join'' is accompanied by a Kronecker delta. For example:
$$
\begin{array}{c}
\strutb{}{}{R} \vspace{-0.4cm} \\ 
\strutb{}{}{G}
\end{array}
=
\sum_{i,j,k}
\begin{array}{l}
\strutb{x_k}{}{M_{jk}(1,\dots,1)^{-1}} \vspace{-0.4cm} \\ 
\strutb{}{x_i}{M_{ij}(e^{h_1},\dots,e^{h_g})}
\end{array}
-
\sum_{i,j,k}
\begin{array}{l}
\strutb{x_k}{}{M_{jk}(1,\dots,1)^{-1}} \vspace{-0.4cm} \\ 
\strutb{}{x_i}{M_{ij}(1,\dots,1)}
\end{array} .
$$
We focus on a particular summand. It will
help in the management of the combinatorics if we label each
of the arcs that appear in the struts:
$$
\left\la  
\left(
\strutb{}{}{R} \right)^n,
\left(
\strutb{}{}{G} \right)^n \right\ra  \rightarrow 
\left\la 
\strutb{r_{1,1}}{r_{1,2}}{}
 \dots\
\strutb{r_{n,1}}{r_{n,2}}{}
,
\strutb{g_{1,1}}{g_{1,2}}{}
\dots\
\strutb{g_{n,1}}{g_{n,2}}{}
\right\ra 
$$
Now we introduce some definitions to assist in the enumeration of the 
polygons that arise from this pairing.

\begin{itemize}
\item An {\em RG-shape of degree $n$} is a collection of even-sided
polygons whose edges are alternately colored red and green. The 
total number of edges is $2n$.
\item
A {\em labeled RG-shape of degree $n$} is an RG-shape with an extra bivalent
vertex in each edge, such that:
\begin{enumerate}
\item{the resulting set of red (resp. green) edges is equipped with 
a bijection to \newline $\{r_{1,1},r_{1,2},\dots,r_{n,1},r_{n,2}\}$ (resp.
$\{g_{1,1},g_{1,2},\dots,g_{n,1},g_{n,2}\}$ ),}
\item{if an edge is labeled, for example $r_{i,1}$ (resp. $g_{i,1}$), 
then its neighboring red (resp. green) edge is labeled $r_{i,2}$
(resp. $g_{i,2}$).}
\end{enumerate}
\item{Let the constant $C_\tau$, for $\tau$ some RG-shape, denote
the number of labeled RG-shapes which recover $\tau$ when their labels
and extra bivalent vertices are forgotten.}
\end{itemize}

Observe that:
\begin{eqnarray*}
\left\la  
\left(
\gau{R}{}{} \right)^n,
\left(
\gau{G}{}{} \right)^n \right\ra    &=&   
\sum_{\mbox{$\tau$ a labeled RG-shape of deg $n$}} \tau \\ &=& 
\sum_{\mbox{$\tau$ an RG-shape of deg $n$}}
C_\tau \tau. 
\end{eqnarray*}
To calculate this constant
describe an RG-shape of degree $n$, $\tau$, by a sequence\break 
$(i_1,\dots,i_n)$
where $i_j$ counts the number of polygons with $2j$ edges. 

\begin{lemma}
$$
C_\tau =
\frac{ n! 2^n n! 2^n}{ (2.1)^{i_1} i_1!(2.2)^{i_2}i_2!
\dots
(2.n)^{i_n}i_n!}.
$$
\end{lemma}

\proof
Fix a copy of $\tau$ with its set of edges ordered (any old how, so as
to distinguish them from each other). There are precisely $2^n n! 2^n n!$
ways of labeling this so as to respect the conditions of the above definition.

But this overcounts by a multiplicative factor of 
the order of the symmetry group of $\tau$ (that is, with the ordering
of the edges forgotten). That symmetry group is precisely the
following group (that is, swaps of identical polygons, and rotations
and flips of individual polygons):
$$
(\S_{i_1} \ltimes (D_1 \times \dots D_1) )
\times \dots \times (\S_{i_n} \ltimes (D_n \times \dots D_n) ),$$
where $D_n$ represents the $n$th dihedral group.

Thus we conclude:
$$
C_\tau = \frac{ 2^n n! 2^n n! }
{\#\{ (\S_{i_1}\ltimes 
(D_1 \times \dots \times D_1)) \times \dots \times
(\S_{i_n}\ltimes 
(D_n \times \dots \times D_n)) \}}\eqno{\qed}
$$

If $\tau_i$ denotes a connected polygon with $2i$ edges, it follows that
\begin{eqnarray*}
\lefteqn{ \left\la  \es\left( -\frac{1}{2} 
\gau{R}{}{} \right),
\es\left( \frac{1}{2} 
\gau{G}{}{} \right) \right\ra } & & \\ &=& 
\sum_{n=0}^{\infty}(-1)^n \frac{1}{2^n} \frac{1}{2^n} \frac{1}{n!}
\frac{1}{n!} 
\left\la  
\left(
\gau{R}{}{} \right)^n,
\left(
\gau{G}{}{} \right)^n \right\ra , \\
& = &
 \sum_n \sum_{i_1,\dots,i_n} \left( \prod_{p=1}^{n}
\left( \frac{(-1)^p}{(2p)} \right)^{i_p} \frac{1}{i_p!} \right)
\tau_1^{i_1}\sqcup \dots \sqcup \tau_n^{i_n}, \\
& = & 
\es\left(
\frac{1}{2} \sum_{p=1}^{\infty}\frac{(-1)^p}{p} \tau_p \right).
\end{eqnarray*}
Substituting the definitions of the ``Red'' and ``Green'' struts into 
$\tau_p$ yields
\begin{eqnarray*}
\tau_p &=& \tr\left( (M(1,\dots,1)^{-1}(M(e^{h_1},\dots,e^{h_g})
- M(1,\dots,1))^p \right) \\ &=&
\tr\left( (M(1,\dots,1)^{-1}M(e^{h_1},\dots,e^{h_g}) -I)^p \right).
\end{eqnarray*}
from which the result follows.
\end{proof}

\section{Infinitesimal versus group-like wrapping relations}
\lbl{sec.wversus}

In this section we show that the infinitesimal wrapping relations $\winf$
imply the group-like wrapping relation $\wgp$. In order to achieve this,
it will be useful to introduce the variant $\winfp_i$ of $\winf_i$
for $i=1,\dots,g$ 
where we allow the diagram $D \in \A(\star_{X \cup \{\pt h\}},\Lloc)$ 
to have {\em several} legs labeled by $\pt h$. In other words, $\winfp_i$ is
generated by the move
$
(M,D_1)\sim(M,D_2)
$
where
$$
D_m =  \conh \left( 
\varphi_{m,i} (D)\sqcup \left( 
\exp_\sqcup \left(-\frac{1}{2} \chiz(M^{-1}
\varphi_{m,i} M) \right)\right)\right)
$$
for $m=1,2$ and $\varphi_{1,i}$ is the identity and $\varphi_{2,i}$ (
abbreviated by $\varphi_i$ below) denotes 
the substitution $\phi_{t_i\rightarrow e^ht_ie^{-h}}$ on beads.

It is obvious that the $\winfp_i$ relation includes the $\winf_i$ relation.
Our goal for the rest of the section is to show that they coincide
(see Proposition \ref{prop.wequal}) and that the $\wgp$-relations are
implied by the $\winf$-relations (see Proposition \ref{prop.winf3}).

Recall the maps $\eta_i: \Lloc\to\Lnew$ of Section \ref{sub.Cohn} that keep
track of the $h$-degree $1$ part of the substitutions
$\varphi_{t_i\to e^{-h}t_ie^h}(a)$ of $a \in \Lloc$. We will introduce an 
extension $\heta_i$ of $\eta_i$
which will allow us to keep track of the $h$-degree $n$ part of 
$\varphi_{t_i\to e^{-h}t_ie^h}(a)$, for all $n$. The following lemma (which
the reader should compare with Lemma \ref{lem.winf1}) summarizes the properties
of $\heta_i$.

\begin{lemma}
\lbl{lem.winf11} For every $i=1, \dots, g$, there exists a unique
$$
\heta_i:\Lnew\to\Lnew
$$ 
that satisfies the following properties:
\begin{eqnarray*}
\heta_i(t_j) & = & \delta_{ij}[t_i,h]=\delta_{ij}(t_ih-ht_i) \\
\heta_i(h) & = & 0 \\
\heta_i(ab) & = & \heta_i(a)b+a \heta_i(b) \\
\heta_i(a+b) & = & \heta_i(a) + \heta_i(b).
\end{eqnarray*}
In particular, the restriction of $\heta_i$ to $\Lloc \subset \Lnew$
equals to $\eta_i$.
\end{lemma}

\begin{proof}
First, we need to show the existence of such a $\heta_i$. A coordinate approach
is to define $\heta_i: \R \to \R$ on the ring $\R$ of Section \ref{sub.Cohn}
which is isomorphic to $\Lnew$. Given the description of $\R$, it is easy
to see that $\heta_i$ exists and is uniquely determined by the above 
properties.

An alternative, coordinate-free definition of $\heta_i$ is the following.
Let $\Lmore$ denote the Cohn localization of the
group-ring $\BQ[F \star \BZ \star \BZ]$, followed by a completion with
respect to $t-1$ and $t'-1$ (where $t$ and $t'$ are generators of $\BZ
\star \BZ$). Let $h=\log t$ and $h'=\log t'$. The map
$\varphi_{t_i\to e^{-h} t_i e^h} : \BQ[F \star \BZ \star \BZ]\to
\BQ[F \star \BZ \star \BZ]$ is $\S$-inverting, thus it extends to a map
of the localization, and further to a map $\Lmore \to \Lmore$ of the 
completion. $\Lmore$ is a $h'$-graded ring, and we can define
$$
\heta_i: \Lnew \to \Lnew
$$
by
$$
\heta_i= \varphi_{h'\to h} \circ \deg_{h'}^1 \circ \,
\varphi_{t_i\to e^{-h} t_i e^h}.
$$
It is easy to see that this definition of $\heta_i$ satisfies the properties
of the lemma, and that the properties characterize $\heta_i$.

Finally, Lemma \ref{lem.winf1} implies that the restriction of $\heta_i$ on 
$\Lloc \subset \Lnew$ equals to $\eta_i$.
\end{proof}

\begin{lemma}
\lbl{lem.etadeg1}
For all $n$ and $i$ we have that 
$$
\frac{1}{n!} \heta_i^n = \degh^n \circ \varphi_{t_i\to e^{-h} t_i e^h}.
$$
considered as maps $\La\to\Lnew$. 
\end{lemma}

\begin{proof}
In view of Lemma \ref{lem.winf1}, it suffices to show that the maps agree
on generators $t_j$ of $F_g$. For $j \neq i$, both maps vanish.
For $j=i$, the following identity
$$
e^{-h} t_i e^h = \exp ( -\mathrm{ad}_h)(t_i)=\sum_{n=0}^\infty \frac{1}{n!}
[[[t_i,h],h],\dots,h]
$$
and Lemma \ref{lem.winf1} show that both maps agree.
The result follows.
\end{proof}

Next, we extend the above lemma for matrices. 

\begin{lemma}
\lbl{lem.etadeg2}
For all $n$ and sites $i$ we have that 
$$
\frac{1}{n!} \heta_i^n = \degh^n \circ \varphi_{t_i\to e^{-h} t_i e^h}.
$$
considered as maps $\mathrm{Mat}(\La\to\BZ)\to\mathrm{Mat}(\Lnew\to\BZ)$. 
\end{lemma}

\begin{proof}
Observe the identity
\begin{equation}
\lbl{eq.eta}
\heta_i(AB)  =  \heta_i(A)B+A \heta_i(B)
\end{equation}
for matrices $A$, $B$ (not necessarily square) that can be multiplied.
It implies that
$$0 = \frac{1}{n!} \heta_i^n(I)=\frac{1}{n!} \heta_i^n(MM^{-1})=
\sum_{k=0}^n \frac{1}{k!} \heta_i^{k}(M) \heta_i^{n-k}(M^{-1}).
$$
The lemma follows by induction on $n$ and the above identity which solves
$\eta_i^n(M)$ in terms of quantities known by induction.
\end{proof}

\begin{proposition}
\lbl{prop.etadeg3}
For all $n$ and sites $i$ we have that 
$$
\frac{1}{n!} \heta_i^n = \degh^n \circ \varphi_{t_i\to e^{-h} t_i e^h}.
$$
considered as maps $\Lloc\to\Lnew$. 
\end{proposition}

\begin{proof}
Consider a matrix presentation of $s \in \Lloc$.  Ie,
$$s=(1,0,\dots,0) M^{-1} \vec{b}$$
for $M \in \mathrm{Mat}(\La\to\BZ)$ and $\vec{b}$ a column vector over $\La$.
Equation \eqref{eq.eta} implies that
$$
\frac{1}{n!} \heta_i^n(s) =
\sum_{k=0}^n \frac{1}{k!} \heta_i^{k}(M^{-1}) \heta_i^{n-k}(\vec{b}).
$$
The proposition now follows from Lemmas \ref{lem.etadeg1} and 
\ref{lem.etadeg2}.
\end{proof}

\begin{proposition}
\lbl{prop.etadeg4}
For all diagrams $D \in \A(\star_{X,\pt h},\Lloc)$, $M \in \Herm(\Lloc\to
\BZ)$, and $i=1,\dots,g$, we have that
$$
\frac{1}{n!} \heta_i^n(D)=
\degh^n\left(
\varphi_{t_i\to e^{-h} t_i e^h}(D) \sqcup \exp_\sqcup 
\left(-\frac{1}{2} \chiz(M^{-1} \varphi_{t_i\to e^{-h} t_i e^h}(M))
\right) \right) .
$$
\end{proposition}

\begin{proof}
The proof follows from a careful combinatorial counting similar to that
of the Wheels Identity of Theorem \ref{thm.wheels}.
Let us precede the proof with an informal examination of what takes
place when we repeatedly apply $\heta_i$ to the empty diagram, $D=\phi$.
The first application results in a term:
$$
\frac{1}{2} \sum_k \psdraw{circleb}{0.2in}
\left( M^{-1}(-\heta_i(M)) 
\right)_{kk}
$$
The second application requires the result of applying
$\heta_i$ to $M^{-1}$. This can be done using the following special case 
of Equation \eqref{eq.eta}:
\begin{equation}
\lbl{niceeqn}
\heta_i( M^{-1} ) = M^{-1} (-\heta_i(M)) M^{-1}.
\end{equation}
So, the second application 
results in the following terms: 
\begin{eqnarray*}
\lefteqn{
\frac{1}{2}\sum_{k,l} \psdraw{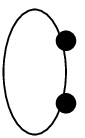}{0.2in}
\begin{array}{c}
\left( M^{-1} (-\heta_i(M)) \right)_{lk} \\
\left( M^{-1} (-\heta_i(M)) \right)_{kl}
\end{array}
+\frac{1}{2} \sum_k
\psdraw{circleb}{0.2in}
\left( M^{-1} \frac{-(\heta_i)^2}{2!}(M)) 
\right)_{kk}} \\
 & & 
 +\frac{(\frac{1}{2})^2}{2!} \sum_{k,l}
\psdraw{circleb}{0.2in}
\left( M^{-1} (-\heta_i(M)) \right)_{kk} \sqcup
\psdraw{circleb}{0.2in}
\left( M^{-1} (-\heta_i(M)) \right)_{ll} 
\end{eqnarray*}
We see that the terms arising from the application of $\heta_i$
fall into three categories.
\begin{eqnarray}
& & 
\heta_i(\phi)  =  \frac{1}{2} \sum_k
\psdraw{circleb}{0.2in}
\left( M^{-1} (-\heta_i(M)) 
\right)_{kk} \lbl{movea}\\
& & 
\psdraw{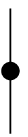}{0.05in} 
\begin{array}{c}
(\heta_i(M^{-1}))_{kl}
\end{array}
 = 
\psdraw{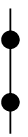}{0.05in}
\begin{array}{c}
(M^{-1})_{ml} \\
( M^{-1} (-\heta_i(M)) )_{km}
\end{array}
\lbl{moveb} \\
& & 
\psdraw{lineb}{0.05in} 
\begin{array}{c}
(\heta_i((\heta_i)^p(M)))_{kl}
\end{array}
=
\psdraw{lineb}{0.05in} 
\begin{array}{c}
((\heta_i)^{p+1}(M))_{kl}
\end{array}
\lbl{movec} 
\end{eqnarray}
Thus, we see that $\heta_i^n(\phi)$ equals to a linear combination of 
diagrams, each consisting of a number of wheel components, each of which is 
of the following form (this can be (over-)determined by a sequence
$(\kappa_1,\dots,\kappa_\mu)$ of positive integers):
\begin{equation*}
\sum_{k_1,\dots,k_\mu}
\psdraw{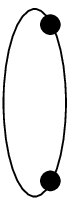}{0.2in}
\begin{array}{c}
\left( M^{-1}(-(\heta_i)^{\kappa_\mu}(M)) \right)_{k_{\mu-1} k_\mu } \\
\dots \\
\left( M^{-1}(-(\heta_i)^{\kappa_1}(M)) \right)_{k_\mu k_1} \\
\end{array}
\end{equation*}
The question, then, is to calculate an expression in terms
of such diagrams, counting ``all ways'' that they arise 
from the application of moves \eqref{movea}-\eqref{movec}. We will
generate a sequence of such moves with the following object, an
{\em edge-list} on a {\em shape}.

The {\em shape} determined by a set of positive
integers $(\mu_1,\dots,\mu_\tau)$ is a set
of oriented polygons, one with $\mu_1$ edges, and so on,
such that each has a distinguished bivalent vertex, its {\em basepoint}.
A {\em labeled shape} is a shape with its edges labeled by 
positive integers.
The diagram $D_\kappa$ {\em determined by some labeled shape} $\kappa$
is the diagram that arises by replacing each polygon by a series
of wheels, in the following manner (to distinguish them from beads, the 
bivalent vertices of the polygons will be drawn by small horizontal segments,
and the base-point by an x).
$$
\psdraw{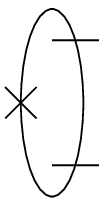}{0.2in} 
\begin{array}{c}
\kappa_\mu \\
\dots \\
\kappa_1
\end{array}
\to
\sum_{k_1,\dots,k_\mu}
\psdraw{circlebb}{0.2in}
\begin{array}{c}
\left( M^{-1}(-(\heta_i)^{\kappa_\mu}(M)) \right)_{k_{\mu-1} k_\mu }
\dots \\
\left( M^{-1}(-(\heta_i)^{\kappa_1}(M)) \right)_{k_\mu k_1}
\end{array}
$$
An {\em edge-list} on a shape is a finite sequence 
of the edges of the shape such that each edge is represented at
least once. (It is convenient to enumerate the edges along
the orientation, starting from the basepoint.)
The labeled shape associated to some edge-list on that shape is 
that shape with each edge labeled by the number of times that edge
appears in the edge-list.

An edge-list $(\beta_1,\dots,\beta_r)$
on a shape will determine a sequence of moves
(\eqref{movea}-\eqref{movec}) as follows. 

\begin{enumerate}
\item If some $\beta_i$ is the first appearance in the edge-list of an edge
from some polygon with $\mu$ edges, 
then that term determines a move \eqref{movea}, where for the purposes
of discussion the introduced
wheel is oriented with $\mu$ beads, 
with one such distinguished as the ``first'',
and with the 
label $M^{-1}(-\heta_i(M))$ sitting on bead number $\beta_i$, and with all
other beads labeled 1. 
\item If some $\beta_i$ is the first appearance in the edge-list of
some edge of a polygon, some other edge of
which has {\em already} appeared in the edge-list, 
then that term determines
a move \eqref{moveb} applied to the factor of $M^{-1}$ that immediately
follows the edge $\beta_i$ on the (oriented) wheel. (That is, the
factor of $M^{-1}(-\heta_i(M))$ that is introduced``slides'' back 
around the polygon from there to the bead number $\beta_i$.)
\item If some edge $\beta_i$ has already appeared in the path
(say, $p$ times) then that step determines the application of
a move \eqref{movec} to the factor $-(\heta_i)^p(M)$ appearing on the bead 
number $\beta_i$.
\end{enumerate}

\begin{remark}
Observe that the result of
these moves will be 
the diagram determined by the labeled shape associated to the edge-list.

Observe also that the number of different
edge-lists on a shape which give rise to a 
given associated labeled shape, e.g.
$$
\psdraw{circlex}{0.2in} 
\begin{array}{c}
\kappa_{1,\mu_1} \\
\dots \\
\kappa_{1,1}
\end{array}
\dots
\psdraw{circlex}{0.2in} 
\begin{array}{c}
\kappa_{\tau,\mu_\tau} \\
\dots \\
\kappa_{\tau,1}
\end{array}
\hspace{0.5cm}
\text{ is }
\hspace{0.5cm}
\frac{
(\kappa_{1,1}+\dots+\kappa_{1,\mu_1}+\dots+\kappa_{\tau,1}+\dots+
\kappa_{\tau,\mu_\tau})!
}
{
\left(
\kappa_{1,1}!\dots \kappa_{1,\mu_1}! 
\right)
\dots 
\left(
\kappa_{\tau,1}! \dots \kappa_{\tau,\mu_\tau}!
\right)
}.
$$
Two edge-lists determine the same sequence of moves if and only if 
they are related by symmetries of the shape:
\begin{enumerate}
\item rotating individual polygons (ie, shifting the ``base-point''),
\item permuting a subset of polygons which all have the same number
of edges.
\end{enumerate}
\end{remark}

We now define a function $\Sym$ from the set of shapes to the positive 
integers as follows.
If a shape $\sigma$ is specified by a sequence
$(\nu_1,\nu_2,\dots,\nu_k)$, where $\nu_1$ counts the number of polygons
with one edge, $\nu_2$ counts the number of polygons with $2$ edges,
and so on, then 
$$
\Sym(\sigma) = \left(\nu_1!\nu_2!\dots\nu_k!\right)
\left(1^{\nu_1}2^{\nu_2}\dots k^{\nu_k}\right).
$$
We are now ready to prove Proposition \ref{prop.etadeg4}.
First consider the case where the given diagram $D= \emptyset$ is empty.
\begin{eqnarray*}
\lefteqn{ \sum_n
\frac{
(\heta_i)^n}{n!}
(\emptyset)  } & & \\
& = &
\sum_{\sigma\ \mbox{a shape}}
\frac{1}{\Sym(\sigma)}
\sum_{\kappa:\ \mbox{a labeling of $\sigma$}}
\frac{1}{
\left(
\kappa_{1,1}!\dots \kappa_{1,\mu_1}! 
\right)
\dots 
\left(
\kappa_{\tau,1}! \dots \kappa_{\tau,\mu_\tau}!
\right)
}
D_\kappa \\
& = &
\es
\left( -\frac{1}{2} \mbox{Tr}
\left( -\sum_{m=1}^{\infty}\frac{1}{m}\left(M^{-1}
\left(-\heta_i(M) + \frac{-(\heta_i)^{2}}{2!}(M)
+ 
\dots
\right) \right)^m \right) \right) \\
& = &
\exp_\sqcup
\left( -\frac{1}{2} \chiz
\left( M^{-1} \phi_{t_i\to e^{-h} t_i e^h}(M)  \right) \right),
\end{eqnarray*}
as needed. In case $D$ is not empty, assume for convenience that there are 
$k$ beads on $D$, and explicitly record their labels
$D(q_1,\dots.q_k)$. Now,
$$
\sum_n
\frac{
(\heta_i)^n}{n!}
(D)  = 
\sum_{n_1+ \dots + n_k + n_{k+1} = n} D\left( 
\frac{(\heta_i)^{n_1}}{n_1!}(q_1),\dots,
\frac{(\heta_i)^{n_k}}{n_k!}(q_k) \right) \sqcup
\frac{(\heta_i)^{n_{k+1}}}{n_{k+1}!}(\emptyset),
$$
and the Proposition follows from the above case of $D = \emptyset$.
\end{proof}

\begin{proposition}
\lbl{prop.wequal}
The $\winf$ and $\winfp$ relations on $\A(\star_X,\Lloc)$ coincide.
\end{proposition}

\begin{proof}
Consider two pairs $(M,s_1)$ and $(M,s_2)$ related by a wrapping move 
$\winfp_i$ as in definition \ref{def.wrapping}.
In other words, for some $s \in \GA(\star_{X,\pt h},\Lloc)$ we have that
\begin{eqnarray*}
s_m & = & \conh \left( 
\varphi_{m,i}(s)\sqcup \left( 
\exp_\sqcup \left(-\frac{1}{2} \chiz(M^{-1}
\varphi_{m,i} M) \right)\right)\right)
\end{eqnarray*}
for $m=1,2$ where $\varphi_{1,i}$ is the identity and 
$\varphi_{2,i}=\phi_{t_i\rightarrow e^{-h}t_ie^h}$.
Consider the case that $s$ has a diagram with $n$ legs
labeled by $\pt h$. Proposition \ref{prop.etadeg4} implies that
$$
s_2=\frac{1}{n!} \conh (\heta_i^n(s)).
$$
Let $s'$ denote the result of gluing all but one of the $n$ $\pt h$ legs
of $\heta_i^{n-1} (s)$ and renaming the resulting one $\pt h'$.
In other words, $s'=\conh(\heta_i^{n-1}(s)(\pt h + \pt h'))$.
Let $\heta_i'$ denote the action of $\heta_i$ for $\pt h'$-labeled legs. 
If $s$ has at least one $\pt h$ leg, it follows that
$$
s_2=\con_{\{h'\}}(\heta_i'(s'))=0 \bmod \winf_i
$$
by the infinitesimal wrapping relations, and it follows that $s_1=0$
by definition. If $s$ has no $\pt h$ leg, then $s_1=s_2$. Thus, in all
cases, $(M,s_1)\stackrel{\winf_i}\sim (M,s_2)$. The converse is obvious;
thus the result follows.
\end{proof}

\begin{proof}(of Proposition \ref{prop.winf3})
It is obvious that the $\winfp$ relation implies the group-like $\wgp$
relation of Section \ref{sub.wrapping}. Since $\winfp=\winf$, the result
follows.
\end{proof}

\section{Infinitesimal basing relations}
\lbl{sec.binf}

This section, although it is not used anywhere in the paper, is added for
completeness. The reader may have noticed that in Section \ref{sec.relations},
we mentioned infinitesimal versions of the wrapping and the $\CA$ relations,
but not of the basing relations. The reason is that we did not need them;
in addition $\intrat$ is not known to preserve them. Nevertheless, we 
introduce them here.

We define
$$
\A(\ostar_X,\Lloc)=\A(\star_X,\Lloc)/\lp \binf_1,\binf_2 \rp ,
$$
where $\binf_{1,X}$ is the subspace of {$X$-flavored 
link relations} that appeared in \cite[Part II, Section 5.2]{A}:
$$
\setlength{\unitlength}{0.03\standardunitlength}
	\begin{array}{c}  \hspace{-1.7mm}
         	\raisebox{-8pt}{\input draws/LinkRel.epc }
         	\hspace{-1.9mm}
	\end{array}
 
$$
and $\binf_{2,X}$ denote subspace generated by an action
$$
F \times \A(\star_X,\Lloc) \to \A(\star_X,\Lloc)
$$
described by pushing some element of
$f\in F_g$ onto all $X$-labeled legs of a diagram. 

We now compare the infinitesimal and group-like basing relations as follows.

\begin{proposition}
\lbl{prop.cbasing}
{\rm (a)} We have an isomorphism
$$
\A(\dots,\ostar_X,\Lloc)=
\A(\dots,\star_X,\Lloc)/\la \binf_X \ra  \stackrel{\simeq}\to
\A(\dots,\tcircle_X,\Lloc).
$$
{\rm (b)} The inclusion
$\GAz(\dots,\star_X,\Lloc) \longrightarrow \A^0(\dots,\star_X,\Lloc)$
maps the $\bgp$-relations to $\binf$-relations, and induces a map
$$
\GAz(\dots,\ostar_X,\Lloc) \longrightarrow \A^0(\dots,\star_X,\Lloc).
$$
\end{proposition}

\begin{proof}
The first part follows from the proof of \cite[Part II, Theorem 3]{A} 
without essential changes. The second part is standard.
\end{proof}

\section{A useful calculation}
\lbl{sub.useful}

In this section, we quote a useful and often appearing computation
that uses the identities of Lemma \ref{lem.identities}.

Let $\varphi: \GA(\star_X,\Lloc)\to \GA(\star_{X \cup \{h\}},\Lloc)$ be a map
that commutes with respect to gluing $X$-colored legs. 

\begin{theorem}
\lbl{thm.varphi}
For $X' \subset X$ and $s \in \IGA$, decomposed as in Equation 
\eqref{eq.Gaussd}, we have that
$$
\intrat dX' (\varphi(s))=
\intrat dX' \left(
\exp\left(
\frac{1}{2}
\sum
\strutb{x_j}{x_i}{\varphi M_{ij}}
\right)
\right) 
\varphi \left(\intrat dX' (s) \right) .
$$
If in addition $\varphi(s)-s$ is $X$-substantial for all 
$s \in \GA(\star_X,\Lloc)$, then
$$
\intrat dX' (\varphi(s))=
\exp\left(-\frac{1}{2}\chiz(M^{-1} \varphi M)\right)
\varphi \left(\intrat dX' (s) \right) .
$$
\end{theorem}

\begin{proof}
Consider the decomposition of $s$ given by Equation \eqref{eq.Gaussd}.
We begin with a ``$\delta$-function trick'' of Equation \eqref{identf} of
\ref{lem.identities}:
$$
s= \es \left( \frac{1}{2} \sum_{i,j}
 \strutb{x_j}{x_i}{M_{ij}} \right) \sqcup R
= \left\la R(\pt y), \, \exp\left( \frac{1}{2} \sum 
\strutb{x_j}{x_i}{M_{ij}}
+ 
\strutb{y_j}{x_i}{}   
\right) \right\ra_Y.
$$
Thus,
\begin{eqnarray*}
\intrat dX' \left(\varphi (s) \right)
& = &
\left\la
\varphi\left(
R(\pt y)\right),
\intrat dX' \left( \exp\left( \frac{1}{2} \sum
\strutb{x_j}{x_i}{\varphi M_{ij}}
+ 
\strutb{y_i}{x_i}{}
\right) \right)
\right\ra_Y .
\end{eqnarray*}
\noindent
Since $\div_{X'} \exp\left(
-\sum
\strutb{\pt x_j}{y_i}{\varphi M^{-1}_{ij}}
\right)=1$, 
the integration by parts Lemma \ref{intbypartslem} implies that

\begin{multline*}
\intrat dX' \left( 
\exp\left( \frac{1}{2} \sum
\strutb{x_j}{x_i}{\varphi M_{ij}}
+ 
\strutb{y_j}{x_i}{}
\right) \right)  = \\ 
\intrat dX' \left( 
\left(
\exp\left(
-\sum
\strutb{\pt x_j}{y_i}{\varphi M^{-1}_{ij}}
\right) \right)\flat_{X'}
\exp\left( \frac{1}{2} \sum
\strutb{x_j}{x_i}{\varphi M_{ij}}   
+ 
\strutb{y_i}{x_i}{}
\right) \right) .
\end{multline*}

\noindent
Completing the square, using
Equation \eqref{identg} of Lemma \ref{lem.identities}, implies that the above 
equals to:
$$
\exp\left( -\frac{1}{2} \sum  
\strutb{y_j}{y_i}{\varphi M_{ij}^{-1}}
 \right)
\,\,\,\,\,
\intrat dX' \left(
\exp\left(
\frac{1}{2}
\sum
\strutb{x_j}{x_i}{\varphi M_{ij}}
\right)
\right).
$$
\noindent
Returning to the expression in question:
\begin{eqnarray*}
\intrat dX' \left(\varphi (s) \right) 
& = & \intrat dX' \left(
\exp\left(
\frac{1}{2}
\sum
\strutb{x_j}{x_i}{\varphi M_{ij}}
\right)
\right) 
  \\  & & \quad \sqcup 
\left\la
\varphi\left(
R(\pt y)\right),
\exp\left( -\frac{1}{2} \sum  
\strutb{y_j}{y_i}{\varphi M_{ij}^{-1}}
 \right) \right\ra_Y. 
\end{eqnarray*}
\noindent
The second factor equals to $\varphi \left(\intrat dX' (s) \right)$,
and concludes the first part of the theorem.

In order to compute the first factor, we need to separate its $X'$-substantial
part. With the added assumption on $\varphi$, we can compute the first
first factor in a manner analogous to the Wheels Identity of Theorem 
\ref{thm.wheels}
\begin{eqnarray*}
\lefteqn{ \intrat dX' \left(
\exp\left(
\frac{1}{2}
\sum
\strutb{x_j}{x_i}{\varphi M_{ij}}
\right)
\right)  } & & \\
& = &
\intrat dX' \left( 
\exp \left(
\frac{1}{2}
\sum \left( 
\strutb{x_j}{x_i}{M_{ij}}
+
\strutb{x_j}{x_i}{\varphi M_{ij}}
-
\strutb{x_j}{x_i}{M_{ij}}
\right) \right) \right) 
\\ & = &
\left\la
\exp \left(
-\frac{1}{2}
\sum
\strutb{\pt x_j}{\pt x_i}{M_{ij}^{-1}}
\right),
\exp \left( \frac{1}{2}
\strutb{x_j}{x_i}{\varphi M_{ij}}
-
\strutb{x_j}{x_i}{M_{ij}}
\right)
\right\ra_{X'} \\
& = &
\exp \left( -\frac{1}{2} \chiz(M^{-1} \varphi M) \right) .
\end{eqnarray*}
\end{proof}

\begin{remark}
\lbl{rem.apply}
The proof of Theorem \ref{thm.varphi} can be applied without changes to
the case of the map $\varphi=\hair_H: 
\A(\star_X,\Lloc)\to \A(\star_{X \cup H})$ in order to show that for
$s \in \IPR_{X'}\GA(\star_X,\Lloc)$, we have that
$$
\int dX' (\hair_H(s))=\exp\left(-\frac{1}{2} \chiz(M^{-1}\hair_H M)\right)
\hair_H \left(\intrat dX' (s) \right) .
$$
\end{remark}

\section{The Aarhus integral preserves group-like basing relations}
\lbl{sub.omissionA}

In this appendix we fill a historical gap on the Aarhus integral
and how it deals with group-like basing relations.

In \cite[Part II, Prop.5.6]{A} it was shown that the Aarhus integral
is invariant under infinitesimal basing relations.
Later on, in \cite[Prop.2.2]{BL}, it was shown that the Aarhus integral
is invariant under group-like basing relations.

Using the results of our paper, one can prove that the Aarhus integral
is invariant under group-like basing relations.
That is,

\begin{theorem}{\rm\cite{A}}\qua
\lbl{thm.cyclicAarhus}
The Aarhus integral $\int$ descends to a map:
$$
\intrat dX' : 
\IPR_{X'}\GA(\star_X)/\la \bgp_{X} \ra
\rightarrow
\GA(\star_{X-X'})/\la \bgp_{X-X'} \ra
$$
where
$
\IPR_{X'}\GA(\star_X)
$
denotes the set of integrable with respect to $X'$ group-like elements
and $\bgp_{X'}$ denotes the group-like basing relation discussed in
Section \ref{sub.groupA}.
\end{theorem}

\begin{proof}
Use the form of the group-like basing relations given in
Definition \ref{def.GAalternative} (see Lemma \ref{lem.b1G}),
and follow the proof of case $\bgp_2$ of Theorem \ref{thm.cyclic},
using Theorem \ref{thm.varphi}. The result follows.
\end{proof}

\end{document}